\documentclass[a4paper]{article}

\makeatletter
\newif\iflabel\labelfalse \let\w@label=\label
\def\label{\global\labeltrue\w@label}
\let\w@eqn=\equation
\def\equation{\global\labelfalse\w@eqn}
\def\endequation{\iflabel\@eeq\else\@eeqw\fi$$\global\@ignoretrue}
\def\@eeq{\eqno \@eqnnum}
\def\@eeqw{\addtocounter{equation}{-1}}
\newskip\bw@\newskip\bws@
\def\dc@@pq{\global \bw@=\belowdisplayskip \global\bws@=\belowdisplayshortskip
   \global\belowdisplayshortskip=3pt plus -3pt
   \global\belowdisplayskip=\belowdisplayshortskip
   \hskip 500pt minus 500pt\relax}
\def\dc@pq{\dc@@pq$$}
\def\fc@pq{$$\hskip 500pt minus 500pt\global
   \belowdisplayshortskip=\bws@\global\belowdisplayskip=\bw@}
\def\coupeq{\dc@pq\fc@pq}
\def\coupepas{\par\noindent\begin{minipage}{\textwidth}}
\def\jusquela{\end{minipage}}
\makeatother

\usepackage{amsfonts,amssymb,amsthm}
\usepackage[totalwidth=17cm,totalheight=24cm]{geometry}

\newcommand{\dr}{\partial}
\newcommand{\ots}{\Omega^{0,1}_{\II}(S)}

\newcommand{\C}{{\mathbb C}}
\newcommand{\R}{{\mathbb R}}
\newcommand{\N}{{\mathbb N}}
\newcommand{\Z}{{\mathbb Z}}
\newcommand{\II}{I\hspace{-0.1cm}I}
\newcommand{\III}{I\hspace{-0.1cm}I\hspace{-0.1cm}I}

\newcommand{\Ib}{\overline{I}}
\newcommand{\IIb}{\overline{\II}}
\newcommand{\hessb}{\overline{\mbox{Hess}}}
\newcommand{\nablab}{\overline{\nabla}}
\newcommand{\tr}{\mbox{tr}}

\newcommand{\dev}{\mbox{dev}}

\newcommand{\SL}{{\mbox{sl}}}

\newcommand{\SO}{\mbox{SO}}

\newcommand{\rp}{\mbox{{\bf R}P}}
\newcommand{\can}{\mbox{can}}

\newcommand{\hess}{\mbox{Hess}}
\newcommand{\dv}{\mbox{div}}

\newcommand{\gammab}{\overline{\gamma}}
\newcommand{\kappab}{\overline{\kappa}}
\newcommand{\sigmab}{\overline{\sigma}}

\newcommand{\Omegab}{\overline{\Omega}}

\newcommand{\fb}{\overline{f}}
\newcommand{\gb}{\overline{g}}
\newcommand{\hb}{\overline{h}}

\newcommand{\ub}{\overline{u}}
\newcommand{\vb}{\overline{v}}
\newcommand{\wb}{\overline{w}}
\newcommand{\xb}{\overline{x}}
\newcommand{\yb}{\overline{y}}
\newcommand{\zb}{\overline{z}}

\newcommand{\Db}{\overline{D}}
\newcommand{\Eb}{\overline{E}}

\newcommand{\Jb}{\overline{J}}

\newcommand{\Mb}{\overline{M}}
\newcommand{\Nb}{\overline{N}}
\newcommand{\Pb}{\overline{P}}
\newcommand{\Rb}{\overline{R}}
\newcommand{\Sb}{\overline{S}}
\newcommand{\cEb}{\overline{\mathcal E}}
\newcommand{\cFb}{\overline{\mathcal F}}
\newcommand{\cHb}{\overline{\mathcal H}}

\newcommand{\omegab}{\overline{\omega}}
\newcommand{\taub}{\overline{\tau}}

\newcommand{\et}{\tilde{e}}

\newcommand{\Kt}{\tilde{K}}
\newcommand{\Mt}{\tilde{M}}

\newcommand{\Sigmat}{\tilde{\Sigma}}
\newcommand{\nablat}{\tilde{\nabla}}

\newtheorem{prop}{Proposition}[section]
\newtheorem{lemma}[prop]{Lemma}
\newtheorem{sublemma}[prop]{Sub-lemma}

\newtheorem{thm}[prop]{Theorem}
\newtheorem{cor}[prop]{Corollary}
\newtheorem{remark}[prop]{Remark}

\newtheorem{df}[prop]{Definition}
\newtheorem{pty}[prop]{Property}
\newtheorem{question}[prop]{Question}

\newcommand{\pg}{\paragraph}

\newenvironment{tmn}[1]{\vskip 0.2cm \noindent{\bf Theorem #1.} \it}{\rm
\vspace{0.2cm}} 
\newenvironment{crn}[1]{\vskip 0.2cm \noindent{\bf Corollary #1.} \it}{\rm
\vspace{0.2cm}} 
\newenvironment{lmn}[1]{\vskip 0.2cm \noindent{\bf Lemma #1.} \it}{\rm
\vspace{0.2cm}} 
\newenvironment{qn}[1]{\vskip 0.2cm \noindent{\bf Question #1.} \it}{\rm
\vspace{0.2cm}} 

\newcommand{\btm}{\begin{thm}}
\newcommand{\etm}{\end{thm}}
\newcommand{\bpt}{\begin{pty}}
\newcommand{\ept}{\end{pty}}
\newcommand{\blm}{\begin{lemma}}
\newcommand{\elm}{\end{lemma}}
\newcommand{\bsl}{\begin{sublemma}}
\newcommand{\esl}{\end{sublemma}}
\newcommand{\bcr}{\begin{cor}}
\newcommand{\ecr}{\end{cor}}
\newcommand{\bdf}{\begin{df}}
\newcommand{\edf}{\end{df}}
\newcommand{\bprop}{\begin{prop}}
\newcommand{\eprop}{\end{prop}}
\newcommand{\bas}{\begin{asser}}
\newcommand{\eas}{\end{asser}}
\newcommand{\beq}{\begin{equation}}
\newcommand{\eeq}{\end{equation}}
\newcommand{\bpv}{\begin{proof}}
\newcommand{\epv}{\end{proof}}
\newcommand{\bpvs}{\begin{sketch}}
\newcommand{\epvs}{\end{sketch}}
\newcommand{\bit}{\begin{itemize}}
\newcommand{\eit}{\end{itemize}}
\newcommand{\bpn}{\begin{pfn}}
\newcommand{\epn}{\end{pfn}}
\newcommand{\btn}{\begin{thn}}
\newcommand{\etn}{\end{thn}}
\newcommand{\bcn}{\begin{crn}}
\newcommand{\ecn}{\end{crn}}
\newcommand{\bqn}{\begin{qn}}
\newcommand{\eqn}{\end{qn}}
\newcommand{\bln}{\begin{lmn}}
\newcommand{\eln}{\end{lmn}}
\newcommand{\brk}{\begin{remark}}
\newcommand{\erk}{\end{remark}}
\newcommand{\bq}{\begin{question}}
\newcommand{\eq}{\end{question}}

\newenvironment{pfn}[1]{\vskip 0.2cm \noindent{\it Proof #1.}}{$\square$
\vspace{0.2cm}}

\newcommand{\cE}{\mathcal{E}}
\newcommand{\cF}{\mathcal{F}}
\newcommand{\cG}{\mathcal{G}}
\newcommand{\cH}{{\mathcal H}}

\newcommand{\cR}{{\mathcal R}}

\newcommand{\cL}{\mathcal{L}}
\newcommand{\isom}{\mathrm{Isom}}
\newcommand{\db}{\overline{\partial}}

\newcommand{\fd}{\dot{f}}
\newcommand{\gd}{\dot{g}}

\newcommand{\Ad}{\dot{A}}
\newcommand{\Bd}{\dot{B}}
\newcommand{\Nd}{\dot{N}}

\newcommand{\phid}{\dot{\phi}}

\newcommand{\IId}{\dot{\II}}

\begin{document}

\title{Hyperbolic manifolds with convex boundary}

\author{Jean-Marc Schlenker\thanks{
Laboratoire Emile Picard, UMR CNRS 5580,
UFR MIG, Universit{\'e} Paul Sabatier,
118 route de Narbonne,
31062 Toulouse Cedex 4,
France.
\texttt{schlenker@picard.ups-tlse.fr; http://picard.ups-tlse.fr/\~{
}schlenker}. }}

\date{May 2002; revised, Aug. 2004 (v5)}

\maketitle

\begin{abstract}

Let $(M, \dr M)$ be a 3-manifold, which carries a hyperbolic metric with
convex boundary. We consider the hyperbolic
metrics on $M$ such that the boundary is smooth and strictly convex. We
show that the induced metrics on the boundary are exactly the metrics
with curvature $K>-1$, and that the third fundamental forms of $\dr M$
are exactly the metrics with curvature $K<1$, for which the
closed geodesics which are contractible in $M$ have length
$L>2\pi$. Each is obtained exactly once.  

Other related results describe existence and uniqueness properties for
other boundary conditions, when the metric which is achieved on $\dr M$
is a linear combination of the first, second and third fundamental
forms.

\bigskip

\begin{center} {\bf R{\'e}sum{\'e}} \end{center}

Soit $(M, \dr M)$ une vari{\'e}t{\'e} de dimension 3, qui admet une
m{\'e}trique hyperbolique {\`a} bord convexe. On
consid{\`e}re l'ensemble des m{\'e}triques hyperboliques sur $M$ pour lesquelles
le bord est lisse et strictement convexe. On montre que les m{\'e}triques
induites sur le bord sont exactement les m{\'e}triques {\`a} courbure $K>-1$, et
que les troisi{\`e}mes formes fondamentales du bord sont exactement les
m{\'e}triques {\`a} courbure $K<1$ dont les g{\'e}od{\'e}siques ferm{\'e}es qui sont
contractiles dans $M$ 
sont de longueur $L>2\pi$. Chacune est obtenue pour une unique m{\'e}trique
hyperbolique sur $M$. 

D'autres r{\'e}sultats, reli{\'e}s {\`a} ceux-ci, d{\'e}crivent des propri{\'e}t{\'e}s
d'existence et d'unicit{\'e} pour d'autres conditions au bord, lorsque la
m{\'e}trique qu'on prescrit sur $\dr M$ est une combinaison lin{\'e}aire des
premi{\`e}re, seconde et troisi{\`e}me formes formes fondamentales. 

\end{abstract}


\tableofcontents

\pg{Convex surfaces in $H^3$}

Consider a compact, convex subset $\Omega$ of $H^3$. Its boundary $\dr
\Omega$ is a convex sphere. Suppose that $\dr \Omega$ is smooth and
strictly convex (i.e. its principal curvatures never vanish); its
induced metric, which is also called its "first fundamental form" and
denoted by $I$, then has curvature $K>-1$.

A well-known theorem, due mainly to Aleksandrov \cite{Al} and Pogorelov
\cite{Po}, states that, conversely, each smooth metric with curvature
$K>-1$ on $S^2$ is induced on a unique smooth, strictly convex surface
in $H^3$. Analogous results also hold in the Euclidean 3-space \cite{N} and
in the 3-sphere. 

We consider here this question, in the more general setting where the
sphere (which can be considered as the boundary of a ball) is replaced
by the boundary of a rather arbitrary 3-manifold $M$.  That is, we will show
that the hyperbolic metrics 
on $M$ with smooth, strictly convex boundary are uniquely determined by
the induced metrics on the boundary, and that, under
mild topological assumptions on $M$, the possible induced
metrics on the boundary are exactly the metrics with curvature
$K>-1$. This statement was conjectured by Thurston, and the existence
part was already proved by Labourie \cite{L4}.

\pg{The third fundamental form}

There is another metric, the third fundamental form, which is defined on
a smooth, strictly convex surface $S$ in $H^3$. To define it, let $N$ be
a unit normal vector field to $S$, and let $\nabla$ be the Levi-Civit{\`a}
connection of $H^3$; the {\bf shape operator} $B:TS\rightarrow TS$ of
$S$ is defined by: $Bx:=-\nabla_xN$. It satisfies the Gauss equation,
$\det(B)=K+1$, and the Codazzi equation, $d^{\nablab}B=0$, where $\nablab$ is the
Levi-Civit{\`a} connection of the induced metric $I$, $K$ is its
curvature, and $d^{\nablab}B(x,y):=(\nablab_xB)y-(\nablab_yB)x$.

The {\bf second fundamental form} of $S$ is defined by:
$$ \forall s\in S, ~\forall x,y\in T_sS, ~\II(x,y)=I(Bx,
y)=I(x,By)~, $$ 
and its {\bf third fundamental form} by:
$$ \forall s\in S, ~\forall x,y\in T_sS, ~\III(x,y)=I(Bx,By)~. $$

The third fundamental form is ``dual'' to the induced metric on
surfaces, in a precise way that is explained in section 1. In
particular, some questions about the third fundamental form of surfaces
in $H^3$ can be translated into questions about isometric embeddings in
a dual space, the de Sitter space, whose definition is also given in
section 1. 

For smooth, strictly convex spheres in $H^3$, the possible
third fundamental forms are the metrics with curvature $K<1$ and closed
geodesics of lengths $L>2\pi$, and each is obtained in exactly one way
(see \cite{cras,these}).  
We will consider here an extension of this result to the setting of
hyperbolic manifolds with convex boundary, and obtain a description of
the possible third fundamental forms of hyperbolic manifolds with
smooth, strictly convex boundary.

\pg{Main results}

In all the paper, we consider a compact, connected 3-manifold $M$ with
boundary. We will suppose that $M$ admits a complete, hyperbolic, convex
co-compact 
metric. This is a topological assumption, and could easily be
stated in purely topological terms (see
e.g. \cite{thurston-notes}). It clearly implies that $M$ also admits a
hyperbolic metric for which the boundary is a smooth, strictly convex
surface, this can be seen by smoothing the set of points at distance $1$ from
the convex core in a complete, convex co-compact hyperbolic metric. 
In addition, we will suppose in most of the paper
--- until the second part of section \ref{se:9} --- that $M$ is
not a solid torus; this is a special case that is treated
separately at the end of section \ref{se:9}. 
Thus each boundary component of $M$ is a compact surface of
genus at least $2$, except in the simpler case where $M$ is a ball. The
two main results are then as follows. 

\btm \label{tm:I}
Let $g$ be a hyperbolic metric on $M$ such that $\dr M$ is smooth and
strictly convex. Then the induced metric $I$ on $\dr M$ has curvature
$K>-1$. Each smooth metric on $\dr M$ with $K>-1$ is induced on $\dr M$
for a unique choice of $g$. 
\etm

\btm \label{tm:III}
Let $g$ be a hyperbolic metric on $M$ such that $\dr M$ is smooth and
strictly convex. Then the third fundamental form $\III$ of $\dr M$ has
curvature $K<1$, and its closed geodesics which are contractible in $M$
have length $L>2\pi$. Each such metric is obtained for a unique choice
of $g$. 
\etm

\pg{Known cases}

Theorems \ref{tm:I} and \ref{tm:III} provide a new light on some
previous results. The existence part of Theorem
\ref{tm:I} is a result of Labourie \cite{L4}. When $M$ is a
ball, Theorem \ref{tm:I} is a consequence of the work, already
mentioned, of Aleksandrov and Pogorelov, while Theorem \ref{tm:III}
already stands in \cite{these}. Moreover, Theorem \ref{tm:III} was
also known in the special case where $M$ is ``fuchsian'' (i.e. it has an
isometric involution fixing a closed surface), see
\cite{iie}. 

When one considers polyhedra instead of smooth surfaces, the analog of
Theorem \ref{tm:I} when $M$ is a ball was obtained by Aleksandrov
\cite{alex}, while the 
analog of Theorem \ref{tm:III} was obtained by Rivin and Hodgson
\cite{Ri,RH}, following work of Andreev \cite{Andreev}.
Going to the category of ideal polyhedra, we find results of Andreev
\cite{Andreev-ideal} and Rivin \cite{rivin-annals,rivin-ideal}. 
Moreover, in this setting, results on the fuchsian case of the analog of
Theorem \ref{tm:III} are hidden in
works on circle packings (see \cite{thurston-notes}, chapter 13 or
\cite{CdeV}) or on similarity 
structures (see \cite{Ri2}). The analog of Theorem \ref{tm:III} when the
topology of $M$ is  
``general'', but when $\dr M$ is locally like an ideal polyhedron, is in
\cite{ideal}. Finally analogs of Theorems \ref{tm:I} 
and \ref{tm:III} also hold for hyperideal polyhedra when $M$ is a ball,
see \cite{shu} and the result of Bao and Bonahon \cite{bao-bonahon}, and
when $M$ is "fuchsian", a result of Rousset \cite{rousset1}. 
Here again the topologically general
case of the analog of Theorem \ref{tm:III} holds, see \cite{hphm}. This
is also related to results on complete smooth surfaces, see \cite{rsc}.
Thus several important results on compact or ideal polyhedra can be seen
as polyhedral versions of Theorems \ref{tm:I} and \ref{tm:III}, in the
special case where $M$ has no topology. 

Another setting where results similar to Theorems \ref{tm:I} and
\ref{tm:III} should hold is for the convex cores of hyperbolic convex
co-compact manifolds. Their induced metrics are hyperbolic, and the fact
that each hyperbolic metric can be obtained is a consequence either of
Labourie's result \cite{L4} or of the
known (weak) form of a conjecture of Sullivan, see \cite{epstein-marden}
(or \cite{bridgeman-canary} when the boundary of $M$ is compressible). The
uniqueness, however, remains elusive. The questions concerning the third
fundamental form can be stated in terms of bending lamination, and
existence results here were obtained recently by Bonahon and Otal, and by
Lecuire \cite{bonahon-otal,lecuire}; the uniqueness is known only when
the bending lamination is supported by 
closed curves, in this case it is proved in \cite{bonahon-otal}.

\pg{Other conditions}

As a consequence of Theorems \ref{tm:I} and \ref{tm:III}, we can fairly
easily find corollaries which are also generalizations. They describe
other metrics, defined from the extrinsic invariants of the boundary,
that one can prescribe, like $I$ and $\III$ respectively in Theorems
\ref{tm:I} and \ref{tm:III}. Since the second fundamental form appears
here, we have to specify that $\II$ is defined with respect to the
interior unit normal of $\dr M$. The geometric meaning of those
statements is not too explicit, and they are of interest mainly if one
considers the problems treated here as dealing with non-linear analysis
problems on manifolds with boundary, namely, finding a solution $g$ of an
elliptic PDE on $M$, $\mbox{ric}_g=-2g$, with some boundary conditions. 
Here again, $M$ is a compact 3-manifold
with boundary which admits a complete, convex co-compact hyperbolic metric.

\btm \label{tm:cor-I}
Let $k_0\in (0,1)$. Let $h$ be a Riemannian metric on $\dr M$. There
exists a hyperbolic metric $g$ on $M$, such that the principal
curvatures of the boundary are between $k_0$ and $1/k_0$, and for which:
$$ I - 2k_0 \II + k_0^2 \III = h $$
if and only if $h$ has curvature $K\geq -1/(1-k_0^2)$. $g$ is then
unique. 
\etm

\btm \label{tm:cor-III}
Let $k_0\in (0,1)$. Let $h$ be a Riemannian metric on $\dr M$. There
exists a hyperbolic metric $g$ on $M$, such that the principal
curvatures of the boundary are between $k_0$ and $1/k_0$, and for which:
$$ k_0^2 I - 2k_0 \II + \III = h $$
if and only if $h$ has curvature $K\leq 1/(1-k_0^2)$, and its closed
geodesics which are contractible in $M$ have length
$L>2\pi\sqrt{1-k_0^2}$. $g$ is then unique. 
\etm

When one takes the limit $k_0\rightarrow 0$ in those statements, one
finds Theorems \ref{tm:I} and \ref{tm:III}. The proofs, which are done in
section 9, are based on a simple trick using normal deformations of the
boundary. 

In addition, one can also prescribe another ``metric'' on the boundary,
the ``horospherical metric'' $I+2\II+\III$, see \cite{horo}. This is done by
very different methods, and has a different flavor. 

Putting all the results above together and adding an elementary scaling
argument, we could state a general result describing existence and
uniqueness properties for boundary conditions where one prescribes a metric
on $\dr M$ of the form:
$$ aI - 2\sqrt{ab}\II + b\III, ~ \mbox{with} ~ a,b\geq 0 ~ \mbox{and}
~ (a,b)\neq (0,0)~. $$ 

On the other hand, one should not be too optimistic about what sort of
boundary conditions one can prescribe. For instance, just by considering
the spheres of radius $r\in (0,\infty)$ one can show that many metrics
made from 
linear combination of the extrinsic invariants of the boundary do not
lead to any infinitesimal rigidity result.

\pg{Outline of the proofs}

We now give a brief outline of the proofs of Theorems \ref{tm:I} and
\ref{tm:III}. Both proofs --- and they are closely related --- use a
deformation argument, and show that a natural map between two spaces is
a homeomorphism.

\bdf \label{df:metriques}
We call:
\begin{itemize}
\item $\cG$ the space of hyperbolic metrics on $M$, with strictly convex,
smooth boundary, considered up to isotopy.
\item $\cH$ the space of $C^\infty$ metrics on $\dr M$ with curvature
  $K>-1$, up to isotopy.
\item $\cH^*$ the space of $C^\infty$ metrics on $\dr M$ with curvature
  $K<1$, with closed geodesics of length $L>2\pi$ when they are
  contractible in $M$, up to
  isotopy. 
\item $\cF:\cG\rightarrow \cH$ the map sending a hyperbolic metric on
  $M$ to the induced metric on the boundary.
\item $\cF^*:\cG\rightarrow \cH^*$ the map sending a hyperbolic metric
  on $M$ to the third fundamental form of the boundary.
\end{itemize}
\edf

The main point will be that $\cF$ and
$\cF^*$ are homeomorphisms. This will follow from the following facts:
\begin{itemize}
\item they are both Fredholm operators of index $0$;
\item they are both locally injective, i.e. their differential is
injective at each point of $\cG$;
\item they are proper;
\item $\cG$ is connected, while both $\cH$ and $\cH^*$ are simply
connected. 
\end{itemize}

The main technical point of the paper is the local injectivity
statement, which can naturally be formulated as infinitesimal rigidity:
given a hyperbolic metric $g$ on $M$ with smooth, strictly convex
boundary, one can not deform it infinitesimally without changing the
induced metric on the boundary (resp. the third fundamental form of the
boundary). There are several possible proofs of this for convex surfaces in
$H^3$, for instance based on an index argument \cite{L1}, and for convex
polyhedra in $H^3$ it can be proved using an argument already used by Legendre
\cite{legendre} and Cauchy \cite{cauchy} in $\R^3$, see \cite{RH}. In our
context, however, things are different since the ``relative position'' of the
boundary components has to be taken into account.

\begin{lmn}{\ref{lm:rigidity}}
Let $g\in \cG$. For any non-trivial first-order deformation $\gd\in
T_g\cG$, the induced first-order variation of the induced metric on $\dr
M$ is non-zero. 
\end{lmn}

\begin{lmn}{\ref{lm:rigidity*}}
Let $g\in \cG$. For any non-trivial first-order deformation $\gd\in
T_g\cG$, the induced first-order variation of the third fundamental form
of $\dr M$ is non-zero. 
\end{lmn}

Those lemmas are the subject of sections 1 through 6. Section 1 starts
with some basic well-known results on hyperbolic 3-space, and also
contains some less standard constructions, like the hyperbolic-de Sitter
duality and the Pogorelov Lemma on isometric deformations of
surfaces. Section 2 
gives details about the geometry of hyperbolic manifolds with convex
boundary and their universal covers, and explains how the infinitesimal
rigidity question can be translated as a Euclidean geometry statement. 
Section 3 is concerned
with some aspects of the rigidity of smooth, strictly convex surfaces in
$\R^3$ (or in $H^3$), and section 4 contains some estimates on geometric
objects associated to a first-order deformation of the metric. 
All those elements are then used in 
sections 5 and 6 to prove Lemmas \ref{lm:rigidity} and
\ref{lm:rigidity*}. 

\medskip

We also have to prove that $\cF$ and $\cF^*$ are proper; this can be
formulated as a compactness property, more precisely as the following
lemmas. 

\begin{lmn}{\ref{lm:compact}}
Let $(g_n)_{n\in \N}$ be a sequence of hyperbolic metrics on $M$ with
smooth, strictly convex boundary. Let $(h_n)_{n\in \N}$ be the sequence of
induced metrics on the boundary. Suppose that $(h_n)$ converges to a
smooth metric $h$ on $\dr M$, with curvature $K>-1$. Then there is a
sub-sequence of $(g_n)$ which converges to a hyperbolic metric $g$ with
smooth, convex boundary. 
\end{lmn}

\begin{lmn}{\ref{lm:compact*}}
Let $(g_n)_{n\in \N}$ be a sequence of hyperbolic metrics on $M$ with
smooth, strictly convex boundary. Let $(h_n)_{n\in \N}$ be the sequence of
third fundamental forms of the boundary. Suppose that $(h_n)$ converges to a
smooth metric $h$ on $\dr M$, with curvature $K<1$ and closed geodesics
of length $L>2\pi$ when they are contractible in $M$. Then there is a
sub-sequence of $(g_n)$ which converges to a hyperbolic metric $g$ with
smooth, convex boundary. 
\end{lmn}

They are proved in section 7, using an approach previously developed in
\cite{L1,these} for isometric embeddings of convex surfaces in
Riemannian (resp. Lorentzian) 3-manifolds. Section 8 follows with the
proof that $\cG$ is connected, and that $\cH$ and $\cH^*$ are simply
connected.  

\medskip

The third ingredient is the fact that $\cF$ and $\cF^*$ are Fredholm:

\begin{lmn}{\ref{lm:indice}}
$\cF$ and $\cF^*$ are Fredholm maps of index $0$.
\end{lmn}

This is related to the elliptic nature of the PDE problem underlying
this situation. The proof is in section 9. The end of the proofs of
Theorems \ref{tm:I} and \ref{tm:III} stand also in section 9, and section 10
contains the proofs of Theorems \ref{tm:cor-I} and \ref{tm:cor-III}, as
well as some remarks and open questions.


\section{Background on hyperbolic geometry} \label{se:1}

This section contains various constructions and results on hyperbolic
3-space. Some are completely classical, and given here only for
reference elsewhere in the text, while other are less 
well-known.

\pg{Hyperbolic 3-space and the de Sitter space}

Hyperbolic 3-space can be constructed as a quadric in the Minkowski
4-space $\R^4_1$, with the induced metric:
$$ H^3 := \{ x\in \R^4_1 ~|~ \langle x,x \rangle =-1 ~ \wedge x_0>0
\}~. $$
But $\R^4_1$ also contains another quadric, which is the called the de
Sitter space of dimension 3:
$$ S^3_1  := \{ x\in \R^4_1 ~|~ \langle x,x \rangle =1\}~. $$
By construction, it is Lorentzian and has an action of $\SO(3,1)$ which
is transitive on orthonormal frames, so it has constant curvature; one
can easily check that its curvature is $1$. $S^3_1$ contains many
space-like totally geodesic 2-planes, each isometric to $S^2$ with its
canonical round metric. Each separates $S^3_1$ into two ``hemispheres'',
each isometric to a "model" which we call $S^3_{1,+}$, and which we will take
to be, in the quadric above, 
the set of points $x\in S^3_1$ with positive first coordinate. 

\pg{Projective models}

Using the definitions of $H^3$ and $S^3_1$ as quadrics, given above, we
can define the projective models of $H^3$ and of the de Sitter
hemisphere $S^3_{1,+}$. We will call $D^3$ the ball of radius $1$ and
center $0$ in $\R^3$.
We will be using the classical projective model of hyperbolic space,
which is well known, and the almost as classical projective model of the
de Sitter space.

\bdf \label{df:klein}
Let:
$$ \begin{array}{rrcl}
\phi_H: & H^3 & \rightarrow & D^3 \\
& x & \mapsto & x/x_0~, \\
&&& \\
\phi_S: & S^3_{1,+} & \rightarrow & \R^3\setminus \overline{D^3} \\
& x & \mapsto & x/x_0~, \\
&&& \\
\phi: & H^3\cup S^3_{1,+} & \rightarrow & \R^3 \\
& x & \mapsto & \left\{ 
\begin{array}{cc}
\phi_H(x) & \mbox{if $x\in H^3$,} \\
\phi_S(x) & \mbox{if $x\in S^3_{1,+}$.}
\end{array} \right.
\end{array} $$
\edf 

The main property of those projective models is that they send geodesic
segments to segments of $\R^3$. $\phi_H$ sends a point at distance
$\rho$ from $x_0$ in $H^3$ to a point at distance $\tanh(\rho)$ from $0$
in $\R^3$. We call
{\bf radial} the direction of the geodesic from $x$ to $x_0$, and {\bf
  lateral} the directions orthogonal to the radial direction. 
Then $d_x\phi_H$:
\begin{itemize}
\item sends the radial direction to the radial direction, and the
  lateral directions to the lateral directions in $\R^3$.
\item "contracts" the lateral directions by a factor $1/\cosh(\rho)$,
  i.e. it multiplies the metric in those directions by $1/\cosh^2(\rho)$.
\item "contracts" the radial direction by a factor $1/\cosh^2(\rho)$,
  i.e. it multiplies the metric in those directions by $1/\cosh^4(\rho)$.
\end{itemize}
The description of $\phi_S$ can be done along similar lines, with the
distance to $x_0$ replaced by the distance to the totally geodesic plane
dual to $x_0$ (see below for the definition of the duality) and the
$\cosh$ essentially replaced by a $\sinh$. 

We will also need the following
rather explicit description of the properties of the geodesic spheres. 

\bprop \label{pr:spheres}
For each $\rho\in \R_+$, call $I_\rho$ and $\II_\rho$ the induced metric
and second fundamental form of the sphere of radius $\rho$ in $H^3$;
for each $t\in (0,1)$, call $\Ib_t$ and $\IIb_t$ the induced metric and
second fundamental form of the sphere of radius $t$ in $\R^3$. Then,
for $t=\tanh(\rho)$, we have:
$$ I_\rho = \cosh^2(\rho) \Ib_t~, ~~ \II_\rho = \cosh^2(\rho)\IIb_t~. $$ 
\eprop

\bpv
Let $\can_{S^2}$ be the canonical metric on $S^2$, then $I_\rho =
\sinh^2(\rho) \can_{S^2}$, $\II_\rho = \sinh(\rho)\cosh(\rho) \can_{S^2}$,
$\Ib_t=t^2 \can_{S^2}$, $\IIb_t = t\can_{S^2}$; the result follows.
\epv

The important point here is that the scaling of $I$ and $\II$ are the
same. The same phenomenon also happens for the de Sitter space.

\bprop \label{pr:spheres-S}
For each $t>1$, let $\Sb_t\subset \R^3$ be the sphere of radius $t$ and
center $0$, let $\rho:=\coth^{-1}(t)$, and let
$S_\rho:=\phi_S^{-1}(\Sb_t)$. Let $I_\rho, \II_\rho, \Ib_t, \IIb_t$ be
the induced metrics and second fundamental forms on $S_\rho$ and $\Sb_t$
respectively. Then:
$$ I_\rho = \sinh^2(\rho) \Ib_t~, ~~ \II_\rho = \sinh^2(\rho)\IIb_t~. $$ 
\eprop

\bpv
It is a simple matter to check that:
$$ I_\rho = \cosh^2(\rho) \can_{S^2}~, ~ \Ib_t = t^2\can_{S^2} =
\coth^2(\rho)\can_{S^2}~, $$
$$ \II_\rho = 
\cosh(\rho)\sinh(\rho) \can_{S^2}~, ~\IIb_t = t\can_{S^2}
=\coth(\rho)\can_{S^2}~, $$ 
and the result follows.
\epv

\paragraph{Principal curvatures and projective transformations}

An important point later on is that given a compact hyperbolic manifold with
smooth, strictly convex boundary, the image in the projective model of its
universal 
cover is a convex set bounded by a surface which has its principal
curvatures bounded between two positive constants. This will follow from
the following Lemma.

\blm \label{lm:reg}
Let $\Omega\subset H^3$ be a complete, non-compact convex domain containing
$x_0$. Suppose that the principal curvatures of $\dr \Omega$ are between
two constants $c_m$ and $c_M$, $c_M>c_m>0$, and that the norm of the
gradient of the distance 
to $x_0$ (restricted to $\dr\Omega$) is 
at most $1-\epsilon$ for some fixed $\epsilon>0$. Let $R>0$. Then
$\phi(\Omega)\subset D^3$ is a convex 
domain and, outside the ball of radius $R$, the principal
curvatures of the boundary are bounded between $c_m$ and 
another constant $c'_M$ depending on $c_M$, $R$ and $\epsilon$.
\elm

The proof of this lemma is based on the fact that, under a projective
map, the second fundamental form of surfaces changes conformally, as
stated in the following proposition. It is of course also valid in
higher dimension, and could be stated in higher codimension too.

\bprop \label{pr:proj-II}
Let $p:M\rightarrow \Mb$ be a projective diffeomorphism between two
3-dimensional space-forms. Let $S\subset M$ be a surface, and let
$\Sb:=p(S)$. Let $N, \Nb$ be the unit normal vectors of $S$ and $\Sb$,
and let $\II$ and $\IIb$ be the second fundamental forms of $S$ and
$\Sb$, respectively. Then:
$$ \IIb = \langle \Nb, p_*N\rangle_{\Mb} p_*\II = \frac{1}{\langle N,
p^*\Nb\rangle_M} p_*\II~. $$
\eprop

\bpv
Let $x\in S$, and let $\xb:=p(x)$. 
Let $P$ and $\Pb$ be the totally geodesic planes tangent to $S$ at $x$
and to $\Sb$ at $\xb$, respectively --- such planes exist since $M$ and
$\Mb$ are space-forms. By definition of a projective transformation,
$\Pb=p(P)$. 

For each $y\in P$ close to $x$, consider the geodesic $g_y$ orthogonal to $P$
starting from $y$, parametrized at speed $1$. Let $u(y)$ be such that
$g_y(u(y))\subset S$. Then $\II$ at $x$ can be defined as the Hessian at $x$
of $u$. Note that the direction of $g_y$ does not really matter, because
$S$ is tangent to $P$ at $x$, and we are only interested in the
second-order behavior of $u$ at $x$; changing the direction of $g_y$
changes only higher-order terms, if the orthogonal projection of $g_y'(0)$
on the direction orthogonal to $P$ has unit norm.

We can take this definition to $\Mb$ using $p$, we obtain that $p_*\II$
is the Hessian at $\xb$ of the function defined as the parameter, at the
intersection with $\Sb$, of the geodesic $\gb_{\yb}$ starting from
$\yb\in \Pb$, with $\gb_{\yb}'(0)$ equal
to the image by $p_*$ of the unit normal vector of $P$. (The function
considered here is different from $u\circ p^{-1}$, because the
parametrization of the geodesics is different, but the Hessian at
$\xb$ remains the same).

Now the definition of $\II$ described above applies also to $\IIb$, and
it shows that $\IIb$ differs from $p_*\II$ only by a factor equal to
$\langle p_*N,\Nb\rangle_{\Mb}$. The first equality follows. The second
equality can be checked by a simple computation, or by applying the
first formula to $p^{-1}$.  
\epv

\begin{proof}[Proof of Lemma \ref{lm:reg}]
The gradient of the restriction to $\dr\Omega$ of the distance $\rho$ to $x_0$
is bounded by $1-\epsilon$. Therefore, the angle of $T\dr\Omega$ with
the radial direction is bounded from below by some angle
$\alpha_\epsilon>0$. 

The differential $d\phi_H$ acts by contracting the lateral directions by
$1/\cosh(\rho)$, and the radial directions by a factor
$1/\cosh^2(\rho)$. Therefore, for $\rho$ large, $T\phi_H(\dr\Omega)$ is
"almost" parallel to the lateral directions; the angle between its
normal vector $\Nb$ and the radial direction is at most
$C_\epsilon/\cosh(\rho)$, for some constant $C_\epsilon>0$. 
If we call $I$ the induced
metric on $\dr\Omega$ and $\Ib$ the induced metric on
$\phi_H(\dr\Omega)$, we have:
$$ \frac{1}{\cosh^2(\rho)}(\phi_{H})_*I\leq \Ib \leq
\frac{C_\epsilon}{\cosh^2(\rho)}(\phi_{H})_*I~. $$
 
Moreover, $(\phi_H)_*N$ is easily estimated by the action of $d\phi_H$;
its lateral component 
has length at most $1/\cosh(\rho)$, and the length of its radial
component is between $1/(C_\epsilon \cosh^2(\rho))$ and $1/\cosh^2(\rho)$. 
Those bounds on the radial and lateral components of $(\phi_H)_*N$,
along with the bound on the lateral component of $\Nb$,
show that:
$$ \frac{1}{C_\epsilon\cosh^2(\rho)} \leq \langle (\phi_H)_*N, \Nb\rangle
\leq\frac{C_\epsilon}{\cosh^2(\rho)}~. $$

Proposition \ref{pr:proj-II} then shows that:
$$ \frac{1}{C_\epsilon\cosh^2(\rho)}(\phi_{H})_*\II\leq \IIb \leq
\frac{C_\epsilon}{\cosh^2(\rho)}(\phi_{H})_*\II~. $$ 
Thus $I$ and $\II$ scale by the same factor, up to a coefficient
$C_\epsilon^2$, and so the principal curvatures of $\dr\Omega$ and
$\phi_H(\dr\Omega)$ are the same, up to the same factor $C^2_\epsilon$. 
\end{proof}

\pg{The hyperbolic-de Sitter duality}

The study of the third fundamental form of surfaces in $H^3$ relies
heavily on an important duality between $H^3$ and the de Sitter space
$S^3_1$. Although this duality has been known for some time (see
e.g. \cite{thurston-notes}, Chapter 2), it was first used in a deep and
fundamentally new way by Rivin and Hodgson \cite{Ri,RH};
other approaches to this duality can be found elsewhere, e.g. in
\cite{shu}.  
It associates to each point $x\in H^3$ a space-like, totally
geodesic plane in $S^3_1$, and to each point $y\in S^3_1$ an oriented
totally geodesic plane in $H^3$. 

It can be defined using the quadric models of $H^3$ and $S^3_1$. Let
$x\in H^3$; define $d_x$ as the line in $\R^4_1$ going through $0$ and
$x$. $d_x$ is a time-like line; call $d_x^*$ the orthogonal space in
$\R^4_1$, which is a space-like 3-plane. So $d_x^*$ intersects $S^3_1$
in a space-like totally geodesic 2-plane, which we call $x^*$, and which
is the dual of $x$. Conversely, given a space-like totally geodesic
plane $p\subset S^3_1$, it is the intersection with $S^3_1$ of a space-like
3-plane $P\ni 0$. Let $d$ be its orthogonal, which is a time-like line;
the dual of $p$ is the intersection $d\cap H^3$.

The same construction works in the opposite direction. Given a point
$y\in S^3_1$, we call $d_y$ the oriented line going through $0$ and $y$,
and $d_y^*$ its orthogonal, which is an oriented time-like 3-plane. Then
$y^*:=d_y^*\cap H^3$ is an oriented totally geodesic plane. Forgetting
the orientation yields a duality between planes in $H^3$ and points in
$S^3_{1,+}$.

We can then define the duality on surfaces. Given a smooth, oriented
surface $S\subset H^3$, its dual $S^*$ is the set of points in $S^3_1$
which are the duals of the oriented planes which are tangent to
$S$. Conversely, given a smooth, space-like surface $\Sigma\subset
S^3_1$, its dual is the set $\Sigma^*$ of points in $H^3$ which are the
duals of the planes tangent to $\Sigma$. 

\bprop \label{pr:duality}
\begin{itemize}
\item If $S\subset H^3$ (resp. $S^3_1$) is smooth and strictly convex
  (resp. smooth, space-like and strictly convex), its dual is smooth and
  strictly convex. 
\item The induced metric on $S^*$ is the third fundamental form of $S$,
  and conversely.
\item For any smooth, strictly convex surface $S\subset H^3$
  (resp. space-like surface $S\subset S^3_1$), $(S^*)^*=S$. 
\end{itemize}
\eprop

The proof of the various points, in the smooth or the polyhedral context, can
be found e.g. in \cite{RH,these,iie,shu,cpt,dap}. 

There is a purely projective definition of this duality, or more
precisely of the closely related duality between $H^3$ and the de Sitter
hemisphere $S^3_{1,+}$; the duality is now between points in $S^3_{1,+}$
and non-oriented planes in $H^3$. To define it, let $p$ be a totally
geodesic plane in $H^3$; its image $\phi_H(p)$ is a disk $d_p$ in $D^3$. The
image by $\phi_S$ of the dual of $p$ is the unique point $p^*$ in $\R^3$
such that the lines going through $p^*$ and $\dr d_p$ are tangent to
$\dr D^3$. Using it, we can --- and we sometimes will --- use the
duality between surfaces in $D^3$ and in $\R^3\setminus D^3$.

The duality
between $H^3$ and $S^3_{1,+}$ is not primarily about surfaces, but about
couples $(x,P)$, where $P\subset T_xH^3$ is a plane containing $x$. 
Given such a
couple, we can define its dual as a couple $(x^*,P^*)$, where
$x^*\subset S^3_{1,+}$ is the dual of $P$, and $P^*\subset T_{x^*}S^3_{1,+}$
is the dual of $x$.

Let $(x,P)$ be such a couple, and let $(x^*, P^*)$ be its dual. Given an
element $u\in so(3,1)$, it defines a Killing vector field $X_u$ on $H^3$, and
a Killing vector field $X^*_u$ on $S^3_1$. Moreover, the first-order
deformation of $(x,P)$ under the Killing field $X_u$ induces, through the
duality, a first-order displacement $Y\in T_{x^*}S^3_{1,+}$ of $x^*$. 

\brk \label{rk:dual}
$Y=X^*(x^*)$. 
\erk

\bpv
Consider $x$, $x^*$ and the normal vectors $N$ and
$N^*$ to $P$ and $P^*$, respectively,
as vectors in $\R^4_1$; then $N=x^*$, while $x=N^*$. Moreover, if we consider
$u$ as a linear map on $\R^4_1$, we see that
$X(x)=u(x)$ as a vector in $\R^4_1$. The first-order
variation of $x^*=N$ under $X_u$ is just $u(N)=u(x^*)$, so that it is equal to
$X^*(x^*)$. 
\epv

\pg{The Pogorelov map}

There is an interesting extension of the projective model of $H^3$ (and
of $S^3_{1,+}$). It was first defined and used by Pogorelov \cite{Po},
at least for $H^3$. Pogorelov actually defined a slightly more
complicated map, between $H^3\times H^3$ and $\R^3\times \R^3$, but the
version which we will use here is simply its linearization along the
diagonal. The reader is refered to \cite{shu} for another, more
geometric approach of those questions. A slightly different but similar
map is also used in \cite{iie,cpt}.

\bdf \label{df:pogo-H}
We define $\Phi_H:TH^3\rightarrow TD^3$ as the map sending $(x,v)\in
TH^3$ to $(\phi_H(x), w)$, where:
\begin{itemize}
\item the lateral component of $w$ is the image by $d\phi_H$ of the
lateral component of $v$.
\item the radial component of $w$ has the same direction and the same
norm as the radial component of $v$.
\end{itemize}
\edf

Here are the main properties of this map. The proof which we give here
is elementary; it might appear simpler to what can be found e.g. in
\cite{shu,iie}, where the arguments are of a more geometric nature.

\blm \label{lm:pogo-H}
Let $S$ be a smooth submanifold in $H^3$, and let $v$ be a vector field
of $H^3$ defined on $S$; then $v$ is an isometric deformation of $S$ if
and only if $\Phi_H(v)$ is an isometric deformation of $\phi_H(S)$.
In particular, if $v$ is a vector field on $H^3$, then $v$ is a Killing
field if and only if $\Phi_H(v)$ is a Killing field of $D^3$.
\elm

\bpv
The second part follows from the first by taking $S=H^3$.

To prove the first part, we have to check that the Lie derivative of the
induced metric on $S$ under $v$ vanishes if and only if the Lie
derivative of the induced metric on $\phi_H(S)$ under $\Phi_H(v)$
vanishes. In other terms, if $x,y$ are vector fields defined on $S$,
if we call $g$ and $\gb$ the metrics on $H^3$ and $\R^3$ respectively,
and $\xb=d\phi_H(x), \yb=d\phi_H(y)$,
then we have to prove that:
$$ (\cL_vg)(x,y)=0\Leftrightarrow (\cL_{\Phi_H(v)}\gb)(\xb, \yb)=0~. $$
We will actually prove that the two sides are proportional:
\beq \label{eq:pogo} 
(\cL_vg)(x,y)=\cosh^2(\rho)(\cL_{\Phi_H(v)}\gb)(\xb, \yb)~. 
\eeq
We decompose $v$ into the radial component, $fN$, where $N$ is the unit
radial vector, and the lateral component $u$. 
Abusing notations a little, we also call $f$ the function $f\circ
\phi_H^{-1}$, defined on the unit ball in $\R^3$. 
Then $\Phi_H(v)=f\Nb +\ub$, where $\Nb$ is the unit radial 
vector in $\R^3$, and $\ub:=d\phi_H(u)$.

By linearity, it is sufficient to prove equation (\ref{eq:pogo}) in the
cases where $x$ and $y$ are non-zero, and each is either radial or
lateral. We consider each case separately, and call $\nabla$ and $\nablab$
the Levi-Civit{\`a} connections of $H^3$ and $\R^3$ respectively.

\pg{1st case: $x$ and $y$ are both radial.} Then:
$$ (\cL_vg)(x,y) = (\cL_ug)(x,y) + (\cL_{fN}g)(x,y)~. $$
But:
$$ (\cL_ug)(x,y) = g(\nabla_xu,y) + g(x,\nabla_yu) = 0~, $$
because both $\nabla_xu$ and $\nabla_yu$ are lateral. 
Moreover:
$$ (\cL_{fN}g)(x,y) = g(\nabla_x fN, y) + g(x, \nabla_y fN) = df(x) g(N,y) +
df(y) g(x,N)~. $$
The same computation applies in $\R^3$. In addition, $df(x)=df(\xb)$ by
definition of $\xb$, so the scaling comes only from $g(x,N)$ versus
$\gb(\xb, N)$. So equation
(\ref{eq:pogo}) holds in this case. 

\pg{2nd case: $x$ and $y$ are both lateral.}
Then, at a point at distance $\rho$ from $x_0$, with $t:=\tanh(\rho)$,
we have with Proposition \ref{pr:spheres}:
$$ (\cL_ug)(x,y) = (\cL_uI_\rho)(x,y) = \cosh^2(\rho)
(\cL_{\ub}\Ib_t)(\xb,\yb) = \cosh^2(\rho)(\cL_{\ub} \gb)(\xb,\yb) ~. $$
Moreover:
$$ (\cL_{fN}g)(x,y) = g(\nabla_x fN, y) + g(x, \nabla_y fN) = -2f \II_\rho (x,y)~,
$$
and since the same computation applies in $\R^3$, we see with
Proposition \ref{pr:spheres} that:
$$ (\cL_{fN}g)(x,y) = \cosh^2(\rho) (\cL_{fN}\gb)(\xb,\yb)~, $$
and:
$$ (\cL_vg(x,y) = \cosh^2(\rho) (\cL_{\vb}\gb)(\xb, \yb)~, $$
so that equation (\ref{eq:pogo}) also holds in this case. 

\pg{3rd case: $x$ is lateral, while $y$ is radial.} 

We choose an arbitrary extension of $x$ and $y$ as vector fields which
remain tangent, resp. orthogonal, to the spheres $\{
\rho=\mbox{const}\}$. Then:
$$ (\cL_ug)(x,y) = u.g(x,y) - g([u,x],y) - g(x,[u,y]) = -g(x,[u,y])~, $$
because $[u,x]$ is tangent to the sphere of radius $\rho$. 
So, since the Lie bracket does not depend on the metric, we have with
Proposition \ref{pr:spheres} that:
$$ (\cL_ug)(x,y) = - g(x,[u,y]) = - \cosh^2(\rho) \gb(\xb,[\ub,\yb]) =
\cosh^2(\rho)(\cL_{\ub}\gb)(\xb,\yb)~. $$  
In addition:
$$ (\cL_{fN}g)(x,y) = g(\nabla_x fN, y) + g(x, \nabla_y fN) = df(x)
g(y,N) = $$ 
$$ = df(x) \cosh^2(\rho) \gb(\Nb,\yb) = \cosh^2(\rho)
(\cL_{f\Nb}\gb)(\xb,\yb)~, $$ 
again because the same computation applies in $\R^3$. As a consequence,
we have again that:
$$ (\cL_vg)(x,y) = \cosh^2(\rho) (\cL_{\vb}\gb)(\xb, \yb)~, $$
and (\ref{eq:pogo}) still holds. 
\epv

\brk
\begin{itemize}
\item the proof only uses Proposition \ref{pr:spheres}, more precisely
the fact that $I$ and $\II$ scale in the same way when one goes from
$H^3$ to $\R^3$ by $\phi_H$.
\item the same proof can be used in different dimensions, and also in
spaces of other signatures. 
\end{itemize}
\erk

In particular, we can then define in the same way a Pogorelov map for the de
Sitter space, as done elsewhere, see e.g. \cite{dap,rcnp}. Note that there is
another related possible map for the de Sitter space, with values in the
Minkowski space, see \cite{gajubov,iie}.

\bdf \label{df:pogo-S}
We define $\Phi_S:TS^3_{1,+}\rightarrow T(\R^3\setminus D^3)$ as the map
sending $(x,v)\in 
TS^3_{1,+}$ to $(\phi_S(x), w)$, where:
\begin{itemize}
\item the lateral component of $w$ is the image by $d\phi_S$ of the
lateral component of $v$.
\item the radial component of $w$ has the same direction and the same
norm as the radial component of $v$.
\end{itemize}
\edf

\blm \label{lm:pogo-S}
Let $S$ be a smooth submanifold in $S^3_{1,+}$, and let $v$ be a vector field
of $S^3_{1,+}$ defined on $S$; then $v$ is an isometric deformation of $S$ if
and only if $\Phi_S(v)$ is an isometric deformation of $\phi_S(S)$.
In particular, if $v$ is a vector field on $S^3_{1,+}$, then $v$ is a Killing
field if and only if $\Phi_S(v)$ is a Killing field of $\R^3$.
\elm

The proof is the same as for Lemma \ref{lm:pogo-H}, using Proposition
\ref{pr:spheres-S} instead of Proposition \ref{pr:spheres}.

Given an element $u\in so(3,1)$, we have already seen that it defines a
Killing vector field $X_u$ on $H^3$ and a Killing field $X^*_u$ on $S^3_1$. An
important remark, which follows directly from the definition of the Pogorelov
maps for $H^3$ and $S^3_1$, is that the images are ``compatible''.

\brk \label{rk:recollement}
$\Phi_H(X_u)$ and $\Phi_S(X^*_u)$ are the restriction to $D^3$ and to
$\R^3\setminus D^3$, respectively, of the same Killing field of $\R^3$.
\erk

\pg{Rigidity as a projective property}

A direct but striking consequence of the previous paragraph is that
infinitesimal rigidity is a purely projective property. 

\blm \label{lm:rig-poly}
Let $P$ be a smooth surface (resp. a polyhedron) in $H^3$
(resp. $S^3_{1,+}$). $P$ is 
infinitesimally rigid if and only if $\phi(P)\subset \R^3$ is
infinitesimally rigid.
\elm

\bpv
Suppose that $P\subset H^3$ is not infinitesimally rigid. 
Let $v$ be a non-trivial infinitesimal isometric deformation of $P$ in
$H^3$. Then $\Phi_H(v)$ is a non-trivial infinitesimal deformation of
$\phi_H(P)$ by Lemma \ref{lm:pogo-H}. The same works in the other
direction, and also in $S^3_{1,+}$ by Lemma \ref{lm:pogo-S}.
\epv

The proof also
works for geometric objects that are more general than compact
polyhedra in $H^3$ or space-like compact polyhedra in $S^3_1$, this is
used for instance in \cite{shu,cpt,dap}. 

One can formulate results like Lemma \ref{lm:rig-poly} as
follows. Consider a polyhedron $P\in \R^3$, 
and consider $\R^3$ as a projective space (one could as well replace
$\R^3$ by $\rp^3$, but this would slightly cloud the issue). One can
then put a constant curvature metric on $\R^3$, which is compatible
with the projective structure: it could be a Euclidean metric, or a
spherical metric (isometric to a hemisphere of $S^3$), or the
structure corresponding to hyperbolic 3-space in a ball and $S^3_{1,+}$
outside. Then the fact that $P$ is, or is not, infinitesimally rigid, is
independent of the constant curvature metric chosen. In particular,
convex polyhedra --- and this is a projective property --- are always
infinitesimally rigid, a result of Legendre \cite{legendre} and Cauchy
\cite{cauchy}.

\pg{Necessary conditions on $\III$}

One should remark that the length condition in Theorem \ref{tm:III} is
necessary. This is a consequence of the next lemma.

\blm \label{lm:2pi}
Let $S$ be a strictly convex surface in $H^3$, and let $S^*$ be the dual
surface in $S^3_1$. The closed geodesics of $S^*$ have length $L>2\pi$. 
\elm

We refer the reader to \cite{CD,RH,these,shu} for various proofs. We can
also remark that there is a dual condition for surfaces in $S^3_1$ which
are not space-like, see \cite{cpt}.


\section{Hyperbolic manifolds with boundary}

After defining, in the previous section, the basic tools that we will
need, we will in this section study the hyperbolic and de Sitter
surfaces that can be associated to a hyperbolic metric on $M$.

\pg{The extension of a manifold with convex boundary}

First remark that, given a hyperbolic manifold with convex boundary,
there is essentially a unique way to embed it isometrically into a
complete hyperbolic manifold. 

\bdf \label{df:extension}
Let $M$ be a compact hyperbolic manifold with convex boundary. There
exists a unique complete, convex co-compact hyperbolic manifold $E(M)$,
which we call the {\bf extension} of $M$, 
in which $M$ can be isometrically embedded so that the natural map
$\pi_1(M)\rightarrow \pi_1(E(M))$ is an isomorphism. 
\edf

As a consequence, $\pi_1M$ has a natural action on $H^3$.

\pg{Surfaces in the projective models}

Let $g$ be a hyperbolic metric with smooth, strictly convex boundary on
$M$. $g$ defines a metric, which we still call $g$, on the universal
cover $\Mt$ of $M$. Since $(\Mt, g)$ is simply connected and has convex
boundary, it admits a unique (up to global isometries of $H^3$)
isometric embedding into 
$H^3$, and the image is a convex subset $\Omega$, with smooth, strictly
convex boundary. Of course $\Omega$ is not compact; its boundary at
infinity is the limit set $\Lambda$ of the action of $\pi_1M$ on $H^3$
which is obtained by considering the extension of the action of
$\pi_1(M)$ on $\Omega$. 

In the sequel, we will call $S:=\dr \Omega$ the boundary of the image of
$\Mt$ in $H^3$. So $S$ is a complete, smooth, strictly convex surface in
$H^3$, with $\dr_\infty S=\Lambda$. We will also call $S^*$ the dual
surface in the de Sitter space $S^3_1$; $S^*$ is a space-like,
metrically complete, smooth, strictly convex surface in $S^3_{1,+}$, and
its boundary at infinity is also $\Lambda$.

\pg{Boundary behavior}

It is then necessary to understand more precisely the behavior of $S$
near the boundary, in particular in the projective model of $H^3$; we
choose the center $x_0$ of the projective model as a point in $\Omega$. 
Recall that a {\bf geodesic ray} starting from $x_0$ is a half-geodesic
$\gamma:\R_+\rightarrow H^3$ with $\gamma(0)=x_0$. 

\bprop \label{pr:dist-omega}
There exists a constant $c_0>1$ such that:
\begin{itemize}
\item for any geodesic ray $\gamma$
starting from $x_0$ and ending on $\Lambda$, the distance from the
points of $\gamma$ to $S$ varies between $1/c_0$ and $c_0$.
\item for any point $x\in \dr \Omega$, there exists a geodesic ray
$\gamma\subset \Omega$ starting from $x_0$ and ending on $\Lambda$ 
which is at distance at most $c_0$ from $x$. 
\end{itemize}
\eprop

\bpv
Since $\dr M$ is compact and strictly convex, its principal curvatures
are bounded from below by a strictly positive constant $\kappa$. An
elementary geometric argument then shows that there exists a constant
$\delta>0$ such that any maximal geodesic segment in $M$ going through a
point $x\in M$ at distance at most $\delta$ from $\dr M$ meets $\dr M$
at distance at most $1/\delta$ from $x$. This implies the existence of a
lower bound on the distance from the geodesic rays to $\dr M$. The upper
bound is obvious since $M$ is compact, so has a bounded diameter. 

For the second point note that each point $x\in \dr M$ is at bounded
distance from a point $x'_0$ in the compact core $C(M)$ of $M$; and each
point $x'_0\in \dr C(M)$ is either in a complete geodesic which is on
$\dr C(M)$, or in an ideal triangle on $\dr C(M)$. In both cases $x'_0$ is
at bounded distance from a geodesic $c$ which is on $\dr
C(M)$, and which lifts 
to a complete geodesic in $\Omega$ with both endpoints on
$\Lambda$. Going to the universal cover, we see that there exists a
constant $C$ such that each point of $\Omega$ is at distance at most $C$
from a complete geodesic $c$ having both endpoints on $\Lambda$.

Let $x_1$ and $x_2$ be the endpoints of $c$ in
$\Lambda$. There are geodesic rays $c_1$ and $c_2$ starting at $x_0$ and
ending on $x_1$ and $x_2$, respectively. A basic property of hyperbolic
triangles is that, for some universal constant $C'$, each point of $c$ is
at distance at most $C'$ from either $c_1$ or $c_2$. Therefore, each
point of $\Omega$ is at distance at most $C+C'$ from a ray starting at
$x_0$ and ending on $\Lambda$.
\epv

The convexity of $\dr M$ then allows us to go from the previous uniform
estimate to a $C^1$ estimate on the behavior of $S$. 

\bcr \label{cr:gradient}
On $S$, the gradient of the distance $\rho$ to $x_0$ is bounded from
above by a constant $\delta_0<1$.
\ecr

\bpv
Let $x\in \dr \Omega$. By Proposition \ref{pr:dist-omega}, there exists
a geodesic ray $\gamma_0$ starting from $x_0$, at distance at most $c_0$
from $x$. Let $P$ be the totally geodesic plane containing 
$\gamma_0$ and $x$. By construction, $P\cap \Omega$ is a convex,
non-compact subset of $P\simeq H^2$.

Let $\gamma$ be the geodesic in $P$ 
tangent to $\dr (P\cap \Omega)$ at $x$. By
convexity, $\dr (P\cap \Omega)$ remains on the side of $\gamma$
containing $\gamma_0$. Let $\gamma_1$ be the geodesic ray starting at
$x_0$ and going through $x$. An elementary fact of plane hyperbolic
geometry is that there exists a constant $\delta_0(c_0)$ such that, if the
angle between $\gamma_1$ and $\gamma$ at $x$ is smaller than
$\delta(c_0)$, then some point of $\gamma$ is at distance smaller than
$1/c_0$ from $\gamma_0$. So some point of $\dr \Omega$ would have to be at
distance smaller than $1/c_0$ from $\gamma_0$, and this would contradict
the first part of Proposition \ref{pr:dist-omega}.
Thus the angle between $\gamma_1$ and $\gamma$ is at least
$\delta_0(c_0)$.
\epv

We will call $\Sb$ the image by
$\phi_H$ of $S$, and we will also call $\Lambda$ --- with some abuse of
notations --- the image on the unit sphere $S^2\subset \R^3$ of the
limit set of $\pi_1M$. By construction, $\Sb\cup
\Lambda=\dr \phi_H(\Omega)$ is a convex, closed surface in $\R^3$ --- in
particular it is 
topologically a sphere --- and we will see below that it has some degree
of smoothness even on $\Lambda$.

Using the projective model to go from $H^3$ to $D^3$, the corollary above
translates into a simple estimate on the gradient of the distance $r$ to $0$
on $\Sb$. We state it in the next
lemma, along with some important estimates on the behavior on $\Sb$ of
the distance along $\Sb$ to the limit set $\Lambda$, which we will call
$\delta$ here and later on. 

\blm \label{lm:slope}
There exists a constant $c_1>1$ such that, at all points of $\Sb$ where
the distance $r$ to $0$ is close enough to $1$, we have:
\begin{enumerate}
\item the angle $\theta$ between the normal to $\Sb$ and the radial
  direction is bounded from above by: $\theta \leq c_1 \sqrt{1-r}$.
\item the distance $\delta$ on $\Sb$ to the limit set $\Lambda$
  satisfies: 
$$ \frac{\sqrt{1-r}}{c_1}\leq d_{\Sb}(x, \Lambda)\leq c_1\sqrt{1-r}~. $$
\end{enumerate}
\elm

\bpv
Consider a point $x\in H^3$ far enough from $x_0$. Let $u\in T_xS$ be
the unit vector on which $d\rho$ is maximal. Then $d\rho(u)\leq 1$. But
$r=1-\tanh(\rho)$, so that $|dr(u)|\leq \cosh^{-2}(\rho)$. 
Let $\xb:=\phi_H(x)$ and $\ub:=d\phi_H(u)\in T_{\xb}\Sb$. Then, by
Corollary \ref{cr:gradient} and the definition of $\phi_H$, there
exists a constant $c>0$ (depending 
on the constant $\delta_0$ of Corollary \ref{cr:gradient}) such that
$\|\ub\|_{D^3} \geq c/\cosh(\rho)$. Moreover, it is clear that $\ub$ is in
the direction of greatest variation of $r$ on $\Sb$. Therefore, for
another constant $c>0$, $\| dr_{|\Sb}\|\leq c/\cosh(\rho)$. Since
$r=\tanh(\rho)$, $1/\cosh(\rho)=\sqrt{1-r^2}$, and  
this proves the first point. 

For the second point, let $\xb\in \Sb$, and let $x:=\phi_H^{-1}(\xb)\in
H^3$. Let $\gamma$ be a geodesic ray
starting at $x_0$ and ending on $\Lambda$, such that the distance from
$x$ to $\gamma$ is minimal. According to Proposition
\ref{pr:dist-omega}, this distance is between $1/c_0$ and $c_0$.

Let $\gammab:=\phi_H(\gamma)$, then $\gammab$ is a segment starting from
$0$ and ending on $\Lambda\subset S^2$, and among those segments it is
at the smallest distance from $\xb$. If the distance $\rho$ from $x_0$
to $x$ is large enough, the properties of the projective model show that:
$$ \frac{1}{2c_0\cosh(\rho)} \leq d(\gammab, \xb)\leq
\frac{2c_0}{\cosh(\rho)}~. $$
Since $1/\cosh(\rho) = \sqrt{1-r^2} =\sqrt{(1-r)(1+r)}$, it follows
that:
$$ \frac{\sqrt{1-r}}{2c_0} \leq d(\gammab, \xb)\leq
4c_0\sqrt{1-r}~, $$
and the estimate on $\delta$ follows.
\epv

An important consequence of the first point is that, as $\xb\rightarrow
\xb_\infty\in \Lambda$ on $\Sb$, the unit normal vector to $\Sb$ at $\xb$
converges to the unit normal vector to $S^2$ at $\xb_\infty$. Thus $\Sb$ is
"tangent" to $S^2$ along $\Lambda$.

\pg{Limit curvatures}

Using the previous estimates and the results of the previous section on
the behavior of the projective model, we can find some higher-order
estimates on $\Sb$, concerning its principal curvatures in the
neighborhood of $\Lambda$.

\blm \label{lm:curvature}
There exists a constant $c_2>1$ such that, at all points of $\Sb$ where
the distance $r$ to $\dr D^3$ is small enough, the principal curvatures
of $\Sb$ are bounded between $c_2$ and $1/c_2$.
\elm

\bpv
Since $\dr M$ is compact and strictly convex, its principal curvatures
are bounded between two constants $c_m$ and $c_M$, with $c_M\geq
c_m>0$. The results then follows from Lemma \ref{lm:reg} and from
Corollary \ref{cr:gradient}.
\epv

\pg{The dual surface}

The results obtained in this section for $\Sb$ are also valid for the
dual surface $\Sb^*$. In particular, we have 
the following analogs of Lemma \ref{lm:slope} and
\ref{lm:curvature*}. 

\blm \label{lm:slope*}
There exists a constant $c_1^*>1$ such that, at all points of $\Sb^*$ where
the distance $r$ to $0$ is close enough to $1$, we have:
\begin{enumerate}
\item the angle $\theta$ between the normal to $\Sb^*$ and the radial
  direction is bounded from above: $\theta \leq c_1^* \sqrt{1-r}$.
\item the distance $\delta$ on $\Sb^*$ to the limit set $\Lambda$
  satisfies: 
$$ \frac{\sqrt{r-1}}{c^*_1}\leq d_{\Sb^*}(x, \Lambda)\leq c_1^*\sqrt{r-1}~. $$
\end{enumerate}
\elm

\bpv[Sketch of the proof]
Let $x\in \Sb$, and let $x^*$ be the dual point in $\Sb^*$, i.e. the
point in $\R^3$ which is the dual of the plane tangent to $\Sb$ at
$x$. By definition of the duality, the position of $x^*$ is determined
by the position of the tangent plane to $\Sb$ at $x$, while the tangent
plane to $\Sb^*$ at $x^*$ is determined by the position of $x$. The
first point is thus a consequence of both points of Lemma
\ref{lm:slope}. 

The second point is then a direct consequence of the first, as in the
proof of Lemma \ref{lm:slope}.
\epv

\blm \label{lm:curvature*}
There exists a constant $c^*_2>1$ such that, at all points of $\Sb^*$ where
the distance $r$ to $\dr D^3$ is small enough, the principal curvatures
of $\Sb^*$ are bounded between $c^*_2$ and $1/c^*_2$.
\elm

\bpv[Sketch of the proof] One possible way to prove this Proposition is
through the estimates of Lemma \ref{lm:curvature}, along with an
explicit computation of the principal curvatures of the dual $\Sigma^*$
of a convex surface $\Sigma\subset D^3$. 

Another way to prove it uses the same method as for
Lemma \ref{lm:curvature}, based on Proposition \ref{pr:proj-II} --- in
which the fact that the target manifold is Riemannian or Lorentzian
plays no role --- and on the de Sitter analog of Lemma \ref{lm:reg}. One
only needs to control the angle between $S^*$ and the radial direction
in $S^3_{1,+}$, which is precisely what Lemma \ref{lm:slope*} does.
\epv


\section{Rigidity of smooth convex surfaces in $\R^3$} \label{se:3}

We now 
consider closely the problem of infinitesimal rigidity of a convex surface in
the Euclidean (or the hyperbolic) 
3-space. In the next sections, we will obtain some estimates on
the deformations associated to a first-order deformation of a hyperbolic
metric on $M$, until, in section 6, we use those estimates and 
the ideas developped here
--- which are mostly classical --- to prove Lemma \ref{lm:rigidity}, dealing
with the infinitesimal rigidity of hyperbolic manifolds with convex boundary. 

It has been known for a long time --- at least since the works of
Cohn-Vossen \cite{CV} and Herglotz \cite{herglotz} --- that the smooth,
strictly convex closed surfaces in $\R^3$ are infinitesimally rigid,
i.e. they admit no smooth infinitesimal deformation which does not
change the metric at first order. The main rigidity result of this
paper, however is based on a small generalization of this statement,
to surfaces which are not so smooth and for deformations which might
also have singularities. It will therefore be necessary to have a good
understanding of some of the ways in which the classical infinitesimal rigidity
statement is proved. Note that, thanks to the statements of section 1,
it does not make much difference to consider a compact (strictly convex)
surface in $\R^3$, in $H^3$ or in $S^3_1$. We will however mostly remain
in $\R^3$, although it is important to keep in mind that almost all applies in
the same way in the hyperbolic 3-space. 

We will describe in this section two formulations of the
infinitesimal rigidity of smooth surfaces, in terms of some elliptic
equations. 
The main point here is that the rigidity of a surface in $\R^3$ or $H^3$ 
goes hand in hand with another problem, which for convex surfaces in $H^3$
describes the infinitesimal deformations of the dual surface in $S^3_1$,
and our main rigidity result (in
section 6) will crucially use this fact.
In this section, we will be less careful about the smoothness of the
functions or vector fields which enter the picture, although the
proofs used here are chosen so that they're not too demanding in terms
of regularity.

In all this section we will consider a smooth, strictly convex, closed
surface $S\subset \R^3$ (or $S\subset H^3$). 
An infinitesimal deformation of $S$ is a vector
field $u$ defined on $S$ --- formally, a section of the pull-back on $S$
by its immersion in $\R^3$ (or $H^3$) of the tangent bundle of $\R^3$ (or
$H^3$). It is an infinitesimal isometry if the induced first-order variation
of the induced metric of $S$ vanishes. $S$ is infinitesimally rigid if it has
no non-trivial infinitesimal isometric deformation. There are three
basic ways of stating this property:
\begin{enumerate}
\item by considering the component of $u$ which is normal to $S$; it is
  of the form $fN$, where $N$ is the unit normal vector field of
  $S$. If $u$ is an isometric deformation, $f$ is solution of a
  second-order elliptic PDE, whose principal symbol is given by the
  second fundamental form $\II$ of $S$.
\item in terms of the component $v$ of $u$ which is tangent to $S$; when
  $u$ is isometric, it is solution of a first-order elliptic PDE.
\item in terms of the first-order variation $\Bd$ of the shape operator
  $B$, which is also 
  solution of a simple PDE which is elliptic and of first order. 
\end{enumerate}
We will consider the formulations (2) and (3); the main remark of this
section is that those two formulations are "adjoint" in a natural
way. More precisely, each is "adjoint" to a "dual" of the other. 
This will be useful since, in the setting of the next section, the
formulation (3) can not directly be used by lack of {\it a priori} smoothness
of the deformations. The "duality" which appears here will allow us to
go from the formulation (2) to the formulation (3) and gain some
smoothness in the process. 

\pg{Extrinsic invariants of a convex surface}

The first point is to understand the first, second and third
fundamental forms of $S$, as well as their Levi-Civit{\`a} connections.
The first proposition is rather classical, see e.g. \cite{these,iie}.

\bprop \label{pr:III}
The Levi-Civit{\`a} connection of $\III$ is:
$$ \nablat_XY = B^{-1}\nabla_X(BY) = \nabla_XY + B^{-1}(\nabla_XB)Y~. $$
\eprop 

\bpv
By definition of the Levi-Civit{\`a} connection, it is sufficient to prove
that the connection $\nablat$ defined in the proposition is compatible with
$\III$ and has no torsion. But, if $X,Y$ and $Z$ are vector fields on $S$:
\begin{eqnarray*}
X.\III(Y,Z) & = & X.I(BY,BZ) \\
& = & I(\nabla_X(BY),BZ) + I(BY,\nabla_X(BZ)) \\
& = & \III(B^{-1}\nabla_X(BY),Z) + \III(Y,B^{-1}\nabla_X(BZ)) \\
& = & \III(\nablat_XY,Z) + \III(Y,\nablat_XZ)~,
\end{eqnarray*}
so that $\nablat$ is compatible with $\III$. 

Moreover:
\begin{eqnarray*}
\nablat_XY -\nablat_YX -[X,Y] & = & B^{-1}\nabla_X(BY) - B^{-1}\nabla_X(BZ) - [X,Y] \\
& = & B^{-1}(\nabla_X (BY)-\nabla_Y(BX)-B[X,Y]) \\
& = & B^{-1} (d^\nabla B)(X,Y) \\
& = & 0~,
\end{eqnarray*}
and thus $\nablat$ has no torsion. 
\epv

\bcr \label{cr:K1}
If $S$ is a surface in $\R^3$, the curvature of $\nablat$ is $\Kt=1$.
\ecr

\bpv
Let $(e_1, e_2)$ be an orthonormal moving frame for $I$, and let
$(\et_1, \et_2):=(B^{-1}e_1, B^{-1}e_2)$ be the corresponding
orthonormal moving frame for $\III$. A direct computation using the
previous proposition shows that the connection 1-forms of $(e_1, e_2)$
(for $I$) and of $(\et_1, \et_2)$ (for $\III$) are the same. Therefore,
their curvature 2-forms are also the same, so that:
$$ K da_I = \Kt da_{\III}~, $$
where $da_I$ and $da_{\III}$ are the area forms of $I$ and $\III$,
respectively. 
Thus the curvature $\Kt$ is such that:
$$ \Kt = \frac{K}{\det(B)}~, $$
so that $\Kt=1$  by the Gauss formula. 
\epv

Another possible proof uses the fact that $\III$ is the pull-back by the
Gauss map of the round metric on $S^2$. 

\bprop \label{pr:II}
The Levi-Civit{\`a} connection of $\II$ is:
$$ \nablab_XY = \frac{1}{2}(\nabla_XY + \nablat_XY) = \nabla_XY +
\frac{1}{2}B^{-1}(\nabla_XB)Y~. $$ 
\eprop 

\bpv
Again we have to prove that $\nablab$ is compatible with $\II$, and without
torsion. Let $X, Y$ and $Z$ be 3 vector fields on $S$. Since $B$ is
self-adjoint for $I$, so is $\nabla_XB$, so that:
\begin{eqnarray*}
X.\II(Y,Z) & = & X.I(BX,Y) \\
& = & I((\nabla_XB)Y,Z) + I(B\nabla_XY,Z) + I(BY,\nabla_XZ) \\
& = & \II(\nabla_XY,Z) + \II(Y,\nabla_XZ) + \frac{1}{2} I((\nabla_XB)Y,Z) + \frac{1}{2}
I(Y,(\nabla_XB)Z) \\
& = & \II(\nablab_XY,Z) + \II(Y,\nablab_XZ) \\ 
\end{eqnarray*}
so that $\nablab$ is compatible with $\II$.

The same computation as in the previous proposition shows that $\nablab$ has
no torsion:
\begin{eqnarray*}
\nablab_XY-\nablab_YX -[X,Y] & = & \nabla_XY -\nabla_YX -\frac{1}{2}
B^{-1}((\nabla_XB)Y-(\nabla_YB)X) -[X,Y] \\
& = & -\frac{1}{2}B^{-1}(d^\nabla B)(X,Y) + (\nabla_XY -\nabla_YX -[X,Y]) \\
& = & 0~. 
\end{eqnarray*}
\epv

A consequence of those propositions is the following formula, which will
be important in the sequel.

\bprop \label{pr:mixte}
Let $X,Y$ and $Z$ be 3 vector fields on $S$. Then:
$$ X.\II(Y,Z) = \II(\nabla_XY,Z) + \II(Y,\nablat_YZ)~. $$
\eprop

\bpv
\begin{eqnarray*}
X.\II(Y,Z) & = & \II(\nablab_XY,Z) + \II(Y,\nablab_XZ) \\
& = & \II(\nablab_XY,Z) - \frac{1}{2}I((\nabla_XB)Y,Z) + \frac{1}{2}I(Y,(\nabla_XB)Z)
+ \II(Y,\nablab_XZ) \\
& = & \II(\nablab_XY,Z) - \frac{1}{2}\II(B^{-1}(\nabla_XB)Y,Z) +
\frac{1}{2}\II(Y,B^{-1}(\nabla_XB)Z) + \II(Y,\nablab_XZ) \\
& = & \II(\nabla_XY,Z) + \II(Y,\nablat_XZ)~.
\end{eqnarray*}
\epv

We will call $\Jb$ the complex structure associated to $\II$. There is a
natural bundle of $(0,1)$-forms with values in 
the tangent bundle of $\Sb$.

\bdf 
We will call $\ots$ the sub-bundle of $T^*S\times TS$ of 1-forms
$\omega$ with values in $TS$ such that, for $X\in TS$, $\omega(\Jb
X)=-\Jb \omega(X)$. 
\edf 

The proof of the main rigidity result uses several differential
operators defined using the connections $\nabla$ and $\nablat$ on $\Sb$, in spite
of the fact that the metric which is mostly used is $\II$. 
Recall that the usual $\db$ operator for $\II$ is:
$$ (\db_{\II}u)(X) := \nablab_Xu + \Jb \nablab_{\Jb X}u~. $$
Several elliptic operators, acting on vector fields on $S$, have a
principal symbol defined by $\II$, in particular the following one,
already considered by Labourie \cite{L1}:

\bdf
We will call:
$$ (\db_Iu)(X) := \nabla_Xu + \Jb \nabla_{\Jb X}u~, $$
$$ (\db_{\III}u)(X) := \nablat_Xu + \Jb \nablat_{\Jb X}u~. $$
\edf

\blm \label{lm:elliptic}
$\db_I$ and $\db_{\III}$ are elliptic operators of index $6$ acting from
$TS$ to $\ots$. 
\elm

\bpv
If $u$ and $X$ are vector fields on $S$, we have:
\begin{eqnarray*}
(\db_Iu)(\Jb X) & = & \nabla_{\Jb X}u + \Jb \nabla_{\Jb^2 X} u \\
& = & -\Jb \nabla_Xu + \nabla_{\Jb X}u \\
& = & -\Jb (\nabla_X u + \Jb \nabla_{\Jb X}u)~,
\end{eqnarray*}
which shows that $\db_I$ has values in $\ots$. 
The same computation shows the same result for $\db_{\III}$.

The definition of $\db_I$ and $\db_{\III}$ shows that their principal
symbol is the same as the principal symbol of the operator obtained by
replacing $\nabla$ (resp. $\nablat$) by $\nablab$; this replacement leads to the
usual $\db$ operator for $\II$, which is well known to be elliptic of
index $6$. The lemma follows by the deformation invariance of the index.
\epv

\pg{Tangent vector fields}

The most natural way to consider the infinitesimal rigidity question for
surfaces is in terms of the infinitesimal deformation vector
fields. Given such a vector field, it is natural to consider its
projection on the surface.

\bprop
Let $u$ be an infinitesimal isometric deformation of $S$, and let $v$ be
its component tangent to $S$. Then the vector field $w:=B^{-1}v$ is
solution of:
\beq \label{eq:dbIII}
\db_{\III} w =0~. 
\eeq
\eprop

\bpv
Let $X$ be a vector field on $S$. 
\begin{eqnarray*}
(\db_{\III} w)X & = & \nablat_Xw + \Jb \nablat_{\Jb X}w \\
& = & B^{-1}\nabla_X(Bw) + \Jb B^{-1} \nabla_{\Jb X} (Bw) \\
& = & B^{-1}\nabla_Xv + \Jb B^{-1} \nabla_{\Jb X}v~.
\end{eqnarray*}
Since $u$ is isometric, there exists $\lambda\in \R$ such that, for all
$x,y\in TS$, $I(\nabla_xv,y)+I(x,\nabla_yv) = 2\lambda \II(x,y)$.
Thus:
\begin{eqnarray*}
\II((\db_{\III} w)X, X) & = & \II(B^{-1}\nabla_Xv + \Jb B^{-1} \nabla_{\Jb X}v,
X) \\
& = & I(\nabla_Xv, X) - I(\nabla_{\Jb X}v, \Jb X) \\
& = & \lambda \II(X,X) - \lambda \II(\Jb X, \Jb X) \\
& = & 0~, \\
\II((\db_{\III} w)X, \Jb X) & = & \II(B^{-1}\nabla_Xv + \Jb B^{-1} \nabla_{\Jb X}v,
\Jb X) \\
& = & I(\nabla_Xv, \Jb X) + I(\nabla_{\Jb X}v, X) \\
& = & 2\lambda \II(X, \Jb X) \\
& = & 0~,
\end{eqnarray*}
and we obtain that $\db_{\III}w = 0$.
\epv

\brk
Let $w$ be a solution of equation (\ref{eq:dbIII}) on $S$; there
exists a unique infinitesimal first-order deformation of $S$ whose
component tangent to $S$ is $Bw$. 
\erk

\bpv
Since $w$ is a solution of equation (\ref{eq:dbIII}), a simple
computation shows that the first-order deformation of the
metric on $S$ induced by $Bw$ is proportional to $\II$, i.e. it is of the form
$\lambda \II$, for some function $\lambda:S\rightarrow \R$. Thus, the vector
field $Bw - \lambda N$ is an isometric first-order deformation of $S$.
\epv

Of course, there are many solutions of equations (\ref{eq:dbIII}); all
"trivial" deformations of $S$ in $\R^3$ lead to one such solution, so
that there is at least a 6-dimensional space of solutions. The point is that
there is no other solution, a fact that we will prove by looking at a dual
formulation of the problem.

\pg{Variations of the shape operator}

The second rather simple way of stating the infinitesimal rigidity of
surfaces is in terms of first-order variations of the shape operator
$B$. We will see that this approach is, in a precise sense, adjoint of
the "tangent vector field" interpretation. We first recall the basics of the
infinitesimal rigidity of convex surfaces in $\R^3$ (or $H^3$).

\bprop
Let $u$ be an isometric first-order deformation of $S$, and let $\Bd$ be
the induced first-order variation of the shape operator $B$ of $S$. Then
$\Bd$ is a solution of:
\beq \label{eq:Bd}
\left\{ 
  \begin{array}{c}
    d^\nabla\Bd = 0 \\
    \tr(B^{-1}\Bd) = 0~.
  \end{array}
\right.
\eeq
\eprop

\bpv
The Gauss and Codazzi equations on $B$ are:
$$ \left\{
    \begin{array}{l}
      d^\nabla B=0 \\
      \det(B)=K~,
    \end{array}
\right. $$
where $K$ is the curvature of the induced metric. Since the deformation
$u$ is isometric, those equations remain true at first order, and the
proposition follows by linearizing them.
\epv

\bprop \label{pr:defos-isom}
Suppose that $S\subset \R^3$. 
Let $\Bd$ be a solution of (\ref{eq:Bd}). There exists an isometric
first-order deformation of $S$, unique up to
global infinitesimal isometries, whose induced deformation of $B$ is
$\Bd$. 
\eprop

This is a direct consequence of the so called "fundamental theorem of
the theory of surfaces", see e.g. \cite{S}, and is also valid for surfaces in
$H^3$. It is however interesting to
check how the proof can be done directly, since we will need the same
argument in a slightly more demanding context in the next sections; the proof
given here is valid in $\R^3$.

We start from our surface $S$ and from an embedding
$\phi:S\rightarrow \R^3$, and we call $N$ the unit normal vector to
$\phi(S)$, which we consider as a $\R^3$-valued function on $S$. The shape
operator $B$ of $S$ can be defined in this context by the equation:
\beq \label{eq:dN} dN = -(d\phi)\circ B~. \eeq

We consider a first-order deformation $\phid$
of $\phi$ which does not change the induced metric at first order, and
call $\Nd$ the first-order deformation of $N$. Then the deformation
$\phid$ acts locally like a rigid motion, so 
we can define at each point $x\in S$ a vector $Y\in \R^3$, which is such
that: 
\beq \label{eq:dphid} d\phid = Y\wedge d\phi~, \eeq
\beq \label{eq:Nd} \Nd = Y\wedge N~. \eeq
Here, and everywhere in the paper, ``$\wedge$'' is the vector product
in $\R^3$. 
Taking the first-order variation of (\ref{eq:dN}) shows that:
$$ d\Nd = -(d\phid)\circ B - (d\phi)\circ \Bd = - Y\wedge (d\phi\circ B)
- d\phi\circ \Bd~, $$
while taking the differential of (\ref{eq:Nd}) shows that:
$$ d\Nd = dY\wedge N + Y\wedge dN = dY\wedge N - Y\wedge d\phi\circ
B~. $$
Comparing those two equations, we have that:
\beq \label{eq:dYN} dY\wedge N = -d\phi\circ \Bd~. \eeq

Moreover, taking the differential of (\ref{eq:dphid}) yields that, for
any two vector fields $x$ and $y$ on $S$:
$$ 0 = d(d\phid)(x,y) = (dY\wedge d\phi)(x,y)~, $$
from which it follows that the image of $dY$ is always contained in the
image of $d\phi$, and therefore orthogonal to $N$. Thus:
$$ dY = - (dY\wedge N)\wedge N = (d\phi\circ \Bd)\wedge N~. $$
Taking the differential of this equations, we obtain that, for any two
vector fields $x,y$ on $S$:
\begin{eqnarray*}
  0 & = & dN(x)\wedge d\phi(\Bd y) - dN(y) \wedge d\phi(\Bd x) +
  N\wedge(x.d\phi(\Bd y) - y.d\phi(\Bd x)) - N\wedge 
  d\phi(\Bd[x,y]) \\
  & = & -d\phi(Bx)\wedge d\phi(\Bd y) + d\phi(By)\wedge d\phi(\Bd x) + \\
  & + & 
  N\wedge d\phi(\nabla_x (\Bd y)-\nabla_y(\Bd x)-\Bd[x,y]) + N\wedge
  ((x.d\phi)(\Bd y) -((y.d\phi)(\Bd x))~.
\end{eqnarray*}
But:
$$ N\wedge ((x.d\phi)(\Bd y) -(y.d\phi)(\Bd x)) = N\wedge ((\II(x,\Bd
y)-\II(y, \Bd x))N) = 0~, $$
because $N\wedge N=0$, so that:
$$ -d\phi(Bx)\wedge d\phi(\Bd y) + d\phi(By)\wedge d\phi(\Bd x) +
  N\wedge d\phi(d^{\nabla}\Bd)(x,y) = 0~, $$
which can be written as, using the complex structure $J$ of $I$, as:
$$ (-I(JBx,\Bd y) + I(JBy,\Bd x))N + N\wedge d\phi(d^{\nabla}\Bd)(x,y) =
0~. $$ 
Taking for $(x,y)$ an orthonormal basis of $TS$ shows that the component
orthogonal to $d\phi(TS)$ vanishes if and only if $\tr(B^{-1}\Bd)=0$,
while the component tangent to $TS$ is zero if and only if
$d^{\nabla}\Bd=0$. So we recover equation (\ref{eq:Bd}). 

\medskip

We can now prove Proposition \ref{pr:defos-isom}. Start from a solution
$\Bd$ of (\ref{eq:Bd}), and set:
$$ \alpha:=(d\phi\circ \Bd)\wedge N~, $$
so that $\alpha$ is a $\R^3$-valued 1-form on $S$. Then the computation
we've just done shows that $d\alpha=0$, precisely because $\Bd$ is a
solution of (\ref{eq:Bd}). So we can integrate $\alpha$ over $S$, and
find a function $Y:S\rightarrow \R^3$ such that $dY=\alpha$. 

We can then define another $\R^3$-valued 1-form $\beta$ on $S$ by:
$$ \beta := Y\wedge d\phi~. $$
Then:
$$ d\beta = dY\wedge d\phi = ((d\phi\circ \Bd)\wedge N)\wedge d\phi~. $$
Let $(x,y)$ be an orthonormal frame in which $\Bd$ is diagonal, we see
by computing $d\beta(x,y)$ that $d\beta=0$, so that we can define a
function $\phid:S\rightarrow \R^3$ such that $d\phid=\beta$. We then
have, by definition of $\beta$:
$$ d\phid = Y\wedge d\phi~, $$
so that $\phid$ is an isometric infinitesimal deformation of
$\phi(S)$. It also clearly follows from the computations made above that
the first-order variation of $B$ is indeed $\Bd$, and this proves
Proposition \ref{pr:defos-isom}.

\medskip

Now the key point, of course, is that there exists no such isometric
infinitesimal deformation, i.e. the smooth, strictly convex surfaces are
rigid. This result has been known for a long time, and the proof we will
give is essentially well-known and classical, see
e.g. \cite{herglotz,iie}. There are 
several other classical proofs, see e.g. \cite{vekua,S}, but this one is well
suited for our needs since it does not demand too much smoothness.

\blm \label{lm:rig-smooth}
Suppose that $S\subset \R^3$.
There is no non-trivial solution $\Bd$ of (\ref{eq:Bd}) on $S$.
\elm

\bpv
We consider the function $f:=r^2/2$ on $\R^3$, where $r$ is the distance
to $0$. Then:
$$ \hess(f) = g_0~, $$
where $g_0$ is the canonical flat metric on $\R^3$.

Let $\hessb(f)$ be the Hessian of the restriction of $f$ to $S$. 
Then:
$$ \hessb(f) = I - df(N) \II~, $$
where $N$ is the unit exterior normal of $S$. Thus, calling $\fd$ the
first-order variation of the restriction of $f$ to $S$:
$$ \hessb(\fd) = - df(N) \IId - \dot{df(N)} \II~. $$

Now let $\omega(X):=d\fd(J\Bd X)$. Then:
\begin{eqnarray*}
  d\omega(X,Y) & = & X.d\fd(J\Bd Y) - Y.d\fd(J\Bd X) - d\fd(J\Bd [X,Y]) \\
  & = & \hessb(\fd)(X, J\Bd Y) - \hessb(\fd)(Y, J\Bd X) + d\fd(J\nabla_X(\Bd Y) -
  J\nabla_Y(\Bd X) - J\Bd[X,Y]) \\
  & = & \hessb(\fd)(X, J\Bd Y) - \hessb(\fd)(Y, J\Bd X) +
  d\fd(J(d^\nabla\Bd)(X,Y)) \\ 
  & = & -df(N)(\IId(X, J\Bd Y) - \IId(Y, J\Bd X)) -
  \dot{df(N)} (\II(X, J\Bd Y) - \II(Y, J\Bd X))~.
\end{eqnarray*}
The second term
vanishes because $\tr(B^{-1}\Bd)=0$, 
while it is simple to check --- for instance by taking for $X$
and $Y$ eigenvectors of $\Bd$ --- that:
$$ d\omega(X,Y) = 2df(N)\det(\Bd)da(X,Y)~. $$
Since $\tr(B^{-1}\Bd)=0$, $\det(B^{-1}\Bd)\leq 0$ and thus $\det(\Bd)\leq 0$ 
everywhere. Since the integral of
$d\omega$ over $S$ is zero, $\det(\Bd)=0$ everywhere, so that $\Bd=0$. 
\epv

\pg{A duality between the two formulations}

An important element used in the next sections to prove a rigidity result
is that the two formulations of the rigidity stated above --- in terms of
a vector field tangent to $S$, or in terms of a bundle morphism of $TS$
describing the first-order variation of $B$ --- are adjoint to each
other. This follows from the next two propositions, which are valid in either
$\R^3$ or $H^3$.

\bprop \label{pr:adjoint}
The adjoint of $\Jb \db_{\III}$ for $\II$ is the operator defined by:
$$
\begin{array}{rrcl}
(\Jb\db_{\III})^*: & \ots & \rightarrow & TS \\
& h & \mapsto & - d^\nabla h(X, \Jb X) = - (\nabla_Xh)(\Jb X) + (\nabla_{\Jb X}h)(X)~,
\end{array}
$$
where $X$ is a unit vector field on $S$. If $\Omega$ is an open subset
with smooth boundary of $S$ and $h$ is a smooth section of $\ots$, then:
$$ \int_\Omega \langle \Jb \db_{\III} v, h\rangle_{\II} da_{\II} =
\int_\Omega \II(v, -d^\nabla h) + \int_{\dr \Omega} \II(v, h(ds))~. $$
Similary, the adjoint of $\Jb \db_I$ for $\II$ is the operator $-d^{\nablat}$
acting on sections of $\ots$, and:
$$ \int_\Omega \langle \Jb \db_{I} v, h\rangle_{\II} da_{\II} =
\int_\Omega \II(v, -d^{\nablat}h) + \int_{\dr \Omega} \II(v, h(ds))~. $$
\eprop

It is implicit in this proposition that 
the expression of the adjoint operators is actually
independent of the choice of the unit vector $X$.

\bpv
Let $h$ be a smooth section of $\ots$. Then, using Proposition
\ref{pr:mixte}: 
\begin{eqnarray*}
  \int_S \II(\Jb\db_{\III}v, h) da_{\II} 
    & = & \int_S \II(\db_{\III}v, -\Jb h)da_{\II} \\
    & = & \int_S \II((\db_{\III}v)(X), h(\Jb X)) da_{\II} \\
    & = & \int_S \II(\nablat_Xv + \Jb\nablat_{\Jb X}v, h(\Jb X)) da_{\II} \\
    & = & \int_S \II(\nablat_Xv, h(\Jb X)) - \II(\nablat_{\Jb X}v, h(X))
    da_{\II} \\ 
    & = & \int_S X.\II(v, h(\Jb X)) - \II(v, \nabla_X(h(\Jb X)))
    - (\Jb X).\II(v,
    h(X)) + \II(v, \nabla_{\Jb X}(h(X)) da_{\II} \\
    & = & \int_S X.\II(v, h(\Jb X)) - (\Jb X).\II(v, h(X)) - \\
    & & - \II(v, h([X, \Jb X])) - \II(v, (\nabla_Xh)(\Jb X)-(\nabla_{\Jb X}h)(X)))
    da_{\II} \\  
    & = & \int_S d\omega(X, \Jb X) - \II(v, (d^\nabla h)(X, \Jb X)) da_{\II}~,
\end{eqnarray*}
with $\omega(X):=\II(v, h(X))$. The result follows for
$(\Jb\db_{\III})^*$. The formula just obtained also shows the result
concerning the integration over $\Omega$. 
A similar computation shows the corresponding results
for $(\Jb\db_I)^*$.
\epv

The point is that those adjoint operators are precisely those which
describe the first-order deformations in terms of the variation of $B$. The
next two propositions are valid if $S\subset \R^3$ or $S\subset H^3$.

\bprop \label{pr:equivalence}
Let $\Bd$ be a bundle morphism from $TS$ to itself, which is self-adjoint for
$I$. $\Bd$ is a solution
of equation (\ref{eq:Bd}) if and only if $h:=B^{-1}\Bd$ is a section of
$\ots$ such that $d^{\nablat}h=0$.
\eprop

\bpv
Clearly, $h$ is self-adjoint for $\II$: if $X,Y$ are two vectors at a
point $s\in S$, then 
$$ \II(h(X),Y) = \II(B^{-1}\Bd X,Y) = I(\Bd X,Y) = I(X, \Bd Y) = \II(X,
B^{-1}\Bd Y) = \II(X, h(Y))~. $$
Moreover $\Bd$ satisfies the second equation of (\ref{eq:Bd}),
i.e. $\tr(B^{-1}\Bd)=0$, if and only if $\tr(h)=0$. But it is a simple
matter of linear algebra that those two conditions are satisfied if and
only if $h$ is a section of $\ots$, i.e. if and only if, for each $X\in
TS$, $h(\Jb X)=-\Jb h(X)$. 

Similarly, the condition that $d^{\nablat}h=0$ can be written as:
\begin{eqnarray*}
  0 & = & (d^{\nablat}h)(X, \Jb X) \\
  & = & \nablat_X(h(\Jb X)) -  \nablat_{\Jb X}(h(X)) - h([X, \Jb X]) \\
  & = & B^{-1}(\nabla_X (\Bd\Jb X) - \nabla_{\Jb X}(\Bd X) - \Bd[X, \Jb X]) \\
  & = & B^{-1}(d^\nabla\Bd)(X, \Jb X)
\end{eqnarray*}
so that it is equivalent to the first equation in (\ref{eq:Bd}). 
\epv

The same proof can be used to show a "dual" statement:

\bprop \label{pr:equivalence*}
Let $\Ad$ be a bundle morphism from $TS$ to itself. $\Ad$ is a solution
of:
\beq \label{eq:Bd*}
\left\{ 
  \begin{array}{c}
    d^{\nablat}\Ad = 0 \\
    \tr(A^{-1}\Ad) = 0
  \end{array}
\right. \eeq
if and only if $h:=A^{-1}\Ad$ is a section of
$\ots$ such that $d^\nabla h=0$.
\eprop

\pg{Rigidity of a surface and its dual}

Let's sum up the remarks we have just made. Each
infinitesimal isometric deformation of $S$ corresponds to a solution of
$\db_{\III}\wb=0$; the adjoint of this equation is $d^\nabla h=0$, and the
solutions of this equation correspond to solutions $\Ad$ of equations
(\ref{eq:Bd*}). 

Suppose now that $S\subset H^3$, and let $S^*\subset S^3_1$ be
the dual surface. Then the induced metric on $S^*$ is the third fundamental
form of $S$ (and conversely) while the second fundamental forms of $S$ and
$S^*$ are the same. Therefore, solutions of equation (\ref{eq:Bd*}), and
solutions of $d^\nabla h=0$ with $h$ a section of $\ots$, describe
isometric first-order deformations of $S^*$.
Conversely, non-trivial infinitesimal
isometric deformations of $S$ also correspond to non-zero solutions
$\Bd$ of (\ref{eq:Bd}), and thus to solutions $d^{\nablat}h=0$, and the
adjoint equation is $\db_I\wb=0$, which also describes infinitesimal isometric
deformations of $S^*$.

\pg{Duality and the Minkowski problem}

The reader has probably by now realized that the duality between the
induced metric and the third fundamental form plays an important role in
the constructions made in this section --- and it will also be central
in the next sections. In the setting of the hyperbolic
3-space, the geometric meaning of the duality is clear: the
third fundamental form of a surface $S\subset H^3$ is the induced metric
on the dual 
surface in the de Sitter space. Then solutions of $\db_{\III}\wb=0$, and
of $d^{\nablat}h=0$ for $h\in  \ots$, correspond to infinitesimal
isometric deformations of the surface $S$, while solutions of
$\db_I\wb^*=0$ and of $d^\nabla h^*=0$ correspond to isometric deformations of
the dual surface $S^*$. 

Although it will not be used in this paper, it is interesting to point
out that this duality retains a geometric meaning in $\R^3$, although
the dual problem is not one of isometric deformations of surfaces. Indeed,
solutions of $d^\nabla h=0$ correspond to variations $\Ad$ of $A:=B^{-1}$
which are solutions of (\ref{eq:Bd*}), and we leave it to the interested
reader to check that the non-existence of those in turns corresponds to
the rigidity in the {\bf Minkowski problem}, which deals with the
existence of convex surfaces with a given push-forward, by the Gauss
map, of its curvature measure.


\section{First-order deformations} \label{se:2b}

\pg{The developing map}

We have seen in section 2 some details on the behavior of $\Sb$ and
$\Sb^*$ near the limit set $\Lambda$. We will now study the
boundary behavior of the vector fields on $\Sb$ and $\Sb^*$
corresponding to the first-order variations $\gd$ of the hyperbolic
metric $g$ on $M$. We will use this in the next two sections to prove the
infinitesimal rigidity of hyperbolic manifolds with convex boundary.

The first point is to remark that one can associate to $\gd$ a closed 
one-form
$\omega$ on $M$ with value in $\SL(2,\C)$. This can be done along the ideas of
Calabi 
\cite{Calabi} and Weil \cite{Weil}, but we will follow here an elementary
approach. $\omega$ will not be canonically defined, but it will not be
necessary. 

To define $\omega$, we first recall the definition of the developing map
of a hyperbolic metric on $M$. 

\bdf 
Let $g$ be a hyperbolic metric with convex boundary on $M$. The
universal cover $\Mt$ of $M$ is a simply connected hyperbolic manifold
with convex boundary, which thus admits a unique isometric embedding
into $H^3$. This defines a map $\dev_g:\Mt\rightarrow H^3$, the {\bf
  developing map} of $g$.
\edf

$\dev_g$ is defined modulo the global isometries of $H^3$. Given a first-order
deformation $\gd$ of $g$, it determines a first-order deformation $\dev'$ of
$\dev$. $\dev'$ determines a vector field $u$ on $H^3$. $\dev'$ is defined
modulo the global Killing vector fields on $\Omega$; from here on we suppose
that $\dev'$ is normalized so that $u$ vanishes at $x_0$. 

\paragraph{The natural $\SL_2(\C)$ bundle over $H^3$}

Questions of infinitesimal deformations of hyperbolic 3-manifolds are
classically associated to a natural $\SL_2(\C)$-bundle over $H^3$, and over
hyperbolic manifolds. We recall some basic constructions concerning this
aspect of the deformations. We consider here a convex subset $\Omega\subset
H^3$. 
We will call $E$ the trivial $\SL_2(\C)$ bundle over $\Omega$. For each $x\in
\Omega$ and each $\kappa\in E_x$, $\kappa$ is an element of $\SL_2(\C)$,
in other terms a 
Killing vector field over $H^3$. 

Given $x\in \Omega$ and $U\in T_x\Omega$, there is a unique Killing field $v$
over 
$\Omega$ which is a "pure translation" of vector $U$ at $x$; 
it is defined uniquely by
the fact that $v(x)=U$ and that $\nabla v=0$ at $x$. There is also a unique
Killing field $v'$ which is a "pure rotation of vector $U$" at $x$: 
it is uniquely defined 
by the fact that $v'(0)=0$ and that $\nabla_wv'=U\wedge w$ for each $w\in
T_x\Omega$. A standard fact is that $\SL_2(\C)$ is the direct sum of those two
kinds of Killing fields, i.e. each Killing field on $\Omega$ has a unique
decomposition as the sum of a pure translation at $x$ and a pure rotation at
$x$. This yields a decomposition: 
$$ E_x = T_x\Omega\oplus T_x\Omega~. $$
This decomposition depends on $x$, since a Killing field which is a pure
translation (or rotation) at a point $x$ is not a pure translation (or
rotation) at another point $y\neq x$. 

So, given a point $x\in \Omega$, we have a natural identification of the vector
space of Killing fields over $H^3$ with $T_x\Omega\oplus T_x\Omega$. We will
often denote a Killing field $\kappa$ by a pair $(\tau_x, \sigma_x)$ ---
depending on 
the choice of a point $x\in \Omega$. Then $\tau_x\in T_x\Omega$ will be the
translation vector 
of the component of $\kappa$ which is a pure translation at $x$ --- so that
$\tau_x=\kappa(x)$ --- while $\sigma_x\in T_x\Omega$ will be the rotation
vector 
of the component of $\kappa$ which is a pure rotation at $x$. 

The bundle $E$ comes equipped with a natural flat connection, which we will
call $D$. Let $\kappa$ be a section of $E$, so that, for each $x\in \Omega$,
$\kappa(x)$ is an element of $\SL_2(\C)$, i.e. a Killing field over $H^3$. 
Then, for each $x\in H^3$, $\kappa(x)=(\tau_x, \sigma_x)\in T_x\Omega\oplus
T_x\Omega$. The connection $D$ is defined in those notations by:
\beq \label{eq:D} 
D_w(\tau_x, \sigma_x) = (\nabla_w\tau_x + w\wedge \sigma_x,
\nabla_w\sigma_x - w\wedge \tau_x)~. \eeq
The key point is the following Property of $D$.

\begin{pty}
Let $\kappa$ be a section of $E$. $D\kappa=0$ if and only if $\kappa$
corresponds, at each point $x\in H^3$, to the same Killing field over $H^3$. 
\end{pty}

From there it follows immediately that $D$ is a flat connection. However $D$
does not preserve the identification $E_x=T_x\Omega\oplus T_x\Omega$. 

Now, given $\dev'$, we wish to define a section $\kappa$ of $E$. We will do it
in terms of the translation and rotation components of $\kappa$ at each
point. 
For each $x\in \Omega$, we call $\tau_x$ the infinitesimal pure
hyperbolic translation at
$x$ of translation vector $u=\dev'$. So $\tau_x$ is a
Killing 
vector field on $H^3$, or, in other terms, an element of $\SL_2(\C)$. Therefore
$\tau$ determines a section over $\Omega$ of the natural $\SL_2(\C)$ bundle.

Consider the vector field $w=u-\tau_x$ on $H^3$. By construction it
vanishes at $x$, so we can consider its differential $dw:T_x\Omega\rightarrow
T_x\Omega$. $dw$ decomposes uniquely as the sum $dw=a+s$, where
$a, s:T_x\Omega\rightarrow T_x\Omega$ are respectively adjoint and
skew-adjoint for 
the hyperbolic metric. There is a unique infinitesimal hyperbolic rotation of
axis containing $x$ and differential $s$, we call it $\sigma_x$ --- and
we also call 
$\sigma_x$ the rotation vector of that infinitesimal isometry. 
Again, $\sigma_x$ is a
Killing field on $H^3$, or an element of $\SL_2(\C)$, and $\sigma$ defines
another section over $\Omega$ of the natural bundle with fiber $\SL_2(\C)$. 
By choosing the right normalization for $\dev'$ above, we decide that
$\sigma_{x_0}=0$. 

In this way, we obtain from $\dev'$ a section $\kappa$ of $E$, which we call
$(\tau, 
\sigma)$ with the notations explained above. At each point $x\in \Omega$,
$\tau_x, \sigma_x\in T_x\Omega$, and $\kappa_x$ is the Killing field over
$H^3$ which is the sum of a pure translation of vector $\tau_x$ and a pure
rotation of rotation vector $\sigma_x$ at $x$.

We can use $\kappa$ and $D$ to define a closed 
1-form $\omega$ with values in the
$\SL_2(\C)$-bundle over $\Omega$, as the exterior differential of $\kappa$:
$$ \forall v\in T\Omega, ~ \omega(v):=(d^D\kappa)(v) = D_v\kappa~. $$
$\kappa$ was defined geometrically, and its definition shows that, if
one adds a Killing field $U$ to $\dev'$, $\kappa$ is replaced by
$\kappa+U$, and $\omega$ does not change. Thus the equivariance of
$\dev'$ implies that 
$\omega$ is invariant under the action of $\pi_1M$, so it is defined over
$M$.

Since $M$ is compact, $\omega$ is bounded: there exists a fixed constant
$C_\omega$ such that, 
for any $x\in \Omega$ and any unit vector $v\in T_x\Omega$, both the
translation and the 
rotation components, at $x$, of $\omega(v)$ are bounded by $C_\omega$. By
"bounded" we mean here that the translation vector of $\tau_x$ at $x$, and the
rotation speed of $\sigma_x$ at $x$, are bounded. 

Using the identification of $E$ with $T\Omega\oplus T\Omega$, $\omega$ can be
written as $\omega=(\omega_\tau, \omega_\sigma)$, where both $\omega_\tau$ and
$\omega_\sigma$ are 1-forms with values in the tangent bundle $T\Omega$. Since
$\omega$ is bounded, both $\omega_\tau$ and $\omega_\sigma$ are bounded. 
Moreover, it is a direct consequence of the definition of $\omega$ that, for
each $w\in T\Omega$: 
\beq \label{eq:nabla} 
\nabla_w\tau = \omega_\tau(w) - w\wedge \sigma, ~ \nabla_w\sigma =
\omega_\sigma(w) + w\wedge \tau~. \eeq

We now sum up some relations between the deformation vector field $u$ and the
section $\kappa$ which we have defined from $\dev'$. 

\bprop
\begin{itemize}
\item Let $u$ be the vector field defined on $S$ by the deformation
  $\gd$. Then $u(x)=\tau_x$. If
  $\gd$ does not change the induced metric 
  on $\dr M$, then $u$ is an infinitesimal isometric deformation of $S$.
\item $\kappa$ defines a vector field $u^*$ on $S^*$. If
  $\gd$ does not change the third fundamental form of $\dr M$, then $u^*$ is
  an infinitesimal isometric deformation of $S^*$. 
\end{itemize}
\eprop

\bpv
$u$ is the vector field which determines the first-order deformation of $S$
under the deformation $\gd$ of $g$, so it is clear that $u$ is an
infinitesimal isometric deformation of $S$ if and only if the induced metric
on $\dr M$ does not change. The definition of $u^*$ follows directly from
the first-order deformation of $S^*$ through the duality with $S$. 
The same argument as for $u$ applies to the first-order
variation of $\III$ under the deformation $u^*$, through the variation
of the induced metric on $S^*$ and the Pogorelov map for $S^3_{1,+}$. 
\epv

\pg{Deformations of Euclidean metrics}

We have seen above how the first-order deformations of a hyperbolic manifold
can be described in terms of a 1-form with values in the canonical
$\SL_2(\C)$-bundle $E$, and how $E$ can be identified with $T\Omega\oplus
T\Omega$ at each point. We will now see how the same formalism applies ---
with minor differences --- to deformations of a Euclidean manifold. The point
will then be that the Pogorelov map allows one to move from the hyperbolic to
the Euclidean formalism.

The space of Killing fields on $\R^3$ is a 6-dimensional vector space, which
we will call $\Eb_0$. We will consider here a domain $\Omegab\subset \R^3$,
which will later be the image in the projective model of $\Omega$. We call
$\Eb$ the trivial $\Eb_0$-bundle over $\Omegab$. 
Given $\xb\in \Omegab$, there is a canonical identification
$\Eb_x=T_x\Omegab\oplus T_x\Omegab$, with the first factor corresponding to
the pure translation component of the Killing field at $\xb$, and the second
factor corresponding to the rotation component at $\xb$. In more explicit
terms, to a Killing field $\kappab$, we associate $(\taub,\sigmab)\in
T_{\xb}\Omegab\oplus T_{\xb}\Omegab$, with:
$$ \taub := \kappab(x)~, $$
and, if $\nablab$ is the flat connection on $\Omegab$, 
for each $w\in T_{\xb}\Omegab$:
$$ \nablab_w\kappab = \sigmab\wedge w~. $$
 
As in the hyperbolic case, there is a flat canonical connection $\Db$ on $\Eb$,
defined as follows. Given a section of $\Eb$, defined by sections $\tau$ and
$\sigma$ of $T\Omegab$ using the identification above, we have for all $w\in
T\Omegab$: 
\beq \label{eq:Db}
\Db_w(\taub,\sigmab) = (\nablab_w \taub + w\wedge \sigmab,
\nablab_w\sigmab)~. \eeq
It has the same basic property as in the hyperbolic setting.

\begin{pty}
Let $\kappab$ be a section of $\Eb$. $\Db\kappab=0$ if and only if $\kappab$
corresponds, at each point $\xb\in \Omegab$, to the same Killing field over
$\Omegab$.  
\end{pty}

\paragraph{The Pogorelov map and infinitesimal isometries}

An important point is that the Pogorelov map gives a way to translate things
from the hyperbolic to the Euclidean formalism. We now consider the domain
$\Omega\subset H^3$ corresponding to the universal cover of $M$, and call
$\Omegab$ its image in the projective model. 

\bdf 
Let $x\in \Omega$ and $\kappa\in E_x$. $\kappa$ is a Killing field on
$\Omegab$, 
and the image of this vector field by the Pogorelov map is a Killing field on
$\Omegab$, which we call $\kappab$. We call $\Psi_H:E\rightarrow \Eb$ the
bundle morphism sending $\kappa\in E_x$ to $\kappab$, considered as an element
of $\Eb_{\phi(x)}$.  
\edf

Of course $\Psi_H$ is just the Pogorelov map, we use a different notation from
$\Phi_H$ since it acts on sections of $E$ instead of vector fields. 
The following statement is then a direct consequence of the properties of the
Pogorelov map. 

\bprop \label{pr:commute}
Let $\kappa$ be a section of $E$. Then:
$$ \Psi_H(D\kappa) = \Db(\Psi_H(\kappa))~. $$
\eprop

\bpv
The Pogorelov map $\Phi_H$ sends hyperbolic Killing fields to Euclidean
Killing fields. Therefore, $\Psi_H$ sends flat sections of $E$ for $D$ to flat
sections of $\Eb$ for $\Db$. It follows that it respects the connections $D$
and $\Db$. 
\epv

It is now time to express explicitely the action of $\Psi_H$ on different
components of a section of $E$. From here on, we call $R$ the unit radial
vector on $\Omega$ --- in 
the direction opposite to $x_0$ --- and $\Rb$ the unit radial vector on
$\Omegab$ --- also in the direction opposite to $0$. 

\blm \label{lm:PsiH}
Let $x\in H^3\setminus \{ x_0\}$,
and let $v\in T_xH^3$ be a unit lateral
vector. Let $\xb:=\phi(x)$,  with $r:=d(\xb,0)$, and let $\vb$ be the
unit lateral vector in the direction of $\Phi_H(v)$. Then:
$$ \Psi_H(R, 0) = (\Rb, 0), ~ \Psi_H(0, R) = (0, \Rb)~, $$
$$ \Psi_H(v, 0) = \left (\left(\sqrt{1-r^2}\right)\vb, -\frac{r}{\sqrt{1-r^2}}
  \Rb\wedge 
  \vb\right), ~ \Psi_H(0,v) = \left( 0, \frac{1}{\sqrt{1-r^2}} \vb\right)~. $$
\elm

\bpv
$(R,0)$ corresponds to an infinitesimal translation at unit speed with axis the
geodesic through $x_0$ and $x$. Therefore $\Psi_H(R,0)$ leaves globally
invariant the segment through $0$ and $\xb$. Checking its behavior at $0$
shows that it is a translation at unit speed along this segment. 

$(0,R)$ is an infinitesimal rotation at unit speed along the geodesic going
through $x_0$ and $x$, so it leaves this geodesic pointwise invariant. So the
same is true for $\Psi_H(0,R)$ and the segment through $0$ and $\xb$, and
again considering the action at $0$ shows that it is the infinitesimal
rotation at unit speed and axis $(0,\xb)$, i.e. $(0, \Rb)$. 

$(v, 0)$ leaves globally invariant the totally geodesic plane containing the
geodesic through $x_0$ and $x$, and tangent to $v$. So $\Psi_H(v,0)$ also
leaves invariant the plane containing $(0,\xb)$ and tangent to
$\vb$. Moreover, $(v,0)$ is a Killing field which is orthogonal to the
geodesic through $x_0$ and $x$, so $\Psi_H(v,0)$ is orthogonal to $(0, \xb)$. 
Therefore $\Psi_H(v,0)$ is the sum of an infinitesimal translation of vector
proportional to $\vb$ and of a rotation of axis parallel to $\Rb\wedge \vb$. 

To identify the two coefficients, we compute its action at $0$ and at $\xb$. At
$0$, it acts by a vector $\cosh(\rho)\vb$, where $\rho=d(x_0,x)$. At $\xb$,
the definition of the Pogorelov map shows that it acts by a vector
$(1/\cosh(\rho)) \vb$. So:
\begin{eqnarray*}
\Psi_H(v,0) & = & \left(\frac{1}{\cosh(\rho)}\vb,
  \frac{1}r\left(\frac{1}{\cosh(\rho)} -\cosh(\rho)
\right)\Rb \wedge \vb\right) \\
  & = & \left(\left(\frac{1}{\cosh(\rho)}\right)\vb,
  \left(-\frac{\sinh^2(\rho)}{\tanh(\rho)\cosh(\rho)} 
\right)\Rb \wedge \vb\right) \\
  & = & \left(\left(\frac{1}{\cosh(\rho)}\right)\vb,
  -\sinh(\rho)\Rb \wedge \vb\right) \\
  & = & \left(\left(\sqrt{1-r^2}\right) \vb, -\frac{r}{\sqrt{1-r^2}}\Rb \wedge
  \vb\right)~.
\end{eqnarray*}

The argument is similar for $\Psi_H(0, v)$; since $(0,v)$ is a rotation with
axis the geodesic containing $x$ and tangent to $v$, $\Psi_H(0,v)$ leaves
globally invariant the segment starting from $\xb$ in the direction of $\vb$,
so it is a rotation having this segment as its axis. To find the rotation
angle, remark that the action of $\Psi(0,v)$ at $0$ is by a vector of norm
$\sinh(\rho)$, so that the rotation speed of $\Psi(0,v)$ is $\sinh(\rho)/r$,
which is equal to $1/\sqrt{1-r^2}$. 
\epv

\pg{Upper bounds on $\omegab$}

We will need later on to estimate the deformation of the Euclidean manifold
$\Omegab$. This deformation is described by the section
$\kappab=\Psi_H(\kappa)$ of $\Eb$. By Proposition \ref{pr:commute}, $\Db
\kappab=\Psi_H(D\kappa)$, so that, if we set $\omegab:=d^{\Db}\kappab$,
we have: $\omegab=\Psi_H(\omega)$. 

From here on, 
given a point $x\in \Omega, x\neq x_0$, and a
vector $u\in T_x\Omega$, we call $u^r$ and $u^\perp$, respectively, the radial
and the lateral components of $u$. We also use the same notation in
$\Omegab$. Given a vector $w\in T\Omegab$, we will also call $\omegab_\tau(w)$
and $\omegab_\sigma(w)$, respectively, the components of $\omegab(w)$ in the
first and second factors of $\Eb=T\Omegab\oplus T\Omegab$. Thus, for instance,
$\omegab^r_\sigma(w)$ is the radial component of the rotation vector of
$\omegab(w)$, etc. 

\blm \label{lm:omegab}
Let $\xb\in \Omegab\setminus \{ 0\}$ with $r=d(0,\xb)>0$, 
and let $\vb\in T_x\Omegab$ be a lateral unit vector. Then:
$$ 
\begin{array}{rclrcl}
\|\omegab^r_\tau(\Rb)\| & \leq & \frac{C_\omega}{1-r},~ &
\|\omegab_\tau^r(\vb)\| & \leq & 
\frac{C_\omega}{\sqrt{1-r}} \\
\|\omegab^r_\sigma(\Rb)\|& \leq & \frac{C_\omega}{1-r},~ &
\|\omegab_\sigma^r(\vb)\|& \leq & 
\frac{C_\omega}{\sqrt{1-r}} \\
\|\omegab^\perp_\tau(\Rb)\|& \leq & \frac{C_\omega}{\sqrt{1-r}},~ &
\|\omegab_\tau^\perp(\vb)\| & \leq & C_\omega \\
\|\omegab^\perp_\sigma(\Rb)\| & \leq & \frac{2C_\omega}{\sqrt{1-r}^3}~, &
\|\omegab_\sigma^\perp(\vb)\| & \leq & \frac{2C_\omega}{1-r}~.
\end{array} $$
\elm

\bpv
We already know that $\omega_\tau$ and $\omega_\sigma$ are both bounded by
$C_\omega$. The upper bounds are consequences of two things:
\begin{itemize}
\item the behavior of the projective map $\phi$, i.e. the fact that
  $\Rb=\cosh^2(\rho)d\phi(R)$ and $\vb=\cosh(\rho)d\phi(v)$, where $v$ is a
  unit vector.
\item the properties of $\Psi_H$, as described in Lemma \ref{lm:PsiH}.
\end{itemize}
In the first and second line, we consider the radial components, for which
there is no scaling by Lemma \ref{lm:PsiH}, so that the only scaling which
comes in is the one coming from $\phi$, which adds a factor
$\cosh^2(\rho)=1/(1-r^2)$ for the radial vector $\Rb$ and a factor
$\cosh(\rho)=1/\sqrt{1-r^2}$ for the lateral vector $\vb$. The estimates
follow using the fact that $1-r^2=(1-r)(1+r)\geq (1-r)$.

On the third line, for $\omegab^\perp_\tau(\Rb)$, there is a factor
$\cosh^2(\rho)$ corresponding to the scaling of $\Rb$, and a factor
$\sqrt{1-r^2}$ coming from Lemma \ref{lm:PsiH}. For
$\omegab^\perp_\tau(\vb)$, there is a factor 
$\cosh(\rho)$ corresponding to the scaling of $\vb$, which compensates a
factor $\sqrt{1-r^2}$ coming from Lemma \ref{lm:PsiH}.

Finally, on the fourth line, one has to be a little more careful, because
Lemma \ref{lm:PsiH} shows that the lateral rotation component of $\omegab$
comes from both the lateral rotation and the lateral translation component of
$\omega$. The term $\omegab^\perp_\sigma(\Rb)$ comes from:
\begin{itemize}
\item the term $\omega^\perp_\sigma(R)$, with a coefficient $\cosh^2(\rho)$
  coming from the scaling of $\Rb$, and a factor $1/\sqrt{1-r^2}$ from Lemma
  \ref{lm:PsiH}. 
\item the term $\omega^\perp_\tau(R)$, with a coefficient $\cosh^2(\rho)$ for
  the scaling of $\Rb$, and a coefficient $r/\sqrt{1-r^2}$ from Lemma
  \ref{lm:PsiH}. 
\end{itemize}
Similarly, the term $\omegab^\perp_\sigma(\vb)$ comes from:
\begin{itemize}
\item the term $\omega^\perp_\sigma(v)$, with a coefficient $\cosh(\rho)$
  coming from the scaling of $\Rb$, and a factor $1/\sqrt{1-r^2}$ from Lemma
  \ref{lm:PsiH}. 
\item the term $\omega^\perp_\tau(v)$, with a coefficient $\cosh(\rho)$ for
  the scaling of $\Rb$, and a coefficient $r/\sqrt{1-r^2}$ from Lemma
  \ref{lm:PsiH}. 
\end{itemize}
\epv

In the next corollary, and in the rest of the paper, we call $C$ a
"generic" constant, whose precise value can change from line to line. In
each line, the value of $C$ might depend on $M$ and on $g$, but not on
more local data. 

\bcr \label{cr:omegab}
Let $\xb\in \Omegab\setminus \{ 0\}$ with $r=d(0,\xb)>0$. Then, at
$\xb$:
$$ \|\taub^r(x)\| \leq C |\log(1-r)|~, ~~ \|\taub^\perp(x)\|\leq
C\sqrt{1-r}~, $$
$$ \|\sigmab^r(x)\| \leq C |\log(1-r)|~, ~~ \|\sigmab^\perp(x)\|\leq
\frac{C}{\sqrt{1-r}}~. $$
\ecr

\bpv 
This is a direct consequence of the previous Lemma, obtained by
integrating the inequalities of the first column on the segment $[0,\xb]$.
\epv

It is interesting to note that the estimate on $\sigmab^\perp$ is better
than what would have been obtained by a straightforward estimate in
$H^3$. This can be explained using a more subtle estimate in $H^3$, very
similar to the one used in the proof of Lemma \ref{lm:holder*} below.


\section{H{\"o}lder bounds}

Before giving, in the next section, the key infinitesimal rigidity Lemma, we
need some precise estimates on the behavior, in the Euclidean ball, of the
deformation vector fields on $\Sb$ and $\Sb^*$. They will follow from a
careful analysis of the Killing vector fields induced by the first-order
deformation $\gd$ of $g$. 

Recall that $u$ is a vector field induced on $S$ by the first-order
deformation $\gd$ of $g$. We call $\ub$ its image by the Pogorelov
transformation $\Phi_H$, so that $\ub$ is an isometric first-order deformation of
the surface $\Sb$ which is the image of $S$ in the projective model. We
also call $\vb$ the orthogonal projection of $\ub$ on $\Sb$. 

\paragraph{Bounds on $\nablab \vb$}

Recall that we call $\delta$ the Euclidean distance, along $\Sb$, to the limit
set $\Lambda$. By Lemma \ref{lm:slope}, $\sqrt{1-r}/c_1\leq \delta\leq
c_1\sqrt{1-r}$ near $\Lambda$, where $c_1>0$ is
some constant. The first lemma we need is an upper bound, near the limit
set $\Lambda$, on the norm of $\nablab \vb$. 

\blm \label{lm:H1-log}
There exists a constant $C>0$ such that, on $\Sb$:
$$ \|\nablab \vb\|\leq C|\log(\delta)|~. $$
\elm

\bpv
We denote by $\Pi_{\xb}$ the orthogonal projection on $\Sb$, seen as a
morphism from 
$T_{\xb}\R^3$ to $T_{\xb}\Sb$, for $\xb\in \Sb$. Let $\xb\in \Sb$, and let
$\wb\in T_{\xb}\Sb$ be a vector of norm at most $1$. Then $\wb=\wb^\perp +
\lambda \Rb$, where $\wb^\perp$ is the lateral component of $\wb$;
$\|\wb^\perp\|\leq 1$ and $|\lambda|\leq C\delta$ by
Lemma \ref{lm:slope}.  

Now $\vb=\Pi_{\xb}(\taub(\xb))$ by definition, so that:
$$ \nablab_{\wb}\vb = (\nablab_{\wb}\Pi)\vb + \Pi(\nablab_{\wb}\taub)~, $$
and we have by (\ref{eq:Db}):
\beq \label{eq:princ}
\nablab_{\wb}\vb = (\nablab_{\wb}\Pi)\vb - \Pi({\wb}\wedge\sigmab) +
\Pi(\omegab_\tau({\wb}))~. 
\eeq

By Corollary \ref{cr:omegab}, we have:
$$ \|\taub(\xb)\|\leq C|\log(1-r)|~. $$
In addition, $\|\nablab_{\wb}\Pi\|\leq C$ by Lemma \ref{lm:curvature}
(recall that, in this proof, $C$ is a generic constant whose value can
change from line to line). Putting both together and using that $\delta$
behaves as $\sqrt{1-r}$, we find that:
\beq \label{eq:i1} 
\|(\nablab_{\wb}\Pi)\vb\|\leq C|\log(\delta)|~. \eeq

For the second term in (\ref{eq:princ}), note that, at $\xb$,
$\sigmab=\sigmab^r+\sigmab^\perp$, where, again  by Corollary \ref{cr:omegab}:
$$ \|\sigmab^r\|\leq C|\log(1-r)|\leq C|\log(\delta)|,~ 
\|\sigmab^\perp \|\leq \frac{C}{\sqrt{1-r}}\leq \frac{C}{\delta}~. $$
Using the bounds $\|{\wb}^\perp\|\leq 1$ and $|\lambda|\leq C\delta$, we
find that:
\begin{eqnarray*}
  \|\Pi({\wb}\wedge \sigmab)\| & \leq & \|\Pi(({\wb}^\perp + \lambda \Rb)\wedge
  (\sigmab^\perp +\sigmab^r)) \| \\
  & \leq & \|\Pi({\wb}^\perp \wedge \sigmab^\perp)\| + |\lambda| \|\Pi(\Rb \wedge
  \sigmab^\perp)\| + \|\Pi({\wb}^\perp \wedge \sigmab^r)\| \\
  & \leq & C + C + C|\log(\delta)|~, 
\end{eqnarray*}
where the upper bound on the first term uses that ${\wb}^\perp\wedge
\sigmab^\perp$ is radial (so that $\Pi$ adds a factor $\delta$) while
the upper bound on the second term uses that $|\lambda|\leq
C\delta$. From this, we see that:
\beq \label{eq:i2} 
\|\Pi({\wb}\wedge \sigmab)\| \leq C|\log(\delta)|~. \eeq

For the third term in equation (\ref{eq:princ}), note that:
$$ \Pi(\omegab_\tau({\wb})) = \Pi(\omegab^\perp_\tau({\wb}^\perp)) + \lambda
\Pi(\omegab^\perp_\tau(\Rb)) +\Pi(\omegab^r_\tau({\wb}^\perp)) + \lambda
\Pi(\omegab^r_\tau(\Rb))~. $$
Lemma \ref{lm:omegab} shows that the first term is bounded by $C$, the
second by $\lambda C/\delta$, and thus by $C$, the third by $\delta
C/\delta=C$ because the projection induces a factor $\delta$, and the
last by $\lambda \delta C/\delta^2$ (again $\Pi$ adds a factor
$\delta$), and thus by $C$. Thus all 4 terms are bounded by $C$.

Using this along with (\ref{eq:i1}) and (\ref{eq:i2}) in equation
(\ref{eq:princ}) proves the Lemma.  
\epv

\paragraph{H{\"o}lder bounds on $\vb$}

We now state the key technical Lemma, about the H{\"o}lder bounds on $\vb$
near $\Lambda$. 

\blm \label{lm:holder}
For all $\epsilon \in (0,1)$, there exists a constant $C_\epsilon>0$ such
that, if $\xb, \yb \in \Sb$ are such that 
$d(\xb,\yb)\leq 1/2$, then $\|\vb(\xb)-\vb(\yb)\|\leq C_\epsilon
d(\xb,\yb)^\epsilon$.
\elm

We will use another point, called $\zb$, which we define as the point of
$[\xb,\yb]$ which is closest to $0$. We call $r_0:=r(\zb)$, and let
$\delta_0:=d(\xb,\yb)$. The key properties
of $\zb$ are the following.

\bprop \label{pr:zb}
\begin{enumerate}
\item if $\zb\neq \xb$, the angle at $\zb$ of the triangle $(0,\zb,\xb)$ is
  between $\pi/2$ and $\pi/2+C\sqrt{1-r_0}$. 
\item $\delta_0\leq C\sqrt{1-r_0}$.
\end{enumerate}
\eprop

\bpv
For the first point, note that, if $\zb\neq \yb$, the angle at $\zb$ of the
triangle $(0,\zb,\xb)$ --- call it $\alpha$ ---
is equal to $\pi/2$ since $\zb$ is at minimal distance
from $0$. Otherwise, $\zb=\yb$ is on $\dr\Omegab$, and, by convexity of
$\Omegab$, the segment $[\zb,\xb]$ is towards the interior of
$\Omegab$. Therefore Lemma \ref{lm:slope} shows that $\alpha$ is at most 
$\pi/2+C\sqrt{1-r_0}$.

For the second point, note that:
$$ r(\xb)^2 = r(\zb)^2 + 2\cos(\alpha) d(\zb,\xb) + d(\zb,\xb)^2\geq r(\zb)^2
+ d(\zb,\xb)^2~. $$ 
Since $r(\xb)<1$, so that:
$$ \sqrt{1 - r(\zb)^2} \geq d(\zb,\xb)~. $$
But $1-r(\zb)^2 \leq 2(1-r(\zb))$, it follows that $d(\zb,\xb)\leq C\sqrt{1 -
  r(\zb)}$. Since the same argument can be used with $\xb$ replaced by $\yb$,
we obtain the result.  
\epv

We now call $(\taub_0,\sigmab_0)$ the (unique) flat section of $\Eb$ which is
such that $\taub_0(\zb)=\taub(\zb)$ and that $\sigmab_0(\zb)=\sigmab(\zb)$; then
we call $\delta\taub :=\taub-\taub_0, \delta\sigmab := \sigmab-\sigmab_0$. 
Lemma \ref{lm:holder} will follow from the following estimates on
$\delta\taub$ and on $\delta\sigmab$.

\bprop \label{pr:delta}
We have the following estimates:
$$ \| \delta\sigmab(\xb)\|\leq \frac{C\delta_0}{1-r_0} +
C\left(\frac{1}{\sqrt{1-r(\xb)}} - \frac{1}{\sqrt{1-r_0}}\right)~, $$
$$ \| \delta\sigmab^r(\xb) \|\leq C\left| \log\left(
    \frac{1-r(\xb)}{1-r_0}\right)\right|~, $$
$$ \| \delta\taub(\xb)\| \leq C\left| \log\left(
    \frac{1-r(\xb)}{1-r_0}\right)\right|~, $$
$$ \| \delta\taub^\perp(\xb)\| \leq C\delta_0 + 
C\left(\sqrt{1-r_0}-\sqrt{1-r(\xb)}\right)~. $$
\eprop

\bpv
We consider two paths, such that the union is a path going from $\zb$ to
$\xb$. The first, $\gamma_\perp:[0,l]\rightarrow D ^3$, is a geodesic segment
on the sphere at distance $r_0$ from $0$, and goes from
$\zb$ to the point of $[0, \xb]$ at distance $r_0$ from $0$. The second,
$\gamma_r:[0,r-r_0]\rightarrow D^3$, goes from $\gamma_\perp(l)$ to $\xb$
along the segment joining them. 

We consider first the upper bound on $\delta\sigmab$. For all $s\in [0,l]$, we
have:
$$ \|\nabla_{\gamma_\perp'(s)}\delta\sigmab \| =
\|\omegab_\sigma(\gamma_\perp'(s))\| \leq \frac{C}{1-r_0} $$
by Lemma \ref{lm:omegab}, so that:
$$ \|\delta\sigmab(\gamma_\perp(l))\|\leq \frac{C\delta_0}{1-r_0}~. $$
For each $t\in [0,r-r_0]$, Lemma \ref{lm:omegab} shows that:
$$ \|\nabla_{\gamma_r'(t)}\delta\sigmab \| = 
\|\omegab_\sigma(\gamma_r'(t))\| \leq \frac{C}{\sqrt{1-r(\gamma_r(t))}^3}~, $$
so that, by integration:
$$ \|\delta\sigmab(\xb)\| \leq \frac{C\delta_0}{1-r_0} + 
C\left(\frac{1}{\sqrt{1-r(\xb)}} -
  \frac{1}{\sqrt{1-r_0}}\right)~. $$

To obtain the bound on $\delta\sigmab^r$, note that, for all $s\in [0,l]$:
$$ \|\nabla_{\gamma_\perp'(s)}\delta\sigmab^r \| \leq
\frac{1}{r_0} \|\delta\sigmab(\gamma_\perp(s))\| + 
\|\omegab^r_\sigma(\gamma_\perp'(s))\| 
\leq \frac{C}{\sqrt{1-r_0}}~. $$ 
Integrating this yields that:
$$ \|\delta\sigmab^r(\gamma_\perp(l))\|\leq C\frac{\delta_0}{\sqrt{1-r_0}}\leq
C~. $$
Moreover, for all $t\in [0,r-r_0]$:
$$ \|\nabla_{\gamma_r'(t)}\delta\sigmab^r \| = 
\|\omegab_\sigma^r(\gamma_r'(t))\| \leq \frac{C}{1-r(\gamma_r(t))}~, $$
and, by integration:
$$ \| \delta\sigmab^r(\gamma_\perp(l))\|\leq C + 
C\left|\log\left(\frac{1-r(\xb)}{1-r_0}\right)\right|~, $$
and the second point follows. 

For the third point, we use that, for all $s\in [0,l]$:
$$ \|\nabla_{\gamma_\perp'(s)}\delta\taub \| \leq
\|\omegab_\tau(\gamma_\perp'(s))\| + \|\delta\sigmab \|
\leq \frac{C}{\sqrt{1-r_0}} + \frac{C}{\sqrt{1-r_0}}~, $$
using the estimates obtained in the proof of point (1). So, by integration:
$$  \|\delta\taub(\gamma_\perp(l))\|\leq C\frac{\delta_0}{\sqrt{1-r_0}}\leq
C~. $$
For all $t\in [0, r-r_0]$:
$$ \|\nabla_{\gamma_r'(t)}\delta\taub \| = 
\|\omegab_\tau(\gamma_r'(t))\| + \|\delta\sigmab\| \leq 
\frac{C}{1-r(\gamma_r(t))} + \frac{C}{\sqrt{1-r(\gamma_r(t))}}~, $$
again by the estimate on $\delta\sigmab$ obtained above, so that, by
integration:
$$ \|\delta\taub(\xb))\|\leq
C\left|\log\left(\frac{1-r(\xb)}{1-r_0}\right)\right|~. $$

For the last point, we use that, for all $s\in [0,l]$:
$$ \|\nabla_{\gamma_\perp'(s)}\delta\taub^\perp \| \leq
\|\omegab_\tau^\perp(\gamma_\perp'(s))\| + \frac{1}{r_0}
\|\delta\taub(\gamma_\perp(s))\| + \|\delta\sigmab^r \|
\leq C + C + C~, $$
where we have used also the upper bound in the proof of point (2) above. 
It follows by integration that:
$$ \|\delta\taub^\perp(\gamma_\perp'(s))\|\leq C\delta_0~. $$
Moreover, for all $t\in [0,r-r_0]$:
$$ \|\nabla_{\gamma_r'(t)}\delta\taub^\perp \| \leq
\|\omegab_\tau^\perp(\gamma_r'(t))\| + \|\delta\sigmab\| \leq 
\frac{C}{\sqrt{1-r(\gamma_r(t))}} + \frac{C}{\sqrt{1-r(\gamma_r(t))}}~, $$
so that, by integration:
$$ \|\delta\taub^\perp(\xb))\|\leq C\delta_0 +
C\left(\sqrt{1-r_0}-\sqrt{1-r(\xb)}\right)~. $$
\epv

Those bounds are basically what we need to prove Lemma \ref{lm:holder} because
of the next proposition.

\bprop \label{pr:bornes-tech}
If $\delta_0$ is small enough, then:
$$ \sqrt{1-r_0}-\sqrt{1-r(\xb)}\leq C\delta_0~, $$ 
$$ \sqrt{1-r(\xb)}\left| \log\left(\frac{1-r(\xb)}{1-r_0}\right)\right|\leq
\delta_0|\log(\delta_0)|~. $$ 
\eprop

\bpv
By the upper bound on the angle at $\zb$, we have:
$$ r(\xb)\leq r_0 + C\delta_0\sqrt{1-r_0} +\delta_0^2~, $$
so that:
\beq \label{eq:ineq}
1-r(\xb) \geq (1-r_0) -C\delta_0\sqrt{1-r_0} - \delta_0^2~. \eeq
If $\sqrt{1-r_0}\leq C_0\delta_0$, then:
$$ \sqrt{1-r_0}-\sqrt{1-r(\xb)}\leq \sqrt{1-r_0}\leq C_0\delta_0~, $$
and:
$$ \sqrt{1-r(\xb)}\left| \log\left(\frac{1-r(\xb)}{1-r_0}\right)\right|\leq
\sqrt{1-r_0} |\log(1-r_0)|\leq C\delta_0|\log(\delta_0)|~, $$
so we only have to consider the case where $\delta_0/\sqrt{1-r_0}$ is small
enough to ensure --- through (\ref{eq:ineq}) --- that $1-r(\xb)\geq(1-r_0)/2$.

Then:
$$ \sqrt{1-r_0}-\sqrt{1-r(\xb)} = \frac{2}{1}\int_{r_0}^{r(\xb)}\frac{ds}{\sqrt{1-s}}
\leq \frac{r(\xb)-r_0}{2\sqrt{1-r(\xb)}} \leq \frac{r(\xb)-r_0}{\sqrt{1-r_0}}~. $$
But the angle at $\zb$ of the triangle $(0,\zb,\xb)$ is at most $C\delta_0$,
and, since $\delta_0\leq \sqrt{1-r_0}/C$, the angle between $[\zb,\xb]$ and
the lateral directions is everywhere at most $C\sqrt{1-r_0}$. Thus:
$$ r(\xb)-r_0 \leq C\delta_0\sqrt{1-r_0}~, $$
and it follows that:
$$ \sqrt{1-r_0}-\sqrt{1-r(\xb)} \leq \frac{r(\xb)-r_0}{\sqrt{1-r_0}} \leq
C\delta_0~. $$
Similarly:
$$ \left| \log\left(\frac{1-r(\xb)}{1-r_0}\right)\right| = 2 \left|
  \log\left(\frac{\sqrt{1-r(\xb)}}{\sqrt{1-r_0}}\right)\right| \leq 
2\int_{\sqrt{1-r(\xb)}}^{\sqrt{1-r_0}} \frac{ds}{s} \leq
  2\frac{\sqrt{1-r_0}-\sqrt{1-r(\xb)}}{\sqrt{1-r(\xb)}} \leq
  \frac{C\delta_0}{\sqrt{1-r(\xb)}}~, $$
which proves the second equation. 
\epv

\bpv[Proof of Lemma \ref{lm:holder}]
Recall that $\vb$ is the orthogonal projection on $\Sb$ of the vector field 
$\taub$. Clearly the orthogonal projection of $\taub_0$ is H\"older, because
$\taub_0$ is a Killing field and 
$\Sb$ is smooth enough by the results of section 2. So the result will follow
from the fact that the orthogonal projection on $T_{\xb}\Sb$ of 
$\delta\taub(\xb)$ is bounded by $C\delta_0|\log(\delta_0)|$, 
with the same bound applying with $\xb$ replaced by $\yb$. 

By definition,
$\delta\taub(\xb)=\delta\taub^\perp(\xb)+\delta\taub^r(\xb)$. The orthogonal
projection on $T_{\xb}\Sb$ of $\delta\taub^\perp(\xb)$ is bounded by
$\|\delta\taub^\perp(\xb)\|$, and so by:
$$ C\delta_0 + 
C\left(\sqrt{1-r_0}-\sqrt{1-r(\xb)}\right) $$
by Proposition \ref{pr:delta}, and thus by $C\delta_0$ by Proposition
\ref{pr:bornes-tech}. The orthogonal
projection on $T_{\xb}\Sb$ of $\delta\taub^r(\xb)$ is bounded by
$C\sqrt{1-r(\xb)}\|\delta\taub^\perp(\xb)\|$ by Lemma \ref{lm:slope}, so by:
$$ C\sqrt{1-r(\xb)} \left| \log\left(
    \frac{1-r(\xb)}{1-r_0}\right)\right| $$
by Proposition \ref{pr:delta}, and therefore by $C\delta_0 |\log(\delta_0)|$
    by Proposition \ref{pr:bornes-tech}. This proves the Lemma. 
\epv

\paragraph{H{\"o}lder estimates on the dual surface}

The argument given here extends with minor modifications to yield the
corresponding estimates on the deformation vector $\vb^*$ on the dual
surface $\Sb^*$. We will outline the proof here, with emphasis on the
small differences with the estimates on $\Sb$.

\bprop \label{pr:geom-dual}
Let $\xb$, and let $\xb^*$ be the dual point on
$\Sb^*$. Let $\theta(\xb^*)$ be the
angle between $T_{\xb^*}\Sb^*$ and the lateral directions, and 
let $(\xb^*-\xb)^\perp$ and
$(\xb^*-\xb)^r$ be the decomposition in lateral and radial components, {\it at
  $\xb$}, of the vector $\xb^*-\xb$. Then, for some constant $C>0$:
\begin{enumerate}
\item $\|(\xb^*-\xb)^r\|\leq C(1-r(\xb))$.
\item $\|(\xb^*-\xb)^\perp\|\leq C\sqrt{1-r(\xb)}$.
\item $|\theta(\xb^*)|\leq C\sqrt{1-r(\xb)}$.
\item the map sending a point of $\Sb$ to its dual in $\Sb^*$ is Lipschitz. 
\end{enumerate}
\eprop

\bpv
The projective definition of the duality shows that $\xb^*$ is in the plane 
dual to $\xb$. This plane
is exactly the set of points $m\in \R^3$ such that the orthogonal projection
of $m$ on the line $(0,\xb)$ is at distance
$1/(1-r(\xb))$ from $0$, this proves the first point. 

For the second point, let $\theta(\xb)$ be the angle between $T_{\xb}\Sb$ and
the lateral directions, then $\theta(\xb)$ is also --- by the projective
definition of the duality --- the angle at $0$ between the direction of $\xb$
and of $\xb^*$. This, along with the proof of the first point, shows the
second point. 

Finally, points (3) is a consequences of Lemma \ref{lm:slope} using the
definition of the duality, and point (4) follows from Lemma \ref{lm:slope} and
from Lemma \ref{lm:curvature}. 
\epv

There is also an analog of Lemma \ref{lm:H1-log}.

\blm \label{lm:H1-log*}
There exists a constant $C>0$ such that, on $\Sb^*$:
$$ \|\nablab \vb^*\|\leq C|\log(\delta)|~. $$
\elm

\bpv
It is a consequence of Remark \ref{rk:dual} that
$\vb^*$ is the orthogonal projection on $T\Sb^*$ of $
\taub(\xb) + \sigmab(\xb) \wedge (\xb^*-\xb)$. The Lemma then follows from
Lemma \ref{lm:omegab}, along with point (4) of the previous proposition. 
\epv

We can now state the dual of Lemma \ref{lm:holder}.

\blm \label{lm:holder*}
For all $\epsilon \in (0,1)$, there exists a constant $C_\epsilon>0$ such
that, if $\xb^*, \yb^* \in \Sb^*$ are such that 
$d(\xb^*,\yb^*)\leq 1/2$, then $\|\vb^*(\xb^*)-\vb^*(\yb^*)\|\leq C_\epsilon
d(\xb^*,\yb^*)^\epsilon$.
\elm

\bpv[Proof of Lemma \ref{lm:holder*}]
We use the same setting as in the proof of Lemma \ref{lm:holder}, i.e. $\zb$
is the point of $[\xb,\yb]$ which is closest to $0$, and we call
$(\taub_0,\sigmab_0)$ the flat section of $\Eb$ such that
$\taub_0(\zb)=\taub(\zb)$ and that $\sigmab_0(\zb)=\sigmab(\zb)$, and call
$\delta\taub =\taub-\taub_0, \delta\sigmab=\sigmab-\sigmab_0$. Clearly the
orthogonal projection on $T\Sb^*$ of the Killing field corresponding to
$(\taub_0,\sigmab_0)$ is H\"older, so the Lemma will follows --- as in the
proof of Lemma \ref{lm:holder} --- from the fact that the orthogonal
projection on $T_{\xb^*}\Sb^*$ of $\delta\taub(\xb) + \delta\sigmab(\xb)\wedge
(\xb^*-\xb)$ is bounded by $C\delta_0|\log(\delta_0)|$.

The bound on the orthogonal projection of $\delta\taub(\xb)$ is obtained
exactly like in the proof of Lemma \ref{lm:holder}, using the angle bound in
Proposition \ref{pr:geom-dual}. To bound the other term, use Proposition
\ref{pr:delta} to obtain that:
\begin{eqnarray*}
\|(\delta\sigmab(\xb)\wedge (\xb^*-\xb))^\perp \| & = &
\|\delta\sigmab^\perp (\xb)\wedge (\xb^*-\xb)^r + \delta\sigmab^r (\xb)\wedge
(\xb^*-\xb)^\perp \| \\
& \leq & \left(\frac{C\delta_0}{1-r_0} +
C\left(\frac{1}{\sqrt{1-r(\xb)}} - \frac{1}{\sqrt{1-r_0}}\right)\right)
(1-r(\xb)) + C\sqrt{1-r(\xb)} \left| \log\left(
    \frac{1-r(\xb)}{1-r_0}\right)\right| \\
& \leq & C\delta_0 + C\left(\sqrt{1-r_0}-\sqrt{1-r(\xb)}\right) +
C\sqrt{1-r(\xb)} \left| \log\left(
    \frac{1-r(\xb)}{1-r_0}\right)\right|~. 
\end{eqnarray*}
Using Proposition \ref{pr:bornes-tech}, we find that:
$$ \|(\delta\sigmab(\xb)\wedge (\xb^*-\xb))^\perp \|\leq
C\delta_0|\log(\delta_0)|~. $$
In a similar way:
\begin{eqnarray*}
\|(\delta\sigmab(\xb)\wedge (\xb^*-\xb))^r\| & = &
\|\delta\sigmab^\perp (\xb)\wedge (\xb^*-\xb)^\perp \| \\
& \leq & \left(\frac{C\delta_0}{1-r_0} +
C\left(\frac{1}{\sqrt{1-r(\xb)}} - \frac{1}{\sqrt{1-r_0}}\right)\right)
\sqrt{1-r(\xb))}~, 
\end{eqnarray*}
so the orthogonal projection of this term on $T_{\xb^*}\Sb^*$ is bounded by:
$$ C\delta_0 + C\left(\sqrt{1-r_0}-\sqrt{1-r(\xb)}\right)~, $$
which is bounded by $C\delta_0|\log(\delta_0)|$ by Proposition 
\ref{pr:bornes-tech}. 
\epv


\section{Infinitesimal rigidity}

We now have all the tools necessary to prove the infinitesimal rigidity
Lemma for hyperbolic manifolds with convex boundary. We first consider
the rigidity with respect to the induced metric, then with respect to
the third fundamental form.

\blm \label{lm:rigidity}
Let $g\in \cG$. For any non-trivial first-order deformation $\gd\in
T_g\cG$, the induced first-order variation of the induced metric on $\dr
M$ is non-zero. 
\elm

We give first a broad outline of the proof, then some additional
estimates, and finally the proof.

\pg{Outline of the proof}

The idea underlying the rigidity lemmas is that, once the problem has
been translated into a rigidity statement in $\R^3$, the known results
concerning the rigidity of convex surfaces in Euclidean 3-space should
hold, or at least their proofs should apply. In practice, things are
not quite as simple, because the vector fields $\ub$ and $\ub^*$
obtained in $\R^3$ are 
indeed infinitesimal deformations of $\Sb$ and $\Sb^*$, but they might
diverge on $\Lambda$; the results of \cite{Po} can therefore not
be applied directly. 

Recall that $\vb$ is the orthogonal projection of $\ub$ on $\Sb$, and that
$\wb:=B^{-1}\vb$. We have seen that $\wb$ 
is solution of $\db_{\III}\wb=0$, and that it satisfies a
H\"older estimate on $\Sb$. This estimate will allow us to prove that 
it is actually
in the Sobolev space $H^1$ on the whole sphere (including the limit set). To
prove that this equation has no solution in $H^1$
other than those coming from the
global infinitesimal isometries, we wish to consider the adjoint equation,
$d^{\nabla}h=0$. 

To have a good interpretation of this adjoint equation, however, it helps to
consider the equation $\db_{\III}\wb=0$ on a hyperbolic, rather than a
Euclidean, surface, since then the adjoint equation describes the
infinitesimal rigidity of the dual surface in the de Sitter space.
Therefore we will use {\it another} projective model of $H^3$, different from
the one used to define $\Sb$ from $S$, but with {\it the same ``center''}; we
decide that in this model the image of $H^3$ is a ball of radius $2$ in
$\R^3$. 

We call $\Sb_h$ the inverse image of $\Sb$ by this model, and $\Lambda_h$ the
set corresponding to the limit set of $M$. Thus $\Sb_h\cup \Lambda_h$ is a
closed, convex hyperbolic surface -- it does not get close to the boundary at
infinity of $H^3$.

The estimates in section 2 indicate that $\Sb$ has
principal curvatures bounded between two positive constants. 
Since $\Sb_h$ is away from the boundary of $H^3$, the projective model we use
does not distort the principal curvavures much, so that $\Sb_h$ also has
principal curvatures bounded between two positive constants. This means
that its second fundamental form is quasi-Lipschitz to its induced
metrics, and thus to the canonical metric on $S^2$. 

Let $\ub_h$ be the image of $\ub$ on $\Sb_h$ by the Pogorelov map, and let
$\vb_h$ be the orthogonal projection of $\ub_h$ on $\Sb_h$. By the basic
property of the 
Pogorelov map, $\ub_h$ is an isometric infinitesimal deformation of $\Sb_h$;
therefore, calling $B_h$ the shape operator of $\Sb_h$, $\wb_h:=B_h^{-1}\vb_h$
is a solution of $\db_{\III}\wb_h=0$ on $\Sb_h$. The
results of section 5 will be used to 
show that $\vb_h$ behaves rather well on $\Lambda_h$,
where it satisfies a H{\"o}lder estimate. 

The key of the argument will be that this fact implies that
$\wb_h$ is in  
the Sobolev space $H^1$ (for the metric $\III$), so that it
is actually not only a pointwise solution of $\db_{\III}\wb_h=0$ on
$\Sb_h$,
but also a global solution on $\Sb_h\cup \Lambda_h$, which is topologically a
sphere. This point uses the
fact that the Hausdorff dimension of $\Lambda$ is strictly less than
$2$ (a result of Sullivan \cite{sullivan}). 

Since $\db_{\III}$ is elliptic, it acts continuously from $H^1$ to $L^2$. 
We then move from $\db_{\III}$ to its adjoint $d^\nabla$; considering the dual
surface (in the de Sitter space) we show that it has no solution. This implies
that the co-kernel of 
$\db_{\III}$ restricted to $H^1$ is empty, and therefore that the
dimension of its kernel is equal to its index, namely $6$. Since there
is a 6-dimensional space of $H^1$ solutions which come from trivial
isometric 
deformations, all solutions are of this kind, and the conclusion follows.

\pg{Analytical properties of $\vb_h$}

We now need some precise informations on the analytical properties of
$\vb_h$, or more precisely on the vector field $\wb_h:=B_h^{-1}\vb_h$. 
Note that
the smoothness properties of $\wb$ will be with respect to the metric
$\III$, and to the connection $\nablat$. 

Note first that $\vb_h$ satisfies the same H\"older estimate, near $\Lambda_h$,
as the estimate on $\vb$ given by Lemma \ref{lm:holder}. This follows from the
estimates of section 5 on the radial and the lateral components of $\ub$;
since we have kept the same ``center'' in the projective model, and since
$\Sb_h$ is compact, the radial and lateral components of $\ub_h$ satisfies
the same estimates (up to constants) as the corresponding components of $\ub$,
and the H\"older estimate on $\vb_h$ follows. 

\bdf 
We call $H^1_{\III}(T\Sb_h)$ the space of sections of $T\Sb_h$ which are in
the Sobolev space $H^1$ for the third fundamental form $\III$ of $\Sb_h$.
\edf

Here $H^1$ is the space of functions $u$ such that the integral of the norm of
$du$ converges. 
We first state an elementary technical proposition useful later on. 

\bprop \label{pr:log}
Let $\delta_{\III}$ be the distance to $\Lambda_h$, for $\III$, on $\Sb_h$. 
Then $\log(\delta_{\III})$ is in $L^2$ for the area form of the third
fundamental form of $\Sb_h$. 
\eprop

\bpv
By a result of Bishop \cite{bishop-minkowski} that the Minkowski dimension
of $\Lambda$ (as a subset of $S^2$) is strictly smaller than $2$ (it was
earlier proved by Sullivan for the Hausdorff dimension). It follows that the
Minkowski dimension of $\Lambda$ in $\Sb$ is also strictly smaller than
$2$. Since the projective map we use does not distort the distance much on
$\Sb$, it follows that the Minkowski dimension of $\Lambda_h$ in
$\Sb_h\cup\Lambda_h$ is also strictly smaller than $2$, and, since the
principal curvature are bounded between two positive constants, the same holds
for the third fundamental form. The result then
follows by an elementary integration argument.
\epv

Note that the proof could also be done using the fact that the Hausdorff
dimension of $\Lambda$ is strictly smaller than $2$, but it is then
technically more tricky.

\bcr \label{cr:wb}
The integral of the square of the norm of $\nablat \wb_h$ over $\Sb_h$
converges, and, for all $\alpha\in (0,1)$, $\wb$ is $C^\alpha$ H{\"o}lder on
$\Sb_h\cup \Lambda_h$. 
\ecr

\bpv
We have seen that $\vb_h$ satisfies the same H\"older estimate as $\vb$ in 
Lemma \ref{lm:holder}; but
$\wb_h:=B_h^{-1}\vb_h$, so that, for any $X\in T\Sb_h$: 
$$ \nablat_X\wb_h = B_h^{-1}\nabla_X(B_h\wb_h) = B_h^{-1}\nabla_X\vb_h~, $$
so (since the eigenvalues of $B_h$ are bounded and $\III$ is quasi-Lipschitz to
$I$) $\wb_h$ is also $C^\alpha$-H{\"o}lder, but for $\III$. 

Moreover, $\|\nabla\vb\|$ is bounded by
$C|\log(\delta)|$ by Lemma \ref{lm:H1-log}, so, since the radial and lateral
components of $\ub_h$ behave as the corresponding components of $\ub$, the
same estimate applies to $\vb_h$; it follows that $\|\nablat
\wb_h\|_{\III}\leq C|\log(\delta_{\III})|$. 
Lemma \ref{pr:log} thus shows that
the integral of the square of the norm of $\nablat \wb_h$ over $\Sb_h$
converges. 
\epv

This corollary will be helpful thanks to the following simple
lemma. It is probably well-known to the specialists --- and can be
extended in various directions --- but we have included a proof for
completeness. 

\blm \label{lm:H1}
Let $\Lambda\subset \Sigma$ be a subset of Hausdorff dimension
$\delta>1$ in a closed surface. 
Let $u:\Sigma\rightarrow \R$ be a $L^2$ function such that:
\begin{itemize}
\item $u$ is everywhere $C^\alpha$, for some $\alpha\in (\delta-1, 1)$. 
\item $u$ is smooth outside $\Lambda$.
\item we have:
$$ \int_{\Sigma\setminus \Lambda} \| du\|^2 da <\infty~. $$
\end{itemize}
Then $u$ is in the Sobolev space $H^1(\Sigma)$. 
\elm

\bpv
The statement is essentially local, so it is sufficient to prove the same
result in $\Omega:=[0,1]\times [0,1]$, with the
support of $u$ in $(0,1)\times (0,1)$. Let $v$ be a smooth vector field
on $\Omega$.  

Choose $\epsilon>0$. By definition of
the Hausdorff dimension of $\Lambda$, 
there exists a covering of $\Lambda$ by geodesic balls $B_1,
\cdots, B_N$,  with $B_i$ of center $x_i$ and radius $r_i\leq \epsilon$,
such that $\sum_i 
r_i^{1+\alpha} \leq \epsilon$.  Let $B:=\Omega\cap (\cup_i B_i)$. We replace
some of the balls $B_i$ 
by sub-domains, and find a covering of $B$ by a finite number of
convex domains $D_i$, with disjoint interiors, so that
$D_i\subset B_i$ for each $i$. Then: 
$$ \forall i\in \{ 1, \cdots, N\}, ~L(\dr D_i) \leq L(\dr B_i) \leq
2\pi r_i~. $$

Let $v$ be a smooth vector field on $S^2$. Then:
\begin{eqnarray*}
\left| \int_{\Omega\setminus B} u \dv(v) + du(v) da \right| & = & \left|
\int_{\Omega\setminus B} \dv(uv) da \right| \\
& = & \left| \int_{\dr B} u\langle v, n\rangle dx \right| \\
& = & \left| \sum_i \int_{\dr D_i} u\langle v, n\rangle dx \right| \\
& \leq & \sum_i \left| \int_{\dr D_i} (u(x_i)+(u(x)-u(x_i)))\langle
v(x_i) + (v(x)-v(x_i)), n\rangle dx \right|~.
\end{eqnarray*}
Note that:
$$ \int_{\dr D_i} \langle v(x_i), n\rangle = \int_{D_i} \div (v(x_i)) da = 0~, $$
so that:
$$ \left| \int_{\Omega\setminus B} u \dv(v) + du(v) da \right| \leq \sum_i
\left| \int_{\dr D_i} u(x_i) \langle v(x_i), n\rangle 
\right| + \coupeq +
\int_{\dr D_i} | u(x)\langle v(x)-v(x_i), n\rangle +
(u(x)-u(x_i)) \langle v(x_i), n\rangle | dx 
\leq 2\pi \sum_i C_1 \max |u| r_i^2 + C_2 \max \| v\| r_i^{1+\alpha}~, $$
where $C_1$ is the maximum of $|Dv|$ and $C_2$ is the H{\"o}lder modulus of
$u$. 

This means that, if $\epsilon$ is small enough, we have for some $C_3>0$:
$$ \left| \int_{\Omega\setminus B} u \dv(v) + du(v) da \right| \leq
C_3 \sum_i r_i^{1+\alpha}\leq C_3 \epsilon~. $$
Taking the limit as $\epsilon \rightarrow 0$, we see that:
$$ \int_{\Omega\setminus \Lambda} u \dv(v) da = - \int_{\Omega\setminus
\Lambda}  du(v) da~. $$
Since the $H^1$ norm of $u$ over $\Omega\setminus \Lambda$ is finite, we
see that:
$$ \left| \int_{\Omega} u \dv(v) da \right| = \left| \int_{\Omega\setminus
\Lambda}  du(v) da \right| \leq \| u\|_{H^1(\Omega\setminus \Lambda)} \|
v\|_{L^2}~. $$
This inequality, which is true for all smooth vector field $v$, shows
that $u$ is in $H^1(\Omega)$, as needed.  
\epv

As a direct consequence of Corollary \ref{cr:wb} and of Lemma \ref{lm:H1}, we
have the following important statement.

\blm \label{lm:wh1}
$\wb_h$ is in $H^1_{\III}(T(\Sb_h\cup\Lambda_h))$, and it is a solution of
$\db_{\III}\wb_h=0$ on $\Sb_h\cup \Lambda_h$. 
\elm

This is important because $H^1_{\III}(T(\Sb_h\cup\Lambda_h))$ 
is a space on which $\db_{\III}$ acts well.

\blm \label{lm:cont-h1}
$\db_{\III}$ acts continously from $H^1_{\III}$ to
$L^2(\Omega^{0,1}_{\II}(\Sb_h\cup\Lambda_h))$.
\elm

\bpv
Let $\zb\in H^1_{\III}$. Let $U\subset \Sb_h\cup\Lambda_h$ 
be an open subset, and let $X$ be
a smooth unit vector field on $U$. Then, by definition of $H^1_{\III}$, both
$\nablat_X\zb$ and $\nablat_{\Jb X}\zb$ are $L^2$ sections of $TU$. So
$\db_{\III}\zb$ is a $L^2$ section of $TU$. Since this is true for any
simply connected open subset $U\subset \Sb_h\cup\Lambda_h$, 
$\db_{\III}\zb$ is a $L^2$
section of $\Omega^{0,1}_{\II}(\Sb_h\cup\Lambda_h)$. 
\epv

\pg{$L^2$ solutions of $d^{\nabla}$}

The other important point of the proof is that there is no $L^2$
solution of 
$d^\nabla h=0$ in $\Omega^{0,1}_{\II}(\Sb_h\cup\Lambda_h)$. 
We first prove the corresponding result for $d^{\nablat}h=0$ on Euclidean
surfaces, following closely the ideas 
used in section 3; we only have to keep closely track of the
smoothness issue. 

\blm \label{lm:L2dD}
Let $S$ be a closed, strictly convex surface in $\R^3$, with principal
curvatures bounded between two strictly positive constants. 
There is no non-trivial solution of $d^{\nablat} h=0$ with $h\in L^2(\ots)$.
\elm

\bpv
We will follow the proof of Proposition \ref{pr:defos-isom} and Lemma
\ref{lm:rig-smooth}, but we will have to be more careful because we have
only a limited degree of smoothness --- actually only the hypothesis
that $h$ is $L^2$.

First note that, according to Proposition \ref{pr:equivalence}, it is
equivalent to show that there is no non-trivial $L^2$ solution of equation
(\ref{eq:Bd}), because, since $B$ has its eigenvalues bounded between
two positive constants, $\Bd$ is in $L^2$ if and only if $B^{-1}\Bd$ is. So
we consider a $L^2$ solution $\Bd$ of (\ref{eq:Bd}), and will prove that
it is $0$.

We now follow the proof of Proposition \ref{pr:defos-isom}. We
consider $S$ as the image of a map $\phi:S\rightarrow \R^3$,
and 
its unit normal vector field as a map $N:S\rightarrow \R^3$. We
define a $\R^3$-valued 1-form $\alpha$ by $\alpha:=(d\phi\circ
\Bd)\wedge N$. Then, as in the proof of Proposition \ref{pr:defos-isom},
$d\alpha=0$ because $\Bd$ is a solution of (\ref{eq:Bd}), and we can
integrate it to obtain a $H^1$ function $Y:S\rightarrow \R^3$. We can
then define another $\R^3$-valued 1-form $\beta$ on $S$ by:
$\beta:=Y\wedge d\phi$, check that it is closed, so that it can be
integrated as a function $\phid:S\rightarrow \R^3$, which defines an
isometric first-order deformation of $\phi$ inducing the
first-order deformation $\Bd$ of $B$.

We now consider the function $f:=r^2/2$ on $\phi(S)$, and its pull-back by
$\phi$ on $S$ (which we still call $f$). The infinitesimal deformation
$\phid$ induces a first-order variation $\fd$ of $f$, which is in $H^1$
because $Y$ is in $H^1$. It is then possible to follow through the proof
of Proposition \ref{pr:defos-isom}, the key point is that, since $\Bd$
is $L^2$, its determinant is $L^1$, so that $d\omega$ is also $L^1$ and
the integration argument works. 
\epv

We now use this result to prove that there is no non-trivial $L^2$ 
solution of $d^\nabla h=0$ on $\Sb_h\cup \Lambda_h$. The proof uses the
duality with the de Sitter space, and the Pogorelov map.

\bcr \label{cr:L2dD}
There is no $L^2$ solution of $d^\nabla h=0$ with $h\in
\Omega^{0,1}_{\II}(\Sb_h\cup\Lambda_h)$. 
\ecr

\bpv
Consider the surface $S^*$ in the de Sitter space which is dual to $\Sb_h\cup
\Lambda_h$. Since $\Sb_h\cup\Lambda_h$ is closed,
strictly convex, and has its second fundamental form bounded 
between two strictly
positive constants, $S^*$ is a closed, space-like, strictly convex surface
with its second fundamental form bounded 
between two positive constant. Since $\Sb_h\cup\Lambda_h$ is a
closed surface in $H^3$, we can choose the projective model in such a way
that the image 
of $S^*$ in this model remains in a compact subset of $\R^3$, and away from
the unit ball containing the projective model of $H^3$.

As pointed above, a solution of $d^\nabla h=0$ on $\Sb_h\cup\Lambda_h$ 
with $h\in L^2(\Omega^{0,1}_{\II}(\Sb_h\cup\Lambda_h))$
translates as a solution of $d^{\nablat}h=0$ on $S^*$, with $h\in
L^2(\Omega_{\II}^{0,1}(S^*))$. If $B^*$ is the shape operator of $S^*$, then
$\Bd^*:=B^*h$ is a solution of equation (\ref{eq:Bd}) on $S^*$, which
describes a first-order infinitesimal deformation of $S^*$.

Consider the image $S^*_e$ of $S^*$ in the projective model of $S^3_1$; as
mentioned above, we choose the projective model so that $S^*_e$ is a closed
surface in $\R^3$. $S^*_e$ remains outside the unit sphere, so the projective
model does not change the principal curvatures much, and therefore the second
fundamental form of $S^*_e$ is bounded between two positive constants. 
Applying the Pogorelov map $\Phi_S$ (see Definition \ref{df:pogo-S})
to the first-order deformation
of $S^*$ obtained above, we get a first-order isometric deformation of
$S^*_e$. Since $S^*_e$ is in a compact set of $\R^3$ away from the unit one
ball, this first-order deformation again corresponds to a solution of equation
(\ref{eq:Bd}) on $S^*_e$ which is in $L^2$. But there is no non-trivial such
solution by the previous lemma. 
\epv

\pg{Proof of the Rigidity Lemma for $I$}

Consider an infinitesimal deformation $\gd$ of the hyperbolic metric on
$M$, which does not change the induced metric on $\dr M$. The construction of
section \ref{se:2b} leads to a vector field $u$ on the boundary $S$ of the
universal cover of $M$, seen as a surface in $H^3$. Using the Pogorelov map,
it translates as an infinitesimal deformation $\ub$ 
of the surface $\Sb$ in $\R^3$. Recall that $\Sb$ is a convex 
surface with second fundamental form bounded between two positive constants,
which is tangent to the unit sphere along the image of the limit set of $M$.

We now take the orthogonal projection $\vb$ of $\ub$ on $\Sb$, and set
$\wb:=B^{-1}\vb$. Then $\wb$ is a solution of $\db_{\III}\wb=0$ on $\Sb$. 
We then call $\ub_h$ be the image of $\ub$ on $\Sb_h$ 
by the Pogorelov transformation, so $\ub_h$ is a first-order isometric
deformation of $\Sb_h$. 
Consider the orthogonal projection $\vb_h$ of $\ub_h$ on
$\Sb_h$, and let $\wb_h:=B_h^{-1}\vb_h$, where $B_h$ is the shape operator of
$\Sb_h$. Then, from Lemma \ref{lm:wh1}, $\wb_h$ is in $H^1$ for $\III$ and is a
solution of $\dr_{\III}\wb_h=0$ on $\Sb_h\cup \Lambda_h$.

By Lemma \ref{lm:elliptic}, $\db_{\III}$ is an elliptic operator of index
$6$ on $\Sb_h\cup\Lambda_h$; it acts (continously by Lemma \ref{lm:cont-h1})
from the space of $H^1$ sections of $T(\Sb_h\cup\Lambda_h)$ 
(for $\III$) to the space
of $L^2$ sections of $\Omega_{\II}^{0,1}(\Sb_h\cup\Lambda_h)$. 
By Proposition \ref{pr:adjoint}, its adjoint (or rather the
adjoint of $\Jb\db_{\III}$) is $-d^\nabla$. But we know (by Corollary
\ref{cr:L2dD}) that the kernel of $d^\nabla$ restricted to the $L^2$
sections on $\Sb_h\cup\Lambda_h$ is $0$. 

Therefore, the co-kernel of $\db_{\III}$ is $0$, and its kernel has dimension
$6$. But there is a 
6-dimensional vector space of trivial solutions, corresponding to the
global Killing vector fields of $H^3$. Therefore, all solutions of
$\db_{\III}$ are trivial. In particular, $\ub_h$ is a trivial deformation, so
that $\ub$  is trivial, and thus
$\gd$ is trivial. This proves Lemma \ref{lm:rigidity}.

\pg{The dual rigidity}

The proof of the infinitesimal rigidity with respect to the induced
metric also applies, with minor modifications, for the third fundamental
form. The most important difference is that one has to use the
hyperbolic-de Sitter duality, and the problem is moved to $S^3_{1,+}$,
and then to the exterior of the unit ball 
in $\R^3$. Here is the analog of Lemma \ref{lm:rigidity}. 

\blm \label{lm:rigidity*}
Let $g\in \cG$. For any non-trivial first-order deformation $\gd\in
T_g\cG$, the induced first-order variation of the third fundamental form
of $\dr M$ is non-zero. 
\elm

The proof follows the same steps as the proof of Lemma
\ref{lm:rigidity}. The key point is of course Lemma \ref{lm:holder*}.
The proof of Lemma \ref{lm:wh1}, with minor changes, also shows the
same result for the dual 
surface. Of course, we call $H^1_{\III}(T\Sb^*)$ the space of sections of
$T\Sb^*$ which are in the Sobolev space $H^1$ for the third fundamental
form of $\Sb^*$. We let $\wb^*:=(B^*)^{-1}\vb^*$, and call $\db_{\III}^*$ the
analog of $\db_{\III}$ on the dual surface $\Sb^*\cup \Lambda$.

\blm \label{lm:solh1*}
$\wb^*$ is in $H^1_{\III}(T\Sb^*)$, and it is a solution of
$\db^*_{\III}\wb^*=0$ on $\Sb^*$.  
\elm

The end of the proof also follows closely the argument given above, with
$\Sb$ replaced by $\Sb^*$; in the final step, one should associate to $\Sb^*$
a closed hyperbolic surface $\Sb^*_h$, as $\Sb_h$ for $\Sb$. 
We leave the details to the reader.


\section{Compactness}

We move here to the proof of the properness of $\cF$ and
$\cF^*$. As mentioned in the introduction, those properties are
equivalent to the following lemmas, which we will prove at the end of
the section. 

\blm \label{lm:compact}
Let $(g_n)_{n\in \N}$ be a sequence of hyperbolic metrics on $M$ with
smooth, convex boundary. Let $(h_n)_{n\in \N}$ be the sequence of
induced metrics on the boundary. Suppose that $(h_n)$ converges to a
smooth metric $h$ on $\dr M$, with curvature $K>-1$. Then there is a
sub-sequence of $(g_n)$ which converges to a hyperbolic metric $g$ with
smooth, convex boundary. 
\elm

\blm \label{lm:compact*}
Let $(g_n)_{n\in \N}$ be a sequence of hyperbolic metrics on $M$ with
smooth, convex boundary. Let $(h_n)_{n\in \N}$ be the sequence of
third fundamental forms of the boundary. Suppose that $(h_n)$ converges to a
smooth metric $h$ on $\dr M$, with curvature $K<1$ and closed
geodesics which are of length $L>2\pi$ when they are contractible in
$M$. Then there is a 
sub-sequence of $(g_n)$ which converges to a hyperbolic metric $g$ with
smooth, convex boundary. 
\elm

\pg{Convex surfaces and pseudo-holomorphic curves}

We will state below some compactness results for isometric embeddings of
convex surfaces in $H^3$ or in $S^3_1$. Although we will not prove those
lemmas here, it should be pointed out that the proofs rest on the use of
a compactness result for pseudo-holomorphic curves. 
Namely, given a surface $\Sigma$
with a Riemannian metric of curvature $K>K_0$ (resp. $K<K_0$) and a
Riemannian (resp. Lorentzian) 
3-manifold $M$ of sectional curvature $K<K_0$ (resp. $K>K_0$), the isometric
embeddings of $\Sigma$ in $M$ lift to pseudo-holomorphic curves in a
natural space of 1-jets of infinitesimal isometric embeddings of
$\Sigma$ in $M$, which comes equipped with an almost-complex structure
defined on a distribution of 4-planes. This idea was outlined by Gromov
in \cite{G1}, and developed in \cite{L1,L-MA,these,cras}.

\pg{Compactness of surfaces in $H^3$}

The compactness lemma for $I$ is based on the following statement.

\blm[\cite{L-MA}] \label{lm:cph-h}
Let $(h_n)_{n\in \N}$ be a sequence of metrics on $D^2$ with curvature
$K>-1$, converging to a smooth metric $h$ with curvature $K>-1$. Suppose
that $d_h(0, \dr D^2)\geq 1$. Let $x_0\in H^3$, and let $(f_n)_{n\in
  \N}$ be a sequence of isometric embeddings, $f_n:(D^2, h_n)\rightarrow
H^3$, with $f_n(0)=x_0$. Then:
\begin{itemize}
\item either there is a subsequence of $(f_n)$ which converges
  $C^\infty$ to an isometric embedding of $(D^2, h)$ in $H^3$.
\item or there exists a geodesic segment $\gamma\ni 0$ on which
  a subsequence of $(f_n)$ converges to an isometry to a geodesic
  segment of $H^3$; and, in addition, the integral of the mean
  curvature $H_n$ on each 
  segment tranverse to $\gamma$ diverges as $n\rightarrow \infty$.
\end{itemize}
\elm

We now come back to the setting of Lemma \ref{lm:compact}, and consider a
sequence $(g_n)_{n\in \N}$ of hyperbolic metrics on $M$. We suppose that
the sequence of induced metrics $(h_n)$ on $\dr M$ converges to a limit
$h$ which is smooth, with curvature $K>-1$. Note that we can apply a
sequence of global isometries of $H^3$ to fix one point of $\dr
\Mt$. Using the next proposition, this means that, after taking a
subsequence, we can suppose that 
each connected component of $\dr M$ has a lift containing a point which
converges in $H^3$. 

\bprop \label{pr:diam}
Under the hypothesis of Lemma \ref{lm:compact}, the diameter of $(M,
g_n)$ remains bounded.
\eprop

\bpv
Since the induced metrics on $\dr M$ converge, the diameter of each
boundary component, for the induced metrics, is bounded. But the
closest-point projection from the boundary components of $M$ to the
corresponding components of the boundary $\dr C(M)$ of the convex core
of $M$ is a contraction, so that the induced metrics on $\dr C(M)$ have
bounded diameter. Since those metrics are hyperbolic, this means that
they remain in a compact subset of Teichm{\"u}ller space. 

Now the map sending the conformal structure at infinity on $\dr_\infty
E(M)$ to the conformal structure of the induced metric on $\dr C(M)$
is proper (see \cite{bridgeman-canary}, and \cite{epstein-marden} when $M$ has
incompressible boundary). So the conformal structures on $\dr_\infty E(M)$
also remain bounded, and, by the Ahlfors-Bers theorem, the sequence of
metrics on $E(M)$ remains in a compact set. In particular, the diameter
of the convex cores of $(M, g_n)$ remain bounded. 

Now an elementary argument shows that, if the maximal distance to $C(M)$
went to infinity on a connected component of $\dr M$, then the diameter
of this connected component would go to infinity. Therefore, the
distance of the points of $\dr M$ to $C(M)$ remains bounded, and the
diameter of $(M, g_n)$ is also bounded. 
\epv

We can now use Lemma \ref{lm:cph-h} to state a compactness result for the
lifts in $H^3$ of the boundary components of $M$. 

\bcr \label{cr:cv-h}
For each connected component $S_i$ of the boundary of the universal 
cover $\Mt$ of $M$, there exists a subsequence of $(g_n)$ and a sequence
$(\rho_n)$ of isometries such that the image by $\rho_n$ of $S_i$
converges $C^\infty$ on compact subsets to an equivariant surface in $H^3$. 
\ecr

\bpv
Choose a point $x_0\in S_i$, and apply Lemma \ref{lm:cph-h} to a disk
$D_0$ of radius $1$ for $h$, centered at $x_0$. 
Suppose that the second alternative of the lemma applies, and consider a
small neighborhood of $\gamma$. The integral of $H_n$ on segments
transverse to $\gamma$ would diverge, and this would contradict the fact
that $\dr M$ --- and thus $S_i$ --- is locally convex and embedded.
So the first case of Lemma \ref{lm:cph-h} applies.

Now choose a point $x_1$ on the boundary of $D_0$, and apply again
Lemma \ref{lm:cph-h} to a disk of radius $1$ centered at $x_1$. The same
argument as above shows that, after taking a subsequence, $S_i$
converges $C^\infty$ on $D_0\cup D_1$. 

Applying recursively the same argument leads to an increasing sequence
of compact subsets of $S_i$ on which the sequence converges $C^\infty$. 
\epv

The proof of the main compactness lemma for the induced metric follows. 

\bpn{of Lemma \ref{lm:compact}}
Corollary \ref{cr:cv-h} shows that, for each end of $E(M)$, the
corresponding boundary component of $\Mt$ converges $C^\infty$, so that
its representation converges. This shows
that each end of $E(M)$ converges, so that the hyperbolic metric on
$E(M)$ converges. Since each 
boundary component of $M$ converges in $E(M)$, the hyperbolic metrics on
$M$ clearly converge. 
\epn

\pg{Compactness in $S^3_1$}

We outline here the proof of the compactness lemma with respect to the
third fundamental form; it follows the proof of Lemma \ref{lm:compact},
with some natural modifications. It is based on a compactness result for
sequences of isometric embeddings of surfaces in Lorentzian 3-manifolds.

\blm[\cite{these,cras}] \label{lm:cph-s}
Let $(h_n)_{n\in \N}$ be a sequence of metrics on $D^2$ with curvature
$K<1$, converging to a smooth metric $h$ with curvature $K<1$. Suppose
that $d_h(0, \dr D^2)\geq 1$. Let $x_0\in S^3_1$, let $P\subset
T_{x_0}S^3_1$ be a space-like plane, and let $(f_n)_{n\in
  \N}$ be a sequence of isometric embeddings, $f_n:(D^2, h_n)\rightarrow
S^3_1$ with $f_n(0)=x_0$ and $\mbox{Im}(d_0f_n)=P$. Then:
\begin{itemize}
\item either there is a subsequence of $(f_n)$ which converges
  $C^\infty$ to an isometric embedding of $(D^2, h)$ in $S^3_1$.
\item or there exists a geodesic segment $\gamma\ni 0$ on which a
  subsequence of 
  $(f_n)$ converges to an isometry to a geodesic segment of $S^3_1$.
\end{itemize}
\elm

The main consequence is that, under the hypothesis of Lemma
\ref{lm:compact*}, the lifts in $H^3$ of each boundary component
converges. 

\bcr \label{cr:cv-s}
Suppose that $(\dr \Mt, h)$ has no closed geodesic of length $2\pi$. Then,
for each connected component $S_i$ of the boundary of the universal 
cover $\Mt$ of $M$, there exists a subsequence of $(g_n)$ and a sequence
$(\rho_n)$ of isometries such that the image by $\rho_n$ of $S_i$
converges $C^\infty$ to a smooth equivariant surface in $H^3$. 
\ecr

\bpv
We basically follow the proof of Corollary \ref{cr:cv-h} (but in $S^3_1$
instead of $H^3$), with a slightly
simpler argument because, if the second alternative in Lemma
\ref{lm:cph-s} happened, there should be a geodesic $\gamma$ of $S_i$
which is sent by the dual of the limit map to a geodesic of $S^3_1$. Since all
geodesics of $S^3_1$ are closed and have length $2\pi$, $\gamma$ should
also be closed and of length $2\pi$, and this is impossible.
\epv

The proof of Lemma \ref{lm:compact*} follows, along the lines of the
proof of Lemma \ref{lm:compact}. 

\bpn{of Lemma \ref{lm:compact*}}
Corollary \ref{cr:cv-s} shows that, for each boundary component of $M$,
the equivariant embedding of the corresponding boundary component of
$\Mt$ converges.

In particular, the induced metrics on all the boundary
components of $M$ converge, and this --- along with Lemma \ref{lm:compact} --- 
proves the result.
\epn


\section{Spaces of metrics}

The last step of the proof is to set up deformation arguments for $\cF$
and for $\cF^*$. The main point of this section is that the spaces of
metrics which enter the picture are topologically simple. 

\blm \label{lm:cnx}
Both $\cH$ and $\cH^*$ are connected and simply connected.
\elm 

The proof is elementary, and it is the same as what is done for the
special case when $\dr M$ is a sphere; for $\cH^*$, see
\cite{RH,these}. 

\bpv
We first prove that $\cH$ is connected. 
Let $h_0, h_1\in \cH$. Since the space of metrics on $\dr M$ is
connected, there is a one-parameter family $(h_t)_{t\in [0,1]}$ of
metrics on $\dr M$ connecting $h_0$ to $h_1$. But then there exists a
constant $C>0$ such that, for all $t\in [0,1]$, $h_t$ has curvature
$K>-C^2/2$. In addition, there is another constant $\epsilon>0$ such that,
for all $t\in [0, \epsilon]\cup [1-\epsilon, 1]$, $h_t$ has curvature
$K>-1$. We define a new family $(\hb_t)_{t\in [0,1]}$ as follows:
\begin{itemize}
\item for $t\in [0, \epsilon]$, $\hb_t=h_t(1+(C-1)t/\epsilon)^2$.
\item for $t\in [\epsilon, 1-\epsilon]$, $\hb_t = C^2h_t$.
\item for $t\in [1-\epsilon, 1]$, $\hb_t=h_t(1+(C-1)(1-t)/\epsilon)^2$. 
\end{itemize}
It is then clear that, for each $t\in [0,1]$, $\hb_t$ has curvature
$K>-1$, and this proves that $\cH$ is connected. 

We proceed similarly to prove that $\cH$ is simply connected. Let
$h:S^1\rightarrow \cH$ be a closed path in $\cH$. Since the space of
metrics (defined up to isotopy) on $\dr M$ is contractible, there is a
map $f$ from $D^2$ to the space of smooth metrics on $\dr M$ such that, for 
each point $x\in S^1$, $f(x)=h(x)$. Then:
\begin{itemize}
\item there exists a constant $C>0$ such that, for each $y\in D^2$,
  $f(y)$ has curvature $K>-C^2/2$.
\item there is another constant $\epsilon>0$ such that, for each $y\in
  D^2$ at distance at most $\epsilon$ from $S^1=\dr D^2$, $f(y)$ has
  curvature $K>-1$.
\end{itemize}
We define another map $\fb$ from $D^2$ to the space of smooth metrics on
$\dr M$ as follows:
\begin{itemize}
\item for $y\in D^2$ at distance at most $\epsilon$ from $\dr D^2$,
$\fb(y)=f(y)(1+(C-1)d(y, \dr D^2)/\epsilon)^2$.
\item for $y\in D^2$ at distance at least $\epsilon$ from $\dr D^2$,
  $\fb(y)=C^2f(y)$. 
\end{itemize}
It is then clear that $\fb$ has values in $\cH$, so that $\cH$ is
simply connected.

\medskip

To proof that $\cH^*$ is connected, we proceed as in the proof that
$\cH$ is connected, with a minor difference; given $h_0, h_1\in \cH^*$,
we choose a path $(h_t)_{t\in [0,1]}$ connecting them, and we pick $C>0$
such that, for each $t\in [0,1]$, $h_t$ has curvature $K<C^2/2$ and
closed geodesics of length $L>4\pi/C$. We also choose $\epsilon>0$ such
that, for $t\in [0,\epsilon]\cup [1-\epsilon, 1]$, $h_t\in \cH^*$. We
then define $(\hb_t)_{t\in [0,1]}$ by scaling exactly as above, and
check that, for all $t\in [0,1]$, $\hb_t\in \cH^*$. So $\cH^*$ is
connected. 

Finally, proving that $\cH^*$ is simply connected is similar to proving
that $\cH$ is simply connected. Given a closed path $h:S^1\rightarrow
\cH^*$, we choose a map $f$ from $D^2$ to the space of smooth
metrics on $\dr M$ which coincides with $f$ on $\dr D^2$. We then pick a
constant $C>0$ such that, for all $y\in D^2$, $f(y)$ has curvature
$K<C^2/2$, and has closed geodesics of length $L>4\pi/C$. The scaling
argument is then the same as for $\cH$ above.
\epv

We can also remark that $\cG$ is connected; this is a direct consequence
of the Ahlfors-Bers theorem, which shows that the space of complete,
convex co-compact metrics on $E(M)$ is connected. 


\section{Proof of the main theorems} \label{se:9}

We now have most of the tools necessary to prove that $\cF$ and $\cF^*$
are homeomorphisms. Since we already know that their differentials are
injective, this will basically follow from the fact that both $\cF$ and
$\cF^*$ are Fredholm of index $0$. This fact is essentially a consequence
of the 
ellipticity of the underlying PDE problems. Actually there are two ways
to consider our geometric problem in an elliptic setting. 

The most natural way is to describe Theorems \ref{tm:I} and \ref{tm:III} as
statements about the boundary value problem of finding a hyperbolic
metric on $M$ with some conditions on $\dr M$. The linearization is then
an elliptic boundary value problem on $(M, \dr M)$, and it should be
possible to prove that its index is $0$ by deforming it to a more
tractable elliptic problem. This is done, for the statement concerning the
induced metric, in \cite{ecb} (in higher dimensions, in the setting of
Einstein metrics of negative curvature). 

We will rather use another description, which seems somehow simpler to
expound, and also better suited to the study of the third fundamental
form. The point will be to remark that, given a hyperbolic manifold with
convex boundary, each boundary component has an equivariant embedding in
$H^3$. Such equivariant embeddings for the various boundary components
usually do not "fit" so as to be obtained from a hyperbolic metric on
$M$, but, since the space of complete hyperbolic metrics on $E(M)$ is
finite dimensional, the obstruction space is also finite
dimensional. Thus the index of the problem concerning hyperbolic metrics
on $M$ can be computed rather simply from the index of the problems
concerning equivariant embeddings of surfaces in $H^3$, which are
indicated by the Riemann-Roch theorem.

\pg{The Nash-Moser theorem}

There is a technical detail of some importance that should be pointed
out; namely, the basic analytical tool that we will use here is the
Nash-Moser inverse function theorem. One could think that a simpler
approach using Banach spaces should be sufficient, and should lead to
results in the $C^r$ or $H^s$ categories, rather than only in the
$C^\infty$ category which appears in the statements of Theorems
\ref{tm:I} and \ref{tm:III}. 

The reason the $C^\infty$ category, and the Nash-Moser theorem, are
necessary, is that there is a "loss of derivatives" phenomenon at work
here, which is typical of some isometric embedding questions, and which
also appears in questions on Einstein manifolds with boundary
\cite{ecb}. Namely, when one deforms a surface $\Sigma$ in $H^3$ through
a vector field $v$, the induced
variation of the induced metric has two components:
\begin{itemize}
\item the normal component of $v$ induces a deformation proportional to
  the second fundamental form of $S$, which is of order $0$ with respect
  to $v$;
\item the tangent component of $v$ induces a deformation which is of
  order $1$ with respect to $v$.
\end{itemize}
Since those two kinds of deformations have different orders in $v$, a
Banach-space approach runs into difficulties; the $C^\infty$ category is
required in our proof, as well as the Nash-Moser Theorem. It remains however
likely that an analogous statement is true in the $C^s$ or the $H^s$
categories; proving it would demand additional efforts. 

\pg{Main lemma}

The key point is the following statement. 

\blm \label{lm:indice}
$\cF$ and $\cF^*$ are local homeomorphisms.
\elm

We first consider the map $\cF$. We will decompose it as follows:
$$ \cF:\cG \stackrel{\cF_1}{\longrightarrow} \cE
\stackrel{\cF_2}{\longrightarrow} \cH\times \cR
\stackrel{\Pi_1}{\longrightarrow} \cH~, $$
where:
\begin{itemize}
\item $\dr_1M, \dr_2M, \cdots, \dr_NM$ are the connected components of
  $\dr M$. 
\item For each $i\in \{ 1,\cdots, N\}$, $\cE_i$ is the space of smooth,
  strictly convex equivariant embeddings (see the definition below)
  of $\dr_iM$ in $H^3$, considered
  modulo right composition by an isotopy of $\dr_iM$ and modulo 
  the global isometries of $H^3$, and $\cE:=\prod_{i=1}^{N}\cE_i$.  
\item $\cF_1$ is the map sending a hyperbolic metric $g\in \cG$ to the
  equivariant embeddings of the $\dr_iM$ corresponding to the boundary components
  of $\Mt$, seen as a convex subset of $H^3$.  
\item For each $i\in \{ 1,\cdots, N\}$, $\cH_i$ is the space of smooth metrics
  on $\dr_iM$ with curvature $K>-1$, considered up to isotopy, 
  and $\cH:=\prod_{i=1}^{N}\cH_i$.
\item For each $i\in \{ 1,\cdots, N\}$, $\cR_i$ is the space of irreducible
  representations of $\pi_1(\dr_iM)$ in $PSL(2,\C)$, considered up to global
  isometry (see below) and $\cR:=\prod_{i=1}^{N}\cR_i$.
\item $\cF_2$ is the natural map sending an $N$-uple of equivariant embeddings
  of the $\dr_iM$ to the $N$-uples of their induced metrics and of their
  representations. 
\item $\Pi_1:\cH\times \cR\rightarrow \cH$ and $\Pi_2:\cH\times \cR\rightarrow
  \cR$ are the projections on the first and on the second factor, respectively.
\end{itemize}
We will give some details on various elements of this decomposition before
proving Lemma \ref{lm:indice}.

\pg{Representations}

For each $i\in \{ 1, \cdots, N\}$, denote by $g_i$ the genus of $\dr_iM$. 

\bdf \label{df:R}
For each $i\in \{ 1, \cdots, N\}$, we call $\cR_i$ the space of irreducible
representations of $\pi_1(\dr_i M)$ in $\isom(H^3)$, modulo the global
isometries, i.e.:
$$ \cR_i:= \mbox{Hom}_{irr}(\pi_1(\dr_iM), \isom(H^3))/\isom(H^3)~. $$
We also define $\cR$
as their product, $\cR:=\prod_{i=1}^{N}\cR_i$. We call $\cR_M$ the
subset of $\cR$ of elements which come from an element $g\in \cG$,
i.e. a hyperbolic metric with smooth, strictly convex boundary on $M$. 
\edf

The dimension of those representation spaces can be computed simply
since the $\pi_1(\dr_iM)$ are finitely generated groups. We only outline the
proof here, and refer the reader to e.g. \cite{kapovich} for details.

\bprop \label{pr:dimR}
For each $i\in \{ 1, \cdots, N\}$, $\dim \cR_i=12g_i-12$. So $\dim \cR =
\sum_{i=1}^{N} (12g_i-12)$. $\cR_M$ is a submanifold of $\cR$ of
dimension $\sum_{i=1}^N (6g_i-6)$. 
\eprop

\bpv[Heuristics of the proof]
For each $i\in \{ 1, \cdots, N\}$, $\pi_1(\dr_iM)$ has a presentation
with $2g_i$ generators and one relation. Since $\dim \isom(H^3)=6$, the
dimension of $\mbox{Hom}(\pi_1(\dr_iM), \isom(H^3))$ is
$12g_i-6$. Quotienting by the global isometries removes $6$ dimensions. 

The elements $\rho\in \cR_M$ are in one-to-one correspondence with the
complete, convex co-compact hyperbolic metrics on $M$. By the
Ahlfors-Bers theorem (see e.g. \cite{ahlfors}) those metrics are in
one-to-one correspondence with the conformal metrics on $\dr M$. But
this space is the product of the Teichm{\"u}ller spaces of the $\dr_iM$, so
it has dimension $\sum_{i=1}^N (6g_i-6)$.
\epv

\pg{Equivariant embeddings}

We will also need the notion of equivariant embedding of a
surface in $H^3$. In particular, given a hyperbolic metric on $M$, each
connected component of $\dr M$ inherits a natural equivariant embedding
in $H^3$. 

\bdf \label{df:equiv}
Let $\Sigma$ be a compact surface. An equivariant embedding of $\Sigma$
in $H^3$ is a couple $(\phi, \rho)$, where:
\begin{itemize}
\item $\phi:\Sigmat\rightarrow H^3$ is an embedding.
\item $\rho:\pi_1\Sigma \rightarrow \isom(H^3)$ is a morphism.
\item for each $x\in\Sigmat$ and each $\gamma\in \pi_1\Sigma$:
$$ \phi(\gamma x) = \rho(\gamma)\phi(x)~. $$
\end{itemize}
$\rho$ is the {\bf representation} of the equivariant embedding $(\phi,
\rho)$. 
\edf

As already mentioned, given a hyperbolic metric $g\in \cG$ on $M$, its
developing map defines a family of equivariant embeddings of the
boundary components of $M$ into $H^3$. 

\bdf \label{df:E}
For each $i\in \{ 1, \cdots, N\}$, we call $\cE_i$ the space of smooth,
strictly convex equivariant embeddings of $\dr_i M$ in
$H^3$, considered modulo right composition with an isotopy, and modulo the global
isometries. $\cE:=\prod_{i=1}^{N}\cE_i$.
\edf

As already mentioned, given a hyperbolic metric $g\in \cG$ on $M$, its
developing map defines a family of equivariant embeddings of the
boundary components of $M$ into $H^3$.

\bdf \label{df:F1}
We call:
$$
\cF_1: \cG\rightarrow \cE
$$
the map sending a hyperbolic metric on $M$ with smooth, strictly convex
boundary, to the N-uple of convex equivariant embeddings of the boundary
components of $M$.
\edf

Given an equivariant embedding of a surface $\Sigma$ in $H^3$, it
induces a metric and a third fundamental form on the universal cover of
$\Sigma$, but also, by the equivariance which is in the definition, on
$\Sigma$ itself. For each $i\in \{ 1, \cdots, N\}$, we call $\cH_i$ the
space of smooth metrics on $\dr_iM$ with curvature $K>-1$, and $\cH^*_i$
the space of smooth metrics on $\dr_iM$ with curvature $K<1$. Then
$\Pi_{i=1}^N \cH_i=\cH$, while $\Pi_{i=1}^N \cH^*_i\subset \cH^*$. 

\bdf \label{df:F2}
For each $i\in \{ 1, \cdots, N\}$, we call $\cF_{2,i}$ the map from
$\cE_i$ to $\cH_i\times \cR_i$ sending an equivariant embedding of
$\dr_i M$ to its induced metric and its representation. Then we call:
$$
\begin{array}{rrcl}
  \cF_2: & \cE & \rightarrow & \cH\times \cR \nonumber \\
  & (\phi_i)_{1\leq i\leq N} & \mapsto & (\cF_{2,i}(\phi_i))_{1\leq i\leq
  N}~. \nonumber 
\end{array}
$$
\edf

$\cF_1$ is a map between infinite dimensional spaces, but it is Fredholm
and has a finite index. The index computation basically uses the fact
that $\cR$ is finite dimensional, and the codimension of $\cR_M$ in
$\cR$. More precisely, the proof will use the following facts.

\bprop \label{pr:facts}
\begin{enumerate}
\item $\Pi_2\circ\cF_2:\cE\rightarrow \cR$ is a submersion.
\item Therefore, $\cE_M:=(\Pi_2\circ\cF_2)^{-1}(\cR_M)$ has codimension
  $\sum_i(6g_i-6)$ in $\cE$. 
\item $\cF_1$ is a local homeomorphism between $\cG$ and $\cE_M$. 
\end{enumerate}
\eprop

\bpv
Point (1) is elementary, since each first-order deformation of the
representation of an equivariant embedding of a compact surface can be
obtained by a first-order deformation of the equivariant embedding. Point (2)
is a consequence since $\cR_M$ has codimension $\sum_i(6g_i-6)$ in $\cR$. 

Point (3) is almost as elementary. Clearly, only elements of $\cE_M$ are in
the image of $\cF_1$, this is by definition of $\cE_M$. Conversely,
given an $N$-uple of equivariant convex embeddings of the $\dr_iM$ such that
their representations come from a complete, convex co-compact hyperbolic
metric on $M$ (the 
definition of $\cR_M$), it is clear that the equivariant embeddings come from
a hyperbolic metric $g\in \cG$ on $M$. The fact that the relation is
one-to-one is a consequence of the fact that the boundary surfaces --- seen as
equivariant surfaces in $H^3$ --- uniquely determine $g$. 
\epv

The next step is that $d\cF_2$ is Fredholm, and its index can be
computed easily using the Riemann-Roch Theorem.

\bprop \label{pr:ind2}
For each $i\in \{ 1, \cdots, N\}$, $d\cF_{2,i}$ is a Fredholm map of
index $-(6g_i-6)$. Therefore, $d\cF_2$ is a Fredholm map of index
$-\sum_{i=1}^{N} (6g_i-6)$. 
\eprop

\bpv
The second part clearly follows from the first. To prove the first part,
we choose $i\in \{ 1, \cdots, N\}$.

Let $(\phi_i, \rho_i)\in \cE_i$, so that $\phi_i$ is a convex equivariant
embedding of $\dr_iM$ in $H^3$, and that $\rho_i$ is its
representation. Let $\cE_i(\rho_i)$ be the subspace of $\cE_i$ of
convex equivariant embeddings of $\dr_iM$ with the same representation
$\rho_i$. The tangent space $T_{(\phi_i,\rho_i)}\cE_i$ is naturally associated to a
space of vector fields along $\phi_i(\tilde{\dr_iM})$ with some equivariance
property under the action of $\rho_i$; the tangent space of $\cE_i(\rho_i)$,
however, is simply associated to the space of vector fields along
$\phi_i(\tilde{\dr_iM})$ which are invariant under the action of $\rho_i$,
i.e. with the space of smooth sections of the pull-back by $\phi_i$ on
$\dr_iM$ of $TH^3$. 

Let $u$ be a smooth section of this bundle, it has a natural
decomposition as $u=\lambda N+v$, where $N$ corresponds (under the
pull-back map) to the unit normal vector of $\phi_i(\tilde{\dr_iM})$,
while $v$ is tangent to $\dr_iM$. The first-order variation of the
induced metric on $\dr_iM$ induced by $\lambda N$ is simply $\lambda
\II$, while the component orthogonal to $\II$ of the first-order
variation of the induced metric associated to $v$ is given by $\db_Iv$,
where the operator $\db_I$, defined in section \ref{se:3}, has the same
principal symbol as the $\db$ operator of the second fundamental form of
$\dr_iM$. 

By the Riemann-Roch theorem (or the index theorem), the $\db$ operator
acting on vector fields on $\dr_iM$ is Fredholm and has index $-(6g_i-6)$. 
So $\db_I$ is Fredholm with the same index. But, from the way $u$ acts
on the induced metric, this means that the map sending a section of
$\phi_i^*(TH^3)$ to the induced variation of the induced metric is also
Fredholm with index $-(6g_i-6)$. This map is the differential at
$\phi_i$ of the restriction of $\cF_{2,i}$ to $\cE_i(\rho_i)$, seen as a
map from $\cE_i(\rho_i)$ to $\cH_i\times \{ \rho_i\}$. But the codimension
of $\cE_i(\rho_i)$ in $\cE_i$ is equal to the codimension of $\cH_i\times 
\{ \rho_i\}$ in $\cH_i\times \cR_i$. So the index of $\cF_{2,i}$ is still
$-(6g_i-6)$. 
\epv

We can now move to the key point. 

\bprop \label{pr:pi}
Let $g\in \cG$. In the neighborhood of $\cF_1(g)\in \cE_M$, the restriction of
$\cF_2$ to $\cE_M$, composed witht he projection $\Pi_1$, 
is a local homeomorphism onto $\cH$. 
\eprop

\bpv
The proof uses the Nash-Moser theorem, and we will use results from
\cite{H}. To keep as close as possible from this text, we first define
versions of $\cE$ and $\cH$ without taking a quotient by the group of
isotopies of $\dr M$. So, for each $i\in \{ 1, \cdots, N\}$, we call $\cEb_i$
the space of smooth, strictly convex equivariant embeddings of $\dr_iM$ into
$H^3$, up to global isometries, 
and let $\cEb:=\sum_{i=1}^N \cEb_i$. We also call $\cEb_M$ the space of
elements of $\cEb$ which project to $\cE_M$. It follows from \cite{H} (see
Corollary II.2.3.2, p. 146) that $\cEb$ is a smooth, tame Fr{\'e}chet manifold. 

In the same way, for each $i\in \{
1, \cdots, N\}$, we call $\cHb_i$ the space of smooth metrics on $\dr_iM$ with
$K>-1$, and let $\cHb:=\sum_{i=1}^N \cHb_i$. $\cHb$ is an open subset in
a smooth, tame Fr{\'e}chet space, by Theorem II. 2.3.1 of \cite{H}, p. 146.

There is a natural map $\cFb_2:\cEb\rightarrow \cHb\times \cR$. We will
change it a little. We now restrict our attention to a neighborhood of some
fixed $g\in \cG$, and of its images in $\cE$ and in $\cH\times \cR$. Since
$\cR_M$ is a submanifold of codimension $\sum_i (6g_i-6)$ of $\cR$, we can
find a smooth submersion $\Sigma:\cR\rightarrow N$, where $N$ is a vector
space of dimension $\sum_i (6g_i-6)$, such that $\cR_M=\Sigma^{-1}(0)$. Let
$\cFb'_2:\cEb\rightarrow \cH\times N$ 
be the map obtained from $\cFb_2$ by composing
it with $\Sigma$ on the $\cR$ factor. As a consequence of Proposition
\ref{pr:ind2}, we see that, at each point, $d\cFb'_2$ is Fredholm with index
$0$. Moreover, it is injective by Lemma \ref{lm:rigidity}.

Let $i\in \{ 1,\cdots, N\}$. 
Given an element of $(h,n)\in T\cH\times N$, with $h=(h_1, \cdots,
h_N)$, one can find its inverse
images by $d\cFb'_2$, which are equivariant vector fields along the
image of the $\dr_i M$, 
in a 2-steps process, repeated for all $i\in \{ 1,\cdots,N\}$:
\begin{itemize}
\item the components of the vector fields normal to the surfaces are
  given by the components of the $h_i$ along the second fundamental
  forms of the $\dr_iM$.
\item their components tangent to the surfaces are solutions of $\db_I
  v=h_i'$, where the $h'_i$ are given by the component of the $h_i$
  orthogonal to the second fundamental forms of the $\dr_iM$. 
\end{itemize}
The second step involves the resolution of an elliptic PDE. 

We can thus apply Theorem II.3.3.3 of \cite{H}, p. 158; it shows that, in the
neighborhood of $\cF_1(g)$, the differentials of $\cFb_2'$ are invertible, and
that the inverses define a smooth, tame family of linear maps. So we can
apply the Nash-Moser inverse function theorem (\cite{H}, Theorem
III.1.1.1, p. 171) and see that $\cFb_2'$ is a local homeomorphism. Restricting
our attention to $\cEb_M$, it follows that the restriction of $\cFb_2'$ to
$\cEb_M$ is a local homeomorphism onto $\cHb$. 

We can now take the quotient on both sides by the isotopy group of $\dr M$,
and obtain that the restriction of $\cF_2$ to $\cE_M$ is a local homeomorphism
onto $\cH$. 
\epv

\pg{Proof that $\cF$ is a local homeomorphism}
It is obviously a consequence of the third point of Proposition \ref{pr:facts}
and of Proposition \ref{pr:pi}.

\pg{The third fundamental form}

The proof given above extends with minor differences to the third
fundamental form. We need a slightly different definition to replace
$\cF_2$. 

\bdf \label{df:F2*}
For each $i\in \{ 1, \cdots, N\}$, we call $\cF^*_{2,i}$ the map from
$\cE_i$ to $\cH^*_i\times \cR_i$ sending a convex equivariant embedding of
$\dr_i M$ to its third fundamental form and its representation. Then:
$$
\begin{array}{rrcl}
  \cF^*_2: & \cE & \rightarrow & \cH^*\times \cR \nonumber \\
  & (\phi_i)_{1\leq i\leq N} & \mapsto & (\cF^*_{2,i}(\phi_i))_{1\leq i\leq
  N}~. \nonumber 
\end{array}
$$
\edf

The following analog of Proposition \ref{pr:ind2} then holds, and its
proof is the same as the proof of Proposition \ref{pr:ind2}, except that
one has to use the hyperbolic-de Sitter duality and to operate in
$S^3_1$. 

\bprop \label{pr:ind2*}
For each $i\in \{ 1, \cdots, N\}$, $\cF^*_{2,i}$ is a Fredholm map of
index $-(6g_i-6)$. Therefore, $\cF^*_2$ is a Fredholm map of index
$-\sum_{i=1}^{N} (6g_i-6)$. 
\eprop

Now, as for $\cF$, $\cF^*=\Pi_1\circ \cF^*_2\circ \cF_1$, and the same
argument as above shows that $\cF^*$ is a local homeomorphism between $\cG$
and $\cH^*$.

\paragraph{Main proofs}
We now have all the tools to prove Theorems \ref{tm:I} and
\ref{tm:III} when $M$ is not a solid torus.

\bpn{of Theorem \ref{tm:I}}
According to Lemma \ref{lm:indice}, $\cF$ is a local homeomorphism
between $\cG$ and $\cH$. Since it is also proper by Lemma
\ref{lm:compact}, it is a covering. But since $\cG$ is connected and
$\cH$ is connected and simply connected (by Lemma \ref{lm:cnx}), $\cF$ is
a global homeomorphism, and Theorem \ref{tm:I} follows.
\epn

\bpn{of Theorem \ref{tm:III}} 
The proof is exactly the same as that of Theorem \ref{tm:I}, using
Lemmas \ref{lm:rigidity*} and \ref{lm:compact*} instead of
\ref{lm:rigidity} and \ref{lm:compact}. 
\epn 

\paragraph{The solid torus}

We now suppose that $M$ is a solid torus. 
The proofs of Theorems \ref{tm:I} and \ref{tm:III}
are then mostly simpler than in the general
case, but a little different. One reason is that the boundary is a torus ---
while in the other cases the boundary components are always surfaces of higher
genus. Another is that the limit set has now only two points, so that the
proof of the infinitesimal rigidity is a little different (but rather
simpler). 

The proof is based on the same series of lemmas as when $M$ is not the solid
torus; we repeat them here, with sketches of proofs only when they differ from
the general case detailed above.

\blm \label{lm:tore:fredholm}
The map $\cF:\cG\rightarrow \cH$ (resp. $\cF^*:\cG\rightarrow \cH^*$) sending
a hyperbolic metric on $M$ to the induced metric (resp. third fundamental
form) of the boundary is Fredholm. 
\elm

This is proved as in Lemma \ref{lm:indice}. There are differences, however, in
the computation of the index (although it is still $0$). The key points are as
follows. 

\blm \label{lm:tore:indice}
\begin{enumerate}
\item The map sending an equivariant smoooth, strictly convex embedding
  of $T^2$ into $H^3$ to its induced metric is Fredholm, with index $2$.
\item The space of representations of $\Z^2$ in $\mbox{SO}(3,1)$ has dimension
  $4$. 
\item The space of representations of $\Z$ in $\mbox{SO}(3,1)$ has dimension
  $2$.
\end{enumerate}
\elm

\bpv
The first point is a direct consequence of the Riemann-Roch theorem, as
in the proof of Proposition \ref{pr:ind2}. The
second is a simple computation using the relationship between the images of
the two generators (they commute) while the third uses the fact that a
hyperbolic isometry without fixed point is uniquely determined, up to
conjugation, by its translation length and rotation angle. 
\epv

\bcr
$\cF$ and $\cF^*$ are Fredholm operators of index $0$.
\ecr

\bpv 
This follows from the previous lemma, because $2=4-2$. 
\epv

\blm \label{lm:tore:rig}
$\cF$ and $\cF^*$ have injective differentials at each point of $\cG$. 
\elm

\bpv
It follows the proof of Lemmas \ref{lm:rigidity} and \ref{lm:rigidity*}. The
main difference is that the limit set $\Lambda$ now has only two points, so
the analytic arguments are made rather easier. 
\epv

\blm \label{lm:tore:compact}
$\cF$ and $\cF^*$ are proper. 
\elm

The proof uses the same tools as for Lemmas \ref{lm:compact} and
\ref{lm:compact*}. 

\blm \label{lm:tore:topo}
$\cG$ is connected, and $\cH$ and $\cH^*$ are simply connected. 
\elm

\bpv
The connectedness of $\cG$ can be obtained in an elementary way, by using the
connectedness of the space of actions without fixed point of $\Z$ on $H^3$. To
show that $\cH$ and $\cH^*$ are simply connected, one can use the same trick
as in the proof of Lemma \ref{lm:cnx}. 
\epv

\medskip

The proofs of Theorems \ref{tm:I} and \ref{tm:III} follow from those lemmas,
when $M$ is a solid torus, just as in the other cases (when $M$ is not a solid
torus) above, because $\cF$ and $\cF^*$ are local covering, which then have to
be homeomorphisms. 



\section{Other conditions on $\dr M$}

We will prove in this section theorems \ref{tm:cor-I} and
\ref{tm:cor-III}. They are simple consequences of Theorems \ref{tm:I}
and \ref{tm:III}, along with some well-known properties
of equidistant surfaces in hyperbolic 3-space, and in particular of the
following statement.

\bprop \label{pr:equidist-I}
Let $S$ be a (complete) convex surface in $H^3$. For $t>0$, let
$S_t$ be the equidistant surface at distance $t$ from $S$. If $I_t$ and
$\III_t$ are the induced metric and third fundamental form,
respectively, of $S_t$, then:
$$ I_t = \cosh^2(t) I_0 + 2\cosh(t)\sinh(t) \II_0 + \sinh^2(t) \III_0~, $$
$$ \III_t = \sinh^2(t) I_0 + 2\cosh(t)\sinh(t) \II_0 + \cosh^2(t) \III_0~. $$
Moreover, the principal curvatures of $S_t$ are bounded between
$\tanh(t)$ and $1/\tanh(t)$. 
\eprop

The proof, which is a simple computation, is left to the reader. As a
consequence, we have the following analogous statement, whose proof we
also leave to the reader.

\bprop \label{pr:equidist-III}
Let $k_0\in (0, 1)$. 
Let $S$ be a (complete) convex surface in $H^3$, with principal
curvatures in $(k_0, 1/k_0)$. For each $t>0$ with $\tanh(t)\leq k_0$, the
set $S_{-t}$ of points at distance $t$ from $S$, on the normal to $S$,
in the direction of the convex side of $S$, is a smooth, convex
surface. Its induced metric and third fundamental form, respectively,
are:
$$ I_{-t} = \cosh^2(t) I_0 - 2\cosh(t)\sinh(t) \II_0 + \sinh^2(t) \III_0~, $$
$$ \III_{-t} = \sinh^2(t) I_0 - 2\cosh(t)\sinh(t) \II_0 + \cosh^2(t)
\III_0~. $$ 
\eprop

We can now prove Theorem \ref{tm:cor-I}, which we recall for the
reader's convenience. 

\begin{tmn}{\ref{tm:cor-I}}
Let $k_0\in (0,1)$. Let $h$ be a Riemannian metric on $\dr M$. There
exists a hyperbolic metric $g$ on $M$, such that the principal
curvatures of the boundary are between $k_0$ and $1/k_0$, and for which:
$$ I - 2k_0 \II + k_0^2 \III = h $$
if and only if $h$ has curvature $K> -1/(1-k_0^2)$. $g$ is then
unique. 
\end{tmn}

\bpv
Let $h$ be a smooth metric on $\dr M$, with curvature $K>
-1/(1-k_0^2)$. Let $t_0:=\tanh^{-1}(k_0)$, and set $\hb:=\cosh^2(t_0)
h$. Since $\cosh^2(t_0)=1/(1-k_0^2)$, $\hb$ has curvature $K>-1$, so by
Theorem \ref{tm:I}, there exists a unique hyperbolic metric $\gb$ with
smooth convex boundary such that the induced metric on $\dr M$ is
$\hb$. 

Now consider the extension $E(M)$ of $M$, and let $g$ be the hyperbolic
metric on $M$ obtained from $\gb$ by replacing each connected component
of the boundary of $M$ by the equidistant surface at distance $t_0$ in
the corresponding end of $E(M)$. Proposition \ref{pr:equidist-I} shows
that, for this metric, the principal curvatures of the boundary are
bounded between $k_0$ and $1/k_0$, and the extrinsic invariants $I$,
$\II$ and $\III$ of the boundary are such that:
$$ \hb = \cosh^2(t_0) I - 2\cosh(t_0)\sinh(t_0) \II
+\sinh^2(t_0)\III~. $$
By definition of the scaling used to define $\hb$ from $h$, we see that:
$$ h = I - 2\tanh(t_0) \II + \tanh^2(t_0)\III = I -2k_0 \II+ k_0^2\III~,
$$
as needed. 

Conversely, let $g'$ be a hyperbolic metric with convex boundary such
that the principal curvatures of the boundary are bounded between $k_0$
and $1/k_0$, and such that, on the boundary, the extrinsic invariants
$I', \II'$ and $\III'$ are such that
$I'-2k_0\II'+k_0^2\III'=h$. Apply proposition \ref{pr:equidist-III}, with
$t_0:=\tanh^{-1}(k_0)$. We obtain a hyperbolic metric $\gb'$ with
smooth, locally convex boundary, such that the induced metric on the
boundary is:
\begin{eqnarray}
I & = & \cosh^2(t_0)^2 I' - 2\cosh(t_0)\sinh(t_0)\II' +
\sinh^2(t_0)\III' \nonumber \\
& = & \cosh^2(t_0) h~. \nonumber
\end{eqnarray}
So $\gb'$ corresponds to the metric $\gb$ constructed above, and
therefore $g'=g$, and this proves the uniqueness part of the statement. 
\epv

The proof of Theorem \ref{tm:cor-III} is analogous, using Theorem
\ref{tm:III} instead of Theorem \ref{tm:I}.

\begin{tmn}{\ref{tm:cor-III}}
Let $k_0\in (0,1)$. Let $h$ be a Riemannian metric on $\dr M$. There
exists a hyperbolic metric $g$ on $M$, such that the principal
curvatures of the boundary are between $k_0$ and $1/k_0$, and for which:
$$ k_0^2 I - 2k_0 \II + \III = h $$
if and only if $h$ has curvature $K< 1/(1-k_0^2)$, and its closed,
contractible 
geodesics have length $L>2\pi\sqrt{1-k_0^2}$. $g$ is then
unique. 
\end{tmn}

\bpv
Let $h$ be a smooth metric on $\dr M$, with curvature $K<
1/(1-k_0^2)$ and with closed geodesics of length $L>2\pi\sqrt{1-k_0^2}$. 
Let $t_0:=\tanh^{-1}(k_0)$, and set $\hb:=\cosh^2(t_0)
h$. Since $\cosh^2(t_0)=1/(1-k_0^2)$, $\hb$ has curvature $K<1$, and its
closed geodesics have length $L>2\pi$. So by
Theorem \ref{tm:III}, there exists a unique hyperbolic metric $\gb$ with
smooth convex boundary such that the third fundamental of $\dr M$ is
$\hb$. 

Now consider the extension $E(M)$ of $M$, and let $g$ be the hyperbolic
metric on $M$ obtained from $\gb$ by replacing each connected component
of the boundary of $M$ by the equidistant surface at distance $t_0$ in
the corresponding end of $E(M)$. Proposition \ref{pr:equidist-I} shows
that, for this metric, the principal curvatures of the boundary are
bounded between $k_0$ and $1/k_0$, and the third fundamental for $\gb'$
of the boundary is:
$$ \hb = \sinh^2(t_0) I - 2\cosh(t_0)\sinh(t_0) \II
+\cosh^2(t_0)\III~. $$
Using the scaling used to define $\hb$ from $h$, we see that:
$$ h = \tanh^2(t_0)I - 2\tanh(t_0) \II + \III = k_0^2I -2k_0 \II+ \III~,
$$
as needed. 

Conversely, let $g'$ be a hyperbolic metric with convex boundary such
that the principal curvatures of the boundary are bounded between $k_0$
and $1/k_0$, and such that, on the boundary, the extrinsic invariants
$I', \II'$ and $\III'$ are such that
$k_0^2I'-2k_0\II'+\III'=h$. Apply proposition \ref{pr:equidist-III}, with
$t_0:=\tanh^{-1}(k_0)$. We obtain a hyperbolic metric $\gb'$ with
smooth, locally convex boundary, such that the third fundamental for of
the boundary is:
\begin{eqnarray}
I & = & \sinh^2(t_0)^2 I' - 2\cosh(t_0)\sinh(t_0)\II' +
\cosh^2(t_0)\III' \nonumber \\
& = & \cosh^2(t_0) h~. \nonumber
\end{eqnarray}
So $\gb'$ corresponds to the metric $\gb$ constructed above, and
therefore $g'=g$, and this proves the uniqueness part of the statement. 
\epv

\pg{Some open questions}

The smoothness assumptions in Theorems \ref{tm:I} and \ref{tm:III} are not
quite satisfactory. It should be possible to state results also in the $C^k$
category. It is not done here because of a well-known "loss of derivatives"
phenomenon, already apparent for isometric embeddings of surfaces. 

As already mentioned in the introduction, the phenomena behind Theorems
\ref{tm:I} and \ref{tm:III} should not be limited to manifolds with
smooth boundaries, but should rather apply to general convex boundaries
with no smoothness assumption, for instance to polyhedral boundaries. 
This is indicated by the case when $M$ is
simply a ball, where general results without smoothness were obtained by
Pogorelov \cite{Po} 
for the induced metric, and by G. Moussong \cite{moussong} for the third
fundamental form. Moreover, in the case of the third fundamental form,
results do hold when the boundary looks locally like an ideal
\cite{ideal} or hyperideal \cite{hphm} polyhedron. 

\medskip

It would in particular be interesting to know whether analogs of
Theorems \ref{tm:I} and \ref{tm:III} hold for the convex cores of
hyperbolic manifolds. The induced metric on the boundary is then
hyperbolic, while the third fundamental form is replaced by the bending
lamination. It is still not known whether an infinitesimal rigidity
result holds for either the induced metric or the bending lamination,
although it is known that those two questions are equivalent
\cite{bonahon-toulouse}. Moreover, it is known that any hyperbolic
metric can be achieved as the induced metric on the boundary of the
convex core of a hyperbolic 3-manifold \cite{L4}, and a characterization
of the possible bending laminations (without uniqueness) was also
obtained recently \cite{bonahon-otal,lecuire}.

\medskip

Some questions also remain concerning smooth surfaces, in
particular those which are not compact but have ``ends''. This is
formally analogous to the hyperideal polyhedra, where one has to
use the condition that, near each end, all faces  are orthogonal to a
hyperbolic plane. In the smooth case, however, it looks like the right
condition might be that the boundary at infinity of each end of infinite
area is a circle
(in the M{\"o}bius geometry of the boundary at infinity). When
one considers constant curvature metrics on the boundary, analogs of
Theorems \ref{tm:I} and \ref{tm:III} hold in
the topologically simplest case (see \cite{rsc}) and, in a wider
setting, for existence results concerning manifolds with a given induced
metric on the boundary \cite{bkc}.

\medskip

It is probably not completely necessary to restrict one's attention to
convex co-compact manifolds in Theorems \ref{tm:I} and \ref{tm:III}. I
would guess that similar results also hold for geometrically finite
manifolds. An interesting and simple 
case would be a manifold with one cusp, and one boundary component which
is topologically a torus.

\pg{Higher dimensions ?}

In my opinion, it would be quite interesting to understand whether the
results described here for hyperbolic manifolds with convex boundary
have some sort of extension in the setting of Einstein manifolds with
boundary, maybe in dimension $4$.
A first step is taken in \cite{ecb}, but it is far from
satisfactory. The main problem is the lack of a strong enough
infinitesimal 
rigidity statements for Einstein manifolds with boundary (note that one
such statement, similar to the one in \cite{ecb}, is obtained in
\cite{sem-era} by very different methods).


\section*{Acknowledgments}

This paper owes much to conversations with Fran{\c{}}{c}ois Labourie, and to
his many illuminating comments and remarks. His help was also crucial in
improving the somewhat cryptic exposition of the previous version of
this text. 

I would also thank Francis Bonahon, Julien Duval, Philippe Eyssidieux,
Vincent Guedj, Gilbert Levitt, Gregory 
McShane, Igor Rivin, Mathias Rousset and Fabrice Planchon for helpful
remarks. Theorems \ref{tm:cor-I} and 
\ref{tm:cor-III} were partly motivated by a question asked by Paul
Yang. The paper was also improved thanks to remarks from an anonymous
referee. 


\bibliographystyle{alpha}

\end{document}